\documentclass[draft]{amsart}
\usepackage{amsfonts,amsmath,amscd,amssymb,fullpage}
\usepackage{amssymb, amsbsy, amsthm, amsmath, amstext, amsopn, verbatim,
multicol}

\numberwithin{equation}{subsection}
\newtheorem*{thm}{Theorem}
\newtheorem*{prop}{Proposition}
\newtheorem*{lem}{Lemma}

\newtheorem*{cor}{Corollary}

\newtheorem*{hypothesis}{Hypothesis}

\theoremstyle{definition}
 \newtheorem*{remark}{Remark}
\newtheorem*{remarks}{Remarks}

\newcommand{\reg}[1]{#1^{\text{reg}}}
\newcommand{\PP}{\mathcal{P}}
\newcommand{\BB}{\mathcal{B}}

\newcommand{\hec}[1]{\mathcal{H}_{#1}}
\newcommand{\ep}{\epsilon}
\newcommand{\caln}{\mathcal{N}}
\newcommand{\calj}{\mathcal{J}}

\newcommand{\TT}{\mathbb{T}}

\newcommand{\h}{\mathfrak{h}}

\newcommand{\LL}{\mathcal{L}}
\newcommand{\CC}{\mathbb{C}}
\newcommand{\OO}{\mathcal{O}}

\newcommand{\NN}{\mathbb{N}}
\newcommand{\Z}{\mathbb{Z}}
\newcommand{\Q}{\mathbb{Q}}
\newcommand{\R}{\mathbb{R}}
\newcommand{\C}{\mathbb{C}}

 \newcommand{\ttt}{\textsf}
\newcommand{\KZ}{\ttt{KZ}}
\newcommand{\irr}[1]{\ttt{Irrep}(#1)}
\newcommand{\hr}{\mathfrak{h}^{\text{reg}}}
\newcommand{\cxy}{\C [\h\oplus \h^*]}
\newcommand{\cx}{\C [\h]}

\newcommand{\hh}{\mathbf{h}}
\newcommand{\EE}{\mathbf{E}}

\newcommand{\WW}{{W}}    
\newcommand{\UU}{U}
\newcommand{\XXX}{\mathbb{X}}
\newcommand{\AAA}{\mathbb{A}}
\newcommand{\JJJ}{\mathbb{J}}
\newcommand{\RRR}{\mathbb{R}}
\newcommand{\dd}{\mathbb{S}}

\newcommand{\bigdot}{ \,{}_{^{^{\bullet}}}}

 \newcommand{\lgr}{\ttt{\text{-}grmod}}
 \newcommand{\lGr}{\ttt{\text{-}Grmod}}
 \newcommand{\lqgr}{\text{-}\ttt{qgr}}  
 \newcommand{\lQgr}{\text{-}\ttt{Qgr}}  
 \newcommand{\ltors}{\text{-}\ttt{tors}}
 \newcommand{\lTors}{\text{-}\ttt{Tors}}
 \newcommand{\lmod}{\text{-}\ttt{mod}}

\DeclareMathOperator{\triv}{\ttt{triv}}
\DeclareMathOperator{\sign}{\ttt{sign}}
\DeclareMathOperator{\Edeg}{\EE\text{-deg}}
\DeclareMathOperator{\hdeg}{\hh\text{-deg}}

\DeclareMathOperator{\hin}{Hilb^n\CC^2}
\DeclareMathOperator{\hi}{Hilb(n)}

\DeclareMathOperator{\tgr}{\ttt{tgr}}

\DeclareMathOperator{\ogr}{\ttt{ogr}}
\DeclareMathOperator{\ord}{\ttt{ord}}

 \DeclareMathOperator{\ten}{\ttt{ten}}
\DeclareMathOperator{\Hom}{Hom}

\DeclareMathOperator{\prj}{Proj}
\DeclareMathOperator{\spc}{Spec}

\DeclareMathOperator{\coh}{\ttt{Coh} }

\DeclareMathOperator{\ad}{ad}
\DeclareMathOperator{\gr}{gr}

\DeclareMathOperator{\qgr}{\ttt{qgr}}  
\DeclareMathOperator{\QGr}{\ttt{Qgr}}  

\begin{document}

\title{Rational Cherednik algebras and Hilbert schemes}
   \author{I. Gordon} \address{Department of Mathematics,
 Glasgow  University, Glasgow G12 8QW, Scotland}
\email{ig@maths.gla.ac.uk}
 \author{J. T. Stafford}
\address{Department of Mathematics, University of Michigan, Ann Arbor,
MI 48109-1109, USA.}
\email{jts@umich.edu}
   \thanks{ The first author was supported by the Nuffield Foundation
    Grant NAL/00625/A and the Leverhulme trust. He would like to thank
    the University of Washington and the University of California at
    Santa Barbara for their hospitality while parts of this paper
    were written. The second author was supported in part by the
NSF through the  grants DMS-9801148 and DMS-0245320.
Part of this research was conducted while he was visiting the
Mittag-Leffler Institute and he would like to thank them for their hospitality
and financial support.}

\dedicatory{Dedicated to Mike Artin
on the occasion of his $70^{\mathrm{th}}$ birthday}

  \begin{abstract} Let $H_c$ be the rational Cherednik algebra of type
  $A_{n-1}$ with spherical subalgebra $U_c=eH_ce$. Then $U_c$ is filtered
  by order of differential operators, with associated graded ring
  $\gr U_c=\C[\h\oplus\h^*]^\WW$ where $\WW$ is the $n$-th symmetric group.
 We construct
  a filtered $\Z$-algebra $B$ such that, under mild conditions on $c$:
\begin{itemize}
\item{} the category $B\lqgr$
  of graded noetherian  $B$-modules modulo torsion is equivalent to $U_c\lmod$;
\item{}  the associated graded
 $\Z$-algebra $\gr B$ has $\gr B\lqgr \simeq \coh \hi$, the category of
 coherent sheaves on the  Hilbert scheme of points in the plane.
  \end{itemize}
This can be regarded as saying that $U_c$ simultaneously gives a noncommutative
deformation of   $\h\oplus\h^*/\WW$ and of its resolution of singularities
$\hi\to \h\oplus\h^*/\WW$. As the companion
paper \cite{GS2} shows, this result is a powerful tool for
studying the representation theory of
$H_c$ and its relationship to $\hi$.
  \end{abstract}

   \maketitle
 \tableofcontents

\clearpage

\section{Introduction}

\subsection{}\label{sec101}
This is the first of two closely related papers on rational
Cherednik algebras.

In their short history, Cherednik algebras have been influential
in a surprising range of subjects: for example they have been used
to answer questions in integrable systems, combinatorics,   and
symplectic quotient singularities (see  \cite{BEGqi, gordc, BFG,
GK}).  In this paper we strengthen the connections between
Cherednik algebras   and geometry by showing that they can be
regarded as noncommutative deformations of  Hilbert schemes of
points in the plane. In the sequel \cite{GS2} this will be used to
show the close relationship between  modules over the Cherednik
algebra and sheaves on the Hilbert scheme as well as to answer
various open problems about these modules.

\subsection{}\label{intro-1.2} Fix $c\in \C$. We assume throughout
the paper that $c\notin \frac{1}{2} + \Z$ and, for simplicity, we will
also assume that  $c\not\in \mathbb{R}_{\leq 0}$ in this introduction,
 see \eqref{morrat-cor} and \eqref{morrat-cor-remark} for the more
 general case.

 Let $H_c= H_{1,c}$ be the
 rational Cherednik algebra of type $A_{n-1}$ with spherical subalgebra
  $U_c = eH_ce$. The formal definition of $H_c$ is given
in \eqref{hc-defn} but one may regard it as a deformation of the twisted group
ring $D(\h)\ast \WW$, where $D(\h)$ is the ring of differential operators on
$\h\cong \C^{n-1}$ with  the natural action of the symmetric group $\WW =
\mathfrak{S}_n$. The algebra  $U_c$ is then the corresponding
deformation of the fixed ring $D(\h)^\WW$.
The algebras $U_c$ and $H_c$ have a natural filtration by
order of differential operators with associated graded rings
$\gr{U_c} \cong \C[\h\oplus \h^\ast]^{\WW}$
and $\gr {H_c} \cong \C[\h\oplus \h^\ast]\ast {\WW}$.
Thus we may also regard $U_c$ as a
 deformation of $\C[\h\oplus \h^\ast]^{\WW}$.
In this introduction we will mostly be concerned with $U_c$,
but since $U_c$ and $H_c$ are Morita equivalent
(Corollary~\ref{morrat-cor}) the results
we prove for $U_c$ also apply to $H_c$.

It is well-known that $\h\oplus\h^*/\WW$ has a crepant resolution  $\hi\to
\h\oplus\h^*/\WW$, where $\hi$ is a variant on the  Hilbert scheme of $n$
points in the plane (see \eqref{hi-defn-sec} for  the formal definition). The
ring $U_c$ has  finite global homological dimension (see
Corollary~\ref{gldim}) and so one should expect that it has the properties of a
smooth  deformation of $\cxy^\WW$; in other words its properties should be more
closely related to those of   $\hi $ than to $\h\oplus\h^*/\WW$. Hints of this
have been reported in \cite{gordc}  and \cite{BEGfd}: finite dimensional
$H_c$-modules deform the sections of some  remarkable sheaves on   $\hi$. The
main aim of this paper is to formalise this idea by  showing that there is a
second way of passing to associated graded objects that maps $U_c\lmod$
precisely to $\coh (\hi)$.

\subsection{}\label{intro-1.3} We take our cue from the
theory of semisimple Lie algebras.
 When $n=2$,  $U_c$ is isomorphic to a factor
of $U(\mathfrak{sl}_2)$
 \cite[Section~8]{EG} and, for all $n$, the properties of
$U_c$ are similar to those of $U(\mathfrak{g})/P$, where
$P$ is  a minimal primitive ideal in the
 enveloping algebra of a complex semisimple Lie algebra $\mathfrak{g}$
  (see, for example \cite{ginz,GGOR,guay}).
 The intuition from the last paragraph  not only holds for
enveloping algebras  but can also be formalised through the
Beilinson-Bernstein
 equivalences of categories. This   gives a diagram
$$\begin{CD}
D_{\mathcal{B}} @< \sim << U(\mathfrak{g})/P \\
@V \gr VV @VV \gr V \\
\OO_{T^*\mathcal{B}} @<\tau << \OO(\mathcal{N})
\end{CD}$$
where $\mathcal{B}=G/B$ is the flag variety,
 the primitive ideal $P$ has trivial central character and
 $\tau: T^*\mathcal{B}\to \mathcal{N}$ is the Springer resolution of the
 nullcone $\mathcal{N}$.
The Morita equivalence from the sheaf of differential operators
$D_{\mathcal{B}}$ to  $  U(\mathfrak{g})/P$ is obtained by taking global
sections under the
identification $U(\mathfrak{g})/P\cong D(\mathcal{B})$.

Ginzburg has raised the question of whether a similar phenomenon
holds for Cherednik algebras (see \cite[Conjecture~1.6]{GK} for a variant on
this conjecture).  In other words,
 can one complete the following diagram?
$$\begin{CD} ?  @< \sim << U_c \\
@V \gr VV @VV \gr V \\
\OO_{\hi}  @< \tau<<\OO(\h\oplus \h^*/\WW)
\end{CD}$$

The main result of the paper gives a positive answer to this question.
Given a graded ring $R$, we write $R\lqgr $ for the quotient category of noetherian
graded $R$-modules modulo those of finite length.

\subsection{Main Theorem}\label{mainthm-intro} {\it There exists a graded ring
$B$, filtered by order of differential operators, such that
\begin{enumerate}
\item
there is an equivalence of categories \,
  $U_c \lmod\simeq B\lqgr$;
\item  there is an equivalence of categories \,
$\gr B\lqgr \simeq \coh(\hi)$.
   \end{enumerate}
}

\subsection{}\label{intro-1.4}
 The construction
of $B$ needs some explanation. For $n>2$, it can be shown that the Hilbert
scheme $\hi$  is not a cotangent bundle, so we cannot use   sheaves of
differential operators as a non-commutative model.
Instead we take as our starting point Haiman's
description of $\hi$ as a blow-up of $\h\oplus\h^*/\WW$ and
 deform this to a non-commutative setting. Set
 $A^0 = \OO(\h\oplus\h^*/\WW)$ with ideal
$I=A^1\delta$, where $\delta $ is the discriminant and $A^1 = \C[\h\oplus
\h^*]^\epsilon$ the module of anti-invariants. Then
  \cite[Proposition~2.6]{haidis} proves that
 $\hi=\ttt{Proj}\, A$ where $A=A^0[tI]$ is the Rees ring of $I$ (see
 Section~\ref{sect-haiman} for the details).

Unfortunately one cannot construct $B$ as an analogous Rees ring
over $U_c$, since  $U_c$  is a  simple ring for generic values of $c$.
We circumvent this problem
by using {\it $\Z$-algebras} (see  Section~\ref{zalg}). Specifically,
the ring $B$ from Theorem~\ref{mainthm-intro} is an algebra
 $B=\bigoplus_{i\geq j\geq 0}B_{ij}$ whose multiplication is  defined in
 matrix fashion: $B_{ij}B_{jk}\subseteq B_{ik}$  but $B_{ij}B_{\ell k}=0$ when
 $j\not= \ell$. The diagonal terms are just $B_{ii}=U_{c+i}$ while the
 off-diagonal terms $B_{ij}$   are given as the appropriate tensor products
 of the $(U_{d+1}, U_{d})$-bimodules $Q_{d}^{d+1} = eH_{d+1}\delta e$.
 The    shift functors $S_d : U_d \lmod
\rightarrow U_{d+1}\lmod$ given by tensoring with  $Q_{d}^{d+1}$
are important operators in the theory of Cherednik
 algebras and have already played a crucial role in  combinatorics
and representation theory; see, for example, \cite{BEGfd,BEGqi, gordc}. A good
way to think of the
functor $S_d$ is as the analogue of the translation functor
\cite{BG} from Lie theory.

In order to have control over  $B$
we need to know that the $Q_{d}^{d+1}$ are progenerators
for all $d\in c+\NN$; equivalently that the $S_d$ are Morita equivalences.
This is a  conjecture from
\cite[Remark~5.17]{GGOR} which we answer with:

 \subsection{Theorem}\label{morrat-intro}
  [Corollary~\ref{morrat-cor}]  {\it The shift functor $S_d$ is a Morita
  equivalence for all $d\in c+\NN$. }

 \medskip
 The significance of this result is that $B$ now  has rather pleasant
 properties; in  particular  Theorem~\ref{mainthm-intro}(1) is an easy
 consequence. For the second assertion of  Theorem~\ref{mainthm-intro},  we
 note that it is easy to obtain a $\Z$-algebra $\widehat{A}  =\bigoplus_{i\geq
 j\geq 0} A_{ij}$ from the  graded algebra $A=\bigoplus_{k\geq 0}I^k$ for which
 $A\lqgr\simeq  \widehat{A}\lqgr$. One simply takes $A_{ij}=I^{i-j}$ for each
 $i,j$. Thus the main step in the proof of Theorem~\ref{mainthm-intro} is
 given by:

  \subsection{Proposition}\label{pre-cohh-intro}
   [Theorem~\ref{main}] {\it Under the filtration
  induced from the order filtration
  of differential operators, $\gr B_{ij}  \cong A_{i-j} = I^{i-j}$
  and so $\gr B\cong \widehat{A}$ as $\Z$-algebras.
  }

\medskip
In this result the inclusion $ I^{i-j}\subseteq \gr B_{ij}$
is straightforward. The opposite inclusion is
much more subtle as it  is difficult to keep close control of
the filtration on $B_{ij}$. Our proof leans heavily on the work of Haiman
in \cite{hai3} and \cite{hai1} surrounding the $n!$ and polygraph
 theorems and
 the strategy is outlined in
more detail in~за\eqref{surjstrat}.

\subsection{Applications}\label{intro-1.10}
Theorem~\ref{mainthm-intro} gives a powerful technique for relating
$H_c$-  or $U_c$-modules to sheaves on $\hi$:
given   a $U_c$-module $M$ with a good filtration $\Lambda$
we obtain a filtered  object $(\widetilde{M},\Lambda)\in  B\lqgr $
by tensoring with $B$ and then  a
coherent sheaf $\Phi_\Lambda (M)\in \coh(\hi)$
by taking the associated graded module.

This process is studied in \cite{GS2}  where we show there  that the subtle
combinatorics  and geometry of $\hi$ is  reflected in the  representation
theory of $U_c$ and $H_c$. Let $\Delta_c(\mu)$  be the {\it standard
$H_c$-module} corresponding  to $\mu\in\irr{\WW}$ (this is the analogue of a
Verma module) with unique simple factor $L_c(\mu)$. These modules have a
natural good filtration $\Lambda$ and we mention  a couple of illustrative results
from \cite{GS2}.

\begin{itemize}
\item Suppose that $c=1/n+k$ for  $k\in \NN$, so that
$L_c(\triv)$ is the unique finite dimensional simple $H_c$-module.
Then $\Phi_\Lambda (eL_c(\triv)) \cong
 \OO_{Z_n}\otimes\mathcal{L}^{k}$, where $Z_n=\tau^{-1}(0)$
  is the {\it punctual Hilbert scheme} and $\mathcal{L}=\OO_{\hi}(1)$
  is the Serre twisting sheaf .
\item For any $c$, the  characteristic cycle of $\Phi_\Lambda
(e\Delta_c(\mu))$ equals
$\sum_{\lambda}K_{\mu\lambda}[Z_\lambda]$, where $K_{\mu\lambda}$
are Kostka numbers  and the $Z_{\lambda }$ are particular
irreducible components of $\tau^{-1}(\h\oplus \{{\bf 0} \}/\WW)$.
 \end{itemize}
The first of these results is used in \cite{GS2} to show that the natural
bigraded structure on
$\gr_\Lambda  (eL_{1/n+k})$ coincides with that on $\mathrm{H}^0(Z_n, \LL^k)$,
thus confirming a conjecture of \cite{BEGfd}. The second of these results
illustrates the subtlety of $\Phi$: if one passes directly from
$U_c$ to $\gr U_c\cong \C[\h\oplus\h^*]^\WW$ then
  $\gr_\Lambda (e\Delta_c(\mu))\cong \C[\h]\otimes \mu$
 for any $\mu$ and $c$. Thus  the support variety of
 $\gr_\Lambda e\Delta_c(\mu)$ is independent of $\mu$.

We prove one such correspondence in this paper.
Let $\PP$ denote the Procesi bundle on $\hi$, the  vector bundle of
rank $n!$ coming from Haiman's $n!$ theorem, see \eqref{PP-defn}.
 Then Corollary~\ref{cohh-subsect} proves:

\subsection{Corollary}\label{cohh-intro} If $eH_c$ is given the order filtration
$\Lambda$, then
$\Phi_\Lambda (eH_c)=\PP$.

\subsection{}\label{sect-1.10}
One reason why  Theorem~\ref{mainthm-intro} provides a
strong bridge between   Hilbert schemes and Cherednik
algebras is that the  construction of $B$  carries within it
key elements of both theories. For instance, we have already mentioned that
the shift functor $S_c$ is an analogue of the translation functor from Lie theory.
It is also the analogue of the shift functor in $\coh(\hi)$ given
by  tensoring with  $\OO_{\hi}(1)$. Indeed,
 given a   $U_c$-module $M$ with a good filtration $\Lambda$,
 it is easy to show that
$\Phi_{\Gamma}(Q^{c+1}_{c} \otimes M)= \OO_{\hi}(1)\otimes\Phi_\Lambda(M)$,
for the appropriate filtration $\Gamma$ (see \cite{GS2}).

Similarly,  Corollary~\ref{cohh-intro}  can be interpreted as saying
 that $H_c$ is a  noncommutative
analogue of the isospectral scheme $X_n$, as defined in
 \eqref{hi-basic-lem2} (see \eqref{cohh-subsect-chat} for further
 details).

\subsection{}\label{intro-1.18} The $\Z$-algebra
has the virtue that it exists whenever one has an analogue of the
translation principle; that is, one has algebras $R_i$ and progenerative
$(R_{i+1},R_{i})$-bimodules $Q_{i,i+1}$ (these
algebras can also be indexed by more general lattices than $\Z$). One can then
construct  a $\Z$-algebra as we have done and Theorem~\ref{mainthm-intro}(1)
will still hold. It is not clear when Theorem~\ref{mainthm-intro}(2) will hold and,
even when it is true, it will undoubtedly be rather subtle.

Hilbert schemes realise crepant resolutions for the symplectic
quotient singularity $(\C^2)^n/G$ whenever $G$ is the wreath
product of a finite subgroup of $SL_2(\C)$ with the symmetric
group $\WW$, see \cite[Theorem~4.2]{wang}.  We believe that our
methods will generalise to the symplectic reflection algebras
$H_{c}(G)=H_{1,c}(G)$ associated with $((\C^2)^n, G)$ to give
non-commutative deformations of those Hilbert schemes. Even when
there is no crepant resolution of such a singularity (by \cite{GK}
this happens for Weyl groups $G$  of types other than A and B) the
$\Z$-algebra associated to $H_c(G)$ will still contain interesting
information, as \cite{gordc} demonstrates.   For a
 Weyl group, the analogue of Theorem~\ref{morrat-intro} is at least
known for sufficiently large values of the defining parameter $c$
\cite[Proposition~4.3]{BEGfd}, but little is known for small
values of $c$.

The translation principle obviously holds for factors of enveloping algebras  of
semisimple Lie algebras and we can prove an analogue  of
Theorem~\ref{mainthm-intro} in this case.  However, the proof uses nontrivial
Lie theoretic results, notably the  Beilinson-Bernstein  equivalence of
categories, and it is unclear whether this approach carries information that
cannot be obtained from that equivalence.
It would be interesting to see if the recent work  \cite{BK,tan} on the
 Beilinson-Bernstein  equivalence for quantised enveloping
algebras can be understood in  a $\Z$-algebra framework.

\subsection{}\label{intro-1.19}
The paper is organised as follows. In Section~\ref{sect-rationalchered}
we recall the needed facts  about rational Cherednik algebras,
while in Section~\ref{shift} we prove Theorem~\ref{morrat-intro}. In
Section~\ref{sect-haiman} we describe
 some of Haiman's work on Hilbert schemes, adapted to the
 variety $\hi$, and use it to describe various  Poincar\'e series
 that will be needed in the proof of Theorem~\ref{mainthm-intro}.
 Section~\ref{zalg} proves the results about   $\Z$-algebras
  that were mentioned earlier in this introduction.
  Section~\ref{sect-filt} is the heart
  of the paper: in it
 we prove  Theorem~\ref{mainthm-intro}(2). This is derived from an analogous
 result about the
 associated graded module of $B_{k0}\otimes_{U_c}eH_c$ that also
 implies Corollary~\ref{cohh-intro}. Section~\ref{sect7} then gives a
 reinterpretation of Theorem~\ref{mainthm-intro}
 in terms of a tensor product filtration of $B_{ij}$.
  In  Appendix~\ref{app-a} we prove the following result that may be of
   independent interest:
{\it   Suppose that $R=\bigoplus_{i\geq 0}R_i$ is an $\NN$-graded
 algebra over a field $k$, with $R_0=k$.
   If $P$ is a right $R$-module that is both graded and projective, then $P$
   is graded-free in the sense that $P$ has a free basis of homogeneous
   elements.} This is a graded analogue of a classic result from \cite{Kap} for
   which we do not know a reference.

\subsection{Acknowledgement} We would like to thank Victor Ginzburg for bringing
his conjecture  to our attention, since it really formed the starting
point for this work.
 We would also  like to thank Tom Nevins and Catharina Stroppel for
   suggesting many improvements to us.

\section{Rational Cherednik algebras}\label{sect-rationalchered}

\subsection{}\label{rcadef}
In this section we define the rational Cherednik algebras (which will
 always be of type $A$ in this paper) and give some of the basic
  properties that will be
needed in the body of the paper.

  Let $\WW=\mathfrak{S}_n$
\label{symmetric-defn} be the {\it symmetric
group} on $n$ letters, regarded as the Weyl group of type $A_{n-1}$ acting
on its $(n-1)$-dimensional representation $\h\subset \C^n$ by permutations.
Let $\mathcal{S}=\{s=(i,j)\ \text{with}\ i<j\}\subset \WW$
\label{involution-defn} denote the reflections,  with reflecting hyperplanes
  $\alpha_s=0$. We make similar
definitions for $\h^*$ and normalise $\alpha_s^{\vee} \in \h$ so
that $\alpha_s(\alpha_s^{\vee}) = 2$.

Given $c\in \C$,  {\it  the rational Cherednik algebra of type $A_{n-1}$}
is the $\C$-algebra $H_c$\label{hc-defn}
generated by the vector
spaces $\h$ and $\h^*$ and the group $\WW$ with defining relations
\begin{eqnarray*}
wxw^{-1} = w(x), \quad wyw^{-1} = w(y), &
&\text{for all }y\in \h, x\in \h^*, w\in \WW \\
x_1x_2 = x_2x_1, \quad y_1y_2 = y_2y_1, &
&\text{for all }y_i\in \h, x_j\in \h^*\\
yx - xy = x(y) - \sum_{s\in \mathcal{S}} c\alpha_s (y) x(\alpha_s^{\vee})s, &
&\text{for all }y\in \h, x\in \h^*.
\end{eqnarray*}

We should note that the definition of the Cherednik algebra is not uniform
throughout the literature. The definition we are using agrees with that in
\cite{BEGqi,BEGfd,EG,guay} but {\it not} that from \cite{GGOR} where our
constant $c$ equals $-k_1$ for their constant $k_1$ (see
\cite[Remark~3.1]{GGOR}).

\subsection{}\label{subsec-3.2} We write the coordinate ring
of an affine variety  $V$ as $\C[V]$. By \cite[Theorem~1.3]{EG},
the subalgebra of $ H_c$ generated by $  \h^*$ can and will be identified with
  $\C[\h]$, while    $  \h$ generates a
copy of  $\C[\h^*]$ inside $H_c$ and the elements $w\in \WW$
span a copy of the group algebra $\C \WW$
in  $H_c$.
Fix once and for all dual bases $\{ x_i \}$ and $\{ y_i \}$
 of $\h^*$ and $\h$ respectively; thus
 $\C[\h]=\C[x_1,\dots,x_{n-1}]$ and
 $\C[\h^*]=\C[y_1,\dots,y_{n-1}]$.

By \cite[Theorem~1.3]{EG} there is a
 Poincar\'{e}-Birkhoff-Witt isomorphism of $\C$-vector spaces
\begin{equation}
\label{PBW}
\C[\h]\otimes_{\C} \C \WW \otimes_{\C} \C[\h^*] \xrightarrow{\sim} \ H_c.
\end{equation}
Filter $H_c$ by\label{filt-defn}
$\ord^0 H_c = \C [\h]\ast W$, $\ord^1H_c = \h + \ord^0H_c$ and
 $\ord^iH_c = (\ord^1H_c)^i$ for $i>1$, and define the
  {\it associated graded ring} to be
 $\ogr H_c=\bigoplus \ogr^n H_c$,
 where $\ogr^n H_c = \ord^n H_c/\ord^{n-1}H_c$. Then \eqref{PBW} is equivalent
 to the assertion that $\ogr H_c$
 is  isomorphic to the skew group ring $\C[\h\oplus \h^*]\ast \WW$
defined by $\sigma f = \sigma(f) \sigma$, for $\sigma\in \WW$ and
$f\in \C[\h\oplus \h^*]$.

\subsection{The Dunkl-Cherednik representation}\label{dunch}
Let $\delta \in \C[\h]$\label{delta-defn} denote the discriminant
polynomial
$\delta = \prod_{s\in \mathcal{S}} \alpha_s$.
Thus $\delta$ transforms under $\WW$ by the sign representation and
  $\hr=\h \setminus \{\delta=0\}$ is
the subset of $\h$ on which the action of $\WW$ is free. By
\cite[Proposition~4.5]{EG} there
is an injective algebra morphism
$\theta_c : H_c \to D (\hr ) \ast \WW,
$\label{theta-defn} where $D(Z)$ denotes the {\it ring of
differential operators} on an affine variety $Z$.
Under $\theta_c$ the elements of $\C[\h]$ are identified with the
multiplication operators while, by \cite[p.280]{EG} and in the notation of
\eqref{subsec-3.2},
 $y_i\in\h$ is sent to the {\it Dunkl operator }
\begin{equation} \label{dunkop}\theta_c (y_i) =
\partial_{i} - \sum_{s\in S} c\alpha_s(y_i)\alpha_s^{-1}(1-s),
\qquad \text{where}\  \partial_{i}= \partial/\partial x_i.\end{equation}

Since $\delta$ acts ad-nilpotently on $D(\hr)\ast\WW$,
the set $\{\delta^n\}$ forms an Ore set in that ring.
As observed in  \cite[p.288]{BEGqi}),
$\theta_c$ becomes an isomorphism  on inverting $\delta$;
that is,
 \begin{equation} \label{locdunk} \reg{H_c} = H_c[\delta^{-1}]
\cong D(\hr)\ast \WW.\end{equation}

For any variety $Z$, there is a natural filtration  on $D(Z)$
by order of operators and this induces a filtration on
$D(\hr)\ast \WW$ and its subalgebras by
defining elements of  $\WW$ to have order zero.
If  $R$ is a subalgebra (or submodule) of $D(\hr)\ast \WW$, we write the
operators of order $\leq n$ as $\ord^n(R)$. \label{order-filt-defn}
When $R=H_c$, $\ord$ is clearly
the same filtration as that defined in \eqref{subsec-3.2}.
The  associated graded ring of $R$ will be written
$\ogr(R)=\bigoplus \ogr^n(R)$, where
$\ogr^n(R)=\ord^n(R)/\ord^{n-1}(R)$, and the resulting graded structure of
$\ogr(R)$ will be called the {\it order} or
{\it $\ogr$ gradation}.\label{ogr-defn}
 (This will be only one of several filtrations used in this paper.)

\subsection{}\label{gradingsec}
The rings of differential operators $D(\h)$ and
$D(\hr)$ also have a graded structure, given by the adjoint action
$[\EE,-]$ of the
 {\it Euler operator } $\EE=\sum x_i\partial_{i}\in D(\h)$.
 \label{Euler-defn} We will call this the {\it Euler grading}
 and write $\Edeg$\label{Euler-deg-defn}
 for the corresponding degree function; thus
 $\Edeg x_i=1$ and $\Edeg \partial_{i}=-1$.
 Since $\EE\in D(\h)^\WW$, $\EE$ commutes
 with $\WW$ in $D(\hr)\ast \WW$
 and so this grading extends to that ring with $\Edeg \WW=0$.
 By inspection, \eqref{dunkop} implies that
 the $y_i$ also have degree $-1$ and so each $H_c$ is also graded under
 $[\EE,-]$ and we continue to call this the Euler grading.

It is well-known and easy to check that the $\EE$-grading is compatible with
the order filtration on $D(\hr)\ast \WW$, in the sense that
$[\EE, \ord^n D(\hr)\ast \WW]\subseteq \ord^n D(\hr)\ast \WW$ for all $n\geq 0$.
We therefore obtain an induced grading, again called the $\EE$-grading,
on the associated graded ring $\ogr D(\hr)\ast \WW\cong
\C[\hr\oplus \h^*]\ast W$. Clearly this is
again given by $\Edeg \h^*=1$ (which we define to mean that
$\Edeg(x)=1$ for every   $0\not=x\in \h^*$) while
 $\Edeg \h=-1$ and $\Edeg \WW=0$.

 One should note that, in general,  $\EE\notin H_c$.
 However, there is a natural element in $H_c$ that has the same adjoint action.
   Indeed, let
\begin{equation} \label{hdefn} \hh = \hh_c =
\frac{1}{2} \sum_{i=1}^{n-1} (x_iy_i + y_ix_i) \in H_c.\end{equation}
\label{boldh-defn}
This is independent of the choice of basis.
By \cite[(2.6)]{BEGqi} we have
\begin{equation}\label{grading}
[\hh,x] = x, \quad [\hh, y] = - y, \quad \text{and} \quad
[\hh, w ] = 0 \qquad\mathrm{for \ all}\ x\in\h^*,y\in \h \ \mathrm{and}\ w\in
\WW.
\end{equation}
    Thus commutation with $\hh$ also induces the Euler
grading on $H_c$.

\subsection{The spherical subalgebra}
\label{shiftfunct}
Let $e\in \C \WW$ be the trivial idempotent and
$e_-\in \C\WW$\label{e-defn} be the sign idempotent; thus
$e = |\WW|^{-1} \sum_{w\in \WW} w$  and
 $e_- = |\WW|^{-1} \sum_{w\in\WW} \text{sign}(w)w.$
The main algebra of study in this paper is not the Cherednik algebra
itself, but its {\it spherical subalgebra}  $\UU_c=eH_ce$\label{spherical-defn} and
the related algebra   $\UU^-_c=e_-H_ce_-$.
We will use frequently and without comment that  $\delta$ is a
$\WW$-anti-invariant
and so $e_-\delta =\delta e$.
Also, as $\Edeg \WW=0$, both   $U_c$ and $U^-_c$ have an induced $\EE$-graded
structure.

\subsection{Partitions}\label{dominance} The rest of this section is devoted
to the definition and basic properties of category $\mathcal{O}_c$. Since
its structure depends upon the
combinatorics of $\WW$-representations, we begin with the relevant notions
from  that theory.

We write a
partition of $n$ as $\mu = (\mu_1\geq \mu_2 \geq \cdots \geq \mu_l > 0)$, with
the understanding that $\mu_i =0$ for $i>l$. The \textit{Ferrers diagram} of
$\mu$ is the set of lattice points\label{d-mu-defn}
 $$d(\mu) = \{ (i,j)\in \NN\times \NN : j <
\mu_{i+1}\}.$$ Following the French style, the diagram is drawn with the
$i$-axis vertical
and the $j$-axis horizontal, so the parts of $\mu$ are the lengths of the rows,
 and $(0,0)$ is the lower left corner. The \textit{arm} $a(x)$ and the
\textit{leg} $l(x)$ of a point $x\in d(\mu)$ denote the number of points
strictly to the right of $x$ and above $x$, respectively. The
 {\it hook length} $h(x)$ is $1+a(x)+l(x)$. For example:
\begin{equation}\label{e:arm-leg-pix}
\mu =(5,5,4,3,1)\qquad
\begin{array}[c]{cccccc}
\bullet &   \hbox to 0pt{\hss $\scriptstyle l(x)$\hss }\\
\cline{2-2}
\bullet &   \multicolumn{1}{|c|}{\bullet }& \bullet \\
\bullet &   \multicolumn{1}{|c|}{\bullet }& \bullet &   \bullet \\
\cline{2-5}
\bullet &   \multicolumn{1}{|c|}{\llap{${}_{x}$}\bullet } & \bullet &
    \bullet&     \multicolumn{1}{c|}{\bullet }& {\scriptstyle a(x)} \\
\cline{2-5}
\llap{${}_{(0,0)}$} \bullet &   \bullet &   \bullet &
\bullet &   \bullet
\end{array}
\qquad a(x) = 3,\quad l(x) = 2, \quad h(x)=6.
\end{equation}
The {\it transpose partition $\mu^t$} is obtained from $\mu$ by exchanging
the rows  and columns of $\mu$.

We will always use the
{\it dominance ordering}\label{dominance-defn}
of partitions as in
\cite[p.7]{MacD};
 thus if $\lambda$ and $\mu$
  are partitions of $n$ then $\lambda\geq \mu$ if and only if
$\sum_{i=1}^k \lambda_i \geq \sum_{i=1}^k\mu_i$
 for all $k\geq 1$.

 Let $\irr{\WW}$\label{irred-defn}
denote the set of simple $\WW$-modules, up to isomorphism. As
usual, irreducible representations of $\WW$ will be parametrised
by partitions of $n$.  We use the ordering on $\irr\WW$ arising
from the dominance ordering; thus, as in \cite[Example~1,
p.116]{MacD}, the  {trivial representation}
$\triv$\label{triv-defn}
 is labelled by $(n)$ while the {sign representation} $\sign$\label{sign-defn}
 is parametrised by $(1^n)$ and so $\triv>\sign$. Note that the
 operation on $\irr{\WW}$ given by tensoring by $\sign$ corresponds to
 the transposition of partitions.

\subsection{Category $\mathcal{O}_c$}\label{subsec-3.7}
 (See \cite{GGOR} and  \cite[Definition~2.4]{BEGqi}.)
  Let $\mathcal{O}_c$\label{cat-O-defn}
 be the abelian category of finitely-generated
$H_c$-modules $M$ which are locally nilpotent for the subalgebra
$\mathbb{C}[\h^*]\subset H_c$. By \cite[Theorem 3]{guay} $\mathcal{O}_c$
is a highest weight category.

Given
$\mu \in \irr{\WW}$, we define $\Delta_c(\mu)$,\label{standard-defn} an
object of $\mathcal{O}_c$ called
the \textit{standard module},
 to be the induced module
$\Delta_c(\mu) = H_c\otimes_{\mathbb{C}[\h^*]\ast \WW} \mu,$
where $\mathbb{C}[\h^*]\ast \WW$ acts on $\mu$ by   $pw\cdot
m = p(0) (w\cdot m)$ for $p\in \mathbb{C}[\h^*]$, $w\in \WW$ and
$m\in \mu$.  It is shown in \cite[Section~2]{BEGqi} that each $\Delta_c(\mu)$
has a unique simple quotient $L_c(\mu)$,\label{L-defn}
 that the set $\{ L_c(\mu) : \mu \in
\irr{W}\}$ provides a complete list of non-isomorphic simple objects in
$\mathcal{O}_c$ and that every object in $\mathcal{O}_c$ has finite
length. Note that it follows from the PBW
Theorem~\ref{PBW} that the
standard module $\Delta_c(\mu)$ is a free left
$\mathbb{C}[\h]$-module of rank $\dim(\mu)$.

\subsection{The $\KZ$ functor} \label{subsec-3.11}
Let $M\in \mathcal{O}_c$. Then its localisation $\reg{M}=M[\delta^{-1}]$
at the powers of $\delta$
is a $W$-equivariant $D$-module on $\hr$ in the sense that $\reg{M}$ is
 a $\WW$-equivariant vector
bundle on $\hr$ with a flat $W$-equivariant connection. On taking the germs of
horizontal sections on $\hr/\WW$ we get a representation of the
braid group $B_n = \pi_1 (\hr/\WW)$. This representation
factors through the {\it Hecke algebra} $\hec{q}$\label{hecke-defn}
of $\WW$ with parameter $q= \exp
(2\pi i c)$
 \cite[Theorem 5.13]{GGOR}. In this way we have the
{\it Knizhnik-Zamolodchikov functor}  \label{KZ-defn}
$\KZ :\mathcal{O}_c \rightarrow
\hec{q}\lmod .$
There is an anti-involution $\iota$ on $\hec{q}$ induced by
$\iota(T_w) = T_{w^{-1}}$.  Given a module $V$ for $\hec{q}$, the space
$V^{\ast} = \Hom_{\hec{q}}(V, \C)$ becomes an $\hec{q}$-module via the rule
$h\cdot f (v) = f(\iota(h)v)$.

The images of the standard modules under $\KZ$ are known,
\cite[Remark~6.9 and Corollary~6.10]{GGOR}.  For $c\in \R_{\geq 0}$ and
$\mu\in \irr{\WW}$
 \begin{equation} \label{kzsp} \KZ (\Delta_c(\mu))
\cong Sp_q(\mu)^{\ast}\end{equation}
where $Sp_q(\mu)$\label{specht-defn}
 is the so-called {\it Specht module}
associated to $\mu$. (The dual module appears since the defining
relations for the rational Cherednik algebra given in \cite{GGOR}
are normalised differently to \eqref{rcadef}; as remarked in
\eqref{rcadef}, our parameter $c$  corresponds to $-k_1$ in
\cite{GGOR}.) Now suppose that $M\in \mathcal{O}_c$ has a
filtration
 $$0 = M_0\subset M_1\subset \cdots \subset M_{t-1} \subset M_t =M$$
such that $M_i/M_{i-1}$
is a standard module for all $1\leq i \leq t$. If $N\in \mathcal{O}_c$ and
$c\notin \frac{1}{2}+\Z$ then \cite[Proposition~5.9]{GGOR}
implies that
\begin{equation} \label{kzhom} \Hom_{H_c}(N,M) = \Hom_{\hec{q}}(\KZ(N),
\KZ(M)).\end{equation}

\section{Morita equivalence of Cherednik algebras} \label{shift}

\subsection{}\label{sec301} A powerful technique in the theory of semisimple
Lie algebras is the  translation principle, given by tensoring with a finite
dimensional module, in part because it gives an equivalence of categories
between the  $\OO$ categories (and the  Harish-Chandra categories)
corresponding to distinct  central characters \cite{BG}.  One interpretation of
this is that tensoring with a module of $\mathfrak{k}$-finite vectors gives a
Morita equivalence  between the corresponding factors of the  enveloping
algebra  \cite[Corollary~4.13]{JS}.

Although it does not involve finite dimensional modules,
there is a natural analogue of this procedure for Cherednik algebras,
given by the Heckman-Opdam shift functors defined in \eqref{shift-defn}.
These functors have proved useful in a number of papers (see,
for example, \cite{BEGqi, BEGfd,gordc}) and for
 particular values of $c$ these functors are known to give
equivalences of categories between   $H_c$, $U_{c}$ and $U_{c+1}$
(see, for example, \cite[Theorem~8.1]{BEGqi} and \cite[Proposition~4.3]{BEGfd}).
It is an open problem to determine precisely when these equivalences exist
\cite[Remark~5.17]{GGOR} and  this question is crucial to our
$\Z$-algebra approach to  Cherednik algebras. We give an
essentially complete answer to this question in Corollary~\ref{morrat-cor}
and Remark~\ref{morrat-cor-remark}.
We also prove   that the equivalence
$H_c\to H_{c+1}$ maps  category $\mathcal{O}_c$ to
  $\mathcal{O}_{c+1}$ and sends the standard module
$\Delta_c(\mu)$  to $\Delta_{c+1}(\mu)$, see Proposition~\ref{shiftonO}.

\subsection{}\label{subsec-4.0}
 Fix $c\in \C$ and keep the notation of \eqref{e-defn}.
 If we identify $H_{c}$  with its
  image in $D (\hr)\ast W$ via the Dunkl operator \eqref{dunkop}
  then, by
 \cite[Proposition~4.1]{BEGfd}, there is an identity
\begin{equation} \label{conj} \UU_c \ =\ \delta^{-1} \UU^-_{c+1}\delta
\ = \ e\delta^{-1} H_{c+1}\delta e.
\end{equation}
In particular, \label{Q-defn}
$Q_c^{c+1}= eH_{c+1}e_-\delta
= eH_{c+1}\delta e$ is a $(\UU_{c+1}, \UU_c)$-bi-submodule of $D(\hr)\ast \WW$.
The  shift functors mentioned above are given by
  \label{shift-defn}
$$S_c: \UU_c \lmod \to \UU_{c+1}\lmod: \qquad
N\mapsto  Q_c^{c+1} \otimes_{\UU_c} N
$$
and
 $$\widetilde{S}_c: H_{c} \lmod \to H_{c+1}\lmod: \qquad
M\mapsto  H_{c+1}e_-\delta \otimes_{\UU_{c}} eM.$$

\subsection{}\label{subsec-4.1}   When $c$ is a positive real number, the Morita
equivalence between  $U_c$ and $ U_{c+1}$  is given by $S_c$ and we begin with
that case. The general case, proved in Corollary~\ref{morrat-cor}, will be an
easy consequence.

\begin{thm} \label{morrat}  Assume that $c\in \R_{\geq 0}$ with  $c\notin
\frac{1}{2} + \mathbb{Z}$.
Then  both shift functors $\widetilde{S}_c: H_{c}\lmod \to H_{c+1}\lmod$ and
$S_c: \UU_c\lmod \to \UU_{c+1}\lmod $ are Morita
equivalences.

Moreover, the idempotent functor $E_c: H_{c} \lmod\to \UU_c\lmod$ given by
$M\mapsto eM$ is a Morita equivalence. \end{thm}

\begin{proof} In order to prove that $S_c$ is an equivalence, we need to
   show that $Q=Q_c^{c+1}$ is a projective
generator for $\UU_{c+1}\lmod$, with endomorphism ring
$\mathrm{End}_{U_{c+1}}(Q) =U_{c}$.
 Arguing as in \cite[Theorem~1.5(iv)]{EG} the dual
$Q^\ast=\Hom_{\UU_{c+1}}(Q,\UU_{c+1})$   is
  $P = \delta^{-1} e_-H_{c+1}e.$
 By the dual basis lemma,  $Q$ is a projective
  $\UU_{c+1}$-module with  $\mathrm{End}_{U_{c+1}}(Q) =
U_{c}$ if and only if
 $PQ=\UU_c$ while $Q$ is a generator if and only if $QP=\UU_{c+1}$.
Substituting in the given formul\ae\ for $Q$ and $P$
shows that we need to prove that
\begin{equation}\label{morrat1}
H_{c+1}e_-H_{c+1} =H_{c+1}\qquad\text{and} \qquad
H_{c+1}eH_{c+1} = H_{c+1}\quad\text{for}\  c\geq 0.
\end{equation}
Similarly, as $H_{c}e$ is a projective left $H_{c}$-module, $E_c$ will
  be a Morita equivalence if we prove that
\begin{equation}\label{morrat11}
H_ceH_c = H_c \quad \text{for}\ c\geq 0.\end{equation}
Since $\widetilde{S}_{c}=E^{-1}_{c+1}\circ S_{c}\circ E_{c}$,
Equations~\ref{morrat1} and \ref{morrat11} will suffice to
prove the theorem.

  The proof of  Theorem~\ref{morrat}  will be through a series of lemmas
  and we begin with the first equality in \eqref{morrat1}.
    Set $d=c+1$; thus $d\in \R_{\geq 1}$, with $d\notin \frac{1}{2}+\Z$.

 \subsection{Reduction to Category $\OO$}\label{subsec-4.2}
If $H_de_-H_d$ is a \textit{proper} two-sided
ideal of $H_d$ it must be contained in a primitive ideal, and
hence, by \cite[Generalized Duflo Theorem]{ginz}, annihilate an
object from category $\mathcal{O}_d$. Thus it is enough to show
that $e_-$ does not annihilate any
simple module belonging to $\mathcal{O}_d$.

To do this we first show in Corollary~\ref{poono}  that  the
composition factors of $\Delta_d(\mu)$ are of the form $L_d(\lambda)$ for
$\lambda \leq \mu$. Under the $\Z$-strings ordering such a result is proved in
\cite{guay} but as we work with the dominance ordering of partitions and representations, as defined in
\eqref{dominance-defn}, this definitely requires work, see also
 \eqref{order-remark}.   We then show that
the lowest weight copy of the sign module for
$\WW$ in $\Delta_d(\mu)$ does not occur in any  standard module
$\Delta_d(\lambda)$ for $\lambda < \mu$. Since $L_d(\mu)$ is the head
 (that is, the unique simple factor module)  of $\Delta_d(\mu)$ it
will follow that $e_-L_d(\mu)\neq 0$.

\subsection{Lemma.}
\label{basiccom}   {\it
Let $c\in \R_{\geq 0}$ with
$c\notin \frac{1}{2}+\mathbb{Z}$.
If  $\mathrm{Hom}_{H_c}(\Delta_c(\lambda),\,\Delta_c(\mu))\not=0$
for $\lambda,\mu\in \irr{\WW}$,
 then $\lambda \leq \mu$ in the dominance ordering.
}

\begin{proof}
Let $S_q = S_q(n,n)$\label{schur-defn}
 be the $q$-Schur algebra defined in \cite[Section~1]{DJ}, where $q =  \exp(2\pi i c)$.
It is conjectured in \cite[Remark~5.17]{GGOR} that $S_q\lmod$
is equivalent to $\OO_c$. We
cannot prove this, but we will show that there is a relationship which implies
the lemma.

For each $\mu\in \irr{\WW}$ there is an $S_q$-module $W_q(\mu)$,
\label{q-weyl-defn}
 called the {\it $q$-Weyl module}. By \cite[Corollary~8.6]{DJ}, there is an
isomorphism
\begin{equation} \label{qsch}\Hom_{\mathcal{H}_q}(Sp_q(\mu), Sp_q(\lambda))
\cong
\Hom_{S_q}(W_q(\lambda),W_q(\mu)). \end{equation}

On the other hand, by \eqref{kzsp} and \eqref{kzhom}  we have
\begin{equation} \label{tf} \Hom_{H_c}(\Delta_c(\lambda), \Delta_c(\mu)) \cong
\Hom_{\mathcal{H}_q}(Sp_q(\lambda)^{\ast}, Sp_q(\mu)^{\ast}) \cong
\Hom_{\hec{q}}(Sp_q(\mu), Sp_q(\lambda)).\end{equation}

Each $W_q(\nu)$ has a simple head $F_q(\nu)$,\label{F-defn}
 \cite[Theorem~4.6]{DJ} and $\{
F_q(\nu): \nu\in \irr{\WW}\}$ is a complete, repetition-free list
of the simple $S_q$-modules up to isomorphism, \cite[Theorem~8.8]{DJ}.
Furthermore, $F_q(\lambda)$ is a composition factor of $W_q(\mu)$
only if $\lambda\leq
\mu$, \cite[Corollary~8.9]{DJ}.  By \eqref{qsch} and \eqref{tf} a non-zero
homomorphism  $\phi: \Delta_c(\lambda)\to \Delta_c(\mu)$ implies the
existence of a non-zero homomorphism $\phi ':W_q(\lambda)\to W_q(\mu)$.  Thus
$F_q(\lambda)$ must be a composition factor of $W_q(\mu)$ and so $\lambda\leq
\mu$.
\end{proof}

\subsection{Corollary} \label{poono}
 {\it Assume that $c\in \R_{\geq 0}$, with
$c\notin \frac{1}{2}+\mathbb{Z}$.
If  $ [\Delta_c(\mu) : L_c(\lambda) ] \neq 0$
for $\lambda,\mu\in \irr{\WW}$, then $ \lambda \leq \mu$
in the dominance ordering.
}

\begin{remarks} (1) For arbitrary $c$ and $\mu$,  the unique
 occurrence of $L_c(\mu)$ as a composition factor of
 $\Delta_c(\mu)$ is as its head---see, for example,
the discussion after Lemma~7 in  \cite[Section~2]{guay}.

(2) Since $\sign$ is minimal in the dominance ordering, the lemma and the
above remark imply that $\Delta_c(\sign)$ is irreducible for all
$c\in \R_{\geq 0}$. This can also be deduced from \cite{guay}.
\end{remarks}

\begin{proof} We argue by induction on $\mu$.
More precisely, suppose that $[\Delta_c(\mu): L_c(\lambda)] \neq 0$ for
some $\mu \neq \lambda$ and that  the lemma holds for any $\nu < \mu$.
(The induction starts since there are only finitely many $\sigma$ with
 $\sigma<\mu$.)  Let $P_c(\lambda)$ \label{P-defn}
denote the projective cover of $\Delta_c(\lambda)$, as in \cite[Section~3.5]{GGOR},
and write $K$ for the kernel of the associated homomorphism
$\phi: P_c(\lambda) \to \Delta_c(\mu)$.
By \cite[Proposition~13]{guay}
there is a $\Delta$-filtration of $P_c(\lambda)$ $$P_c(\nu) = M_0 \supset M_1
\supset \cdots \supset M_t = 0$$ with each factor $M_j/M_{j+1}$  of the form
$\Delta_c(\lambda_j)$ for some $\lambda_j\in \irr{\WW}$.
 Thus there exists $i$ such that
$M_i+K/K \neq 0$ but $M_{i+1}+K/K = 0$. This gives a non-zero composition
$$\psi: \Delta_c(\lambda_i) \cong
M_i / M_{i+1} \longrightarrow (M_i + K)/K \longrightarrow
P_c(\lambda)/K \longrightarrow \Delta_c(\mu).$$
By Lemma~\ref{basiccom}, $\lambda_i \leq\mu$. If $\lambda_i=\mu$ then
the first remark after the statement of the lemma would imply that $\psi$ and
hence $\phi$ are surjective, contradicting the fact that $\lambda\not=\mu$.
Thus $\lambda_i<\mu$.
By BGG reciprocity \cite[Theorem~19]{guay},
$[P_c(\lambda): \Delta_c(\lambda_i)] =
[\Delta_c(\lambda_i): L_c(\lambda)]\neq 0$ and so,
 by induction, $\lambda \leq \lambda_i$. Thus $\lambda < \mu$.
 \end{proof}

\subsection{}\label{order-remark}
A result analogous to Corollary~\ref{poono}
 is proved as part of  the proof of \cite[Proposition~13]{guay}.
  However the $\Z$-strings ordering used in
\cite{guay} is different from the dominance ordering.
An explicit example where the orderings differ
 can be found when $n=8$, by taking $\lambda = (6,1,1)$ and $\mu = (4,4)$.
 In this case $\lambda$ and $\mu$ are incomparable in the dominance ordering,
 but comparable in the $\Z$-strings ordering.

\subsection{The canonical grading on $\mathcal{O}_c$.}
\label{cangrad}  The final ingredient we need for the proof of
Theorem~\ref{morrat} is a canonical  grading on $\OO_c$. Let
 $\hh_c\in H_c$ be defined as
in \eqref{hdefn}. Then, for $M\in \OO_c$ and $\alpha\in\C$, set
$$W_{\alpha}(M) = \{ m\in M : (\hh_c - \alpha)^km = 0
 \text{ for $k\gg
0$}\}.$$\label{cangrad-defn}
By \cite[(2.4.1)]{GGOR} this gives the
\textit{canonical grading} $M = \sum_{\alpha\in \C} W_{\alpha}(M).$

This observation has two useful consequences. First, if
$\theta:M_1\rightarrow M_2$ is
 an $H_c$-module homomorphism with  $M_i\in \OO_c$, then $\theta
(W_{\alpha}(M_1)) \subseteq W_{\alpha}(M_2)$ for each $\alpha\in\C$. Secondly,
if $p\in H_c$ has $\Edeg p = t$, then   \eqref{grading},
implies that   $p\cdot W_{\alpha}(M) \subseteq
W_{\alpha+t}(M)$. Note that the standard module $\Delta_c(\mu)$ is
therefore a lowest weight module
since it is generated as a $\C[\h]$-module by the space $1\otimes \mu$.

\subsection{} \label{fkdeg}
To describe the graded structure of the standard
modules we need a little notation. Recall that the space of coinvariants
$\C[\h]^{\text{co} \WW} =  {\C[\h]}\big/{\C[\h]_+^{\WW}\C[\h]}$
 is a finite dimensional graded
algebra isomorphic as  a $\WW$-module to the regular
representation. As in \cite{op}, the polynomials
\begin{equation} \label{fakedegrees} f_{\mu}(v) = \sum_{i\geq 0}
[\C[\h]^{\text{co} \WW}_i : \mu] v^i \end{equation} are called the
\textit{fake degrees}\label{fake-defn}
 of $\mu\in\irr{\WW}$.
 We define $n(\mu)$ to be the lowest power of $v$ appearing
  in $f_{\mu}(v)$; thus, $f_{\mu}(v) = a
v^{n(\mu) }+ \text{higher order terms.}$ In the notation of \cite{hai3},
$n(\mu)$ is equal to the
{\it partition statistic}  $  \sum_{i}\mu_i (i-1)$ (see   the proof of
 \cite[Theorem~6.4]{babyv}).
Finally, \eqref{fakedegrees} implies that
\begin{equation}\label{fakedegrees2}
 f_{\mu^t}(1)=\dim \mu^t = \dim \mu = f_\mu(1)
\qquad\rm{for}\ \mu\in \irr\WW.\end{equation}

\subsection{}\label{subsec-3.10}  Given a graded $\WW$-module $M=\sum_{\alpha\in\C}
W_\alpha(M)$ we define its {\it graded Poincar\'e series} to
be $$p(M,v,\WW) = \sum_{\alpha\in \C} v^\alpha \sum_{\lambda\in \irr{\WW}}
[W_\alpha(M) : \lambda][\lambda].$$\label{W-poincare}
This  is easily determined for standard modules.

\begin{prop} {\rm (1)} Under the canonical grading, the
subspace $1\otimes \mu$ of $\Delta_c(\mu)$
has weight $m+c(n(\mu) -
n(\mu^{t}))$, where $m=(n-1)/2$.

 {\rm(2)} The Poincar\'{e} series of $\Delta_c(\mu)$ as a graded
$\WW$-module is
\begin{equation}
\label{polystand} p(\Delta_c(\mu), v, \WW) =
v^{m+c(n(\mu)-n(\mu^t))} \frac{\sum_{\lambda}f_{\lambda}(v)
[\lambda\otimes \mu]}{\prod_{i=2}^n (1-v^i)}.\end{equation}
\end{prop}

\begin{proof} (1) We need to compute the action of
$\hh=\frac{1}{2}\sum_{i=1}^{n-1}(x_iy_i+y_ix_i)$ on the
 space $1\otimes \mu$. By the defining relations of $H_c$
 from \eqref{hc-defn},
and the fact that the $\{x_i\}$ and  $\{y_i\}$ are dual bases, we obtain
$$\begin{array}{rl}
\hh = \sum_i x_iy_i + (n-1)/2 -\frac{1}{2}\sum_{s\in \mathcal{S}}
\sum_{i}c\alpha_s(y_i)x_i(\alpha_s^\vee)s
=& \sum x_iy_i + (n-1)/2 -\frac{c}{2}\sum_{s\in \mathcal{S}}
\alpha_s(\alpha_s^\vee)s\\
 \noalign{\vskip 5pt}
 =&\sum x_iy_i + m -c\sum_{s\in \mathcal{S}}s.
 \end{array}$$
 The action of $\sum(1-s)$ on $\lambda\in \irr{\WW}$
 can be derived  from \cite{BM,Lu}. More precisely,
 $\lambda$ is  special by \cite[(4.2.2)]{Lu} and so $n(\lambda)=
 b_\lambda=a_\lambda$
 in the notation of \cite{Lu}. Therefore,
 by \cite[Section~4.21]{BM}  and \cite[Section~4.1 and (5.11.5)]{Lu},
$\sum_s(1-s)$ acts on
$\lambda\in \irr{\WW}$ with weight $N+n(\lambda)-n(\lambda^t)$, where
$N=n(n-1)/2$ is the cardinality of $\mathcal{S}$.
Thus $\sum_ss$ acts on $1\otimes \mu$
 with weight $-(n(\mu) - n(\mu^t))$ and hence
 $\hh$ acts  with weight
$m+c(n(\mu) - n(\mu^t))$.

(2) As graded $\WW$-modules, $\Delta(\mu)\cong (\C[\h]\otimes \mu)[k]$
for $k=m+c(n(\mu) - n(\mu^t))$. The shift arises from the fact that,
  by (1),   the generator $1\otimes \mu$ of $\Delta_c(\mu)$
  lives in degree $k$. The Chevalley-Shephard-Todd Theorem
  implies that, as  graded $\WW$-modules, $\C[\h]\cong
  \C[\h]^\WW \otimes \C[\h]^{\text{co} \WW}$.
  Now $\C[\h]^\WW$ is a polynomial ring with generators in
  degrees $2\leq i\leq n$ and so its
   Poincar\'e polynomial is
$ \prod_{i=2}^n (1-v^{i})^{-1}$. On the other hand,
the coinvariant ring $ \C[\h]^{\text{co}\WW}$ has graded
Poincar\'e polynomial
$\sum_{\lambda} \sum_i [\C[\h]_i^{\text{co} \WW} : \lambda] [\lambda]v^i$.
By definition, this is just $\sum_{\lambda} f_{\lambda}(v)[\lambda]$. Combining
these observations gives  \eqref{polystand}.
\end{proof}

\subsection{Completion of the proof of Theorem~\ref{morrat}}
\label{subsec-4.5}  We first prove that $H_de_-H_d=H_d$
(where $d=c+1$, as before).
Since $\mu\cong \mu^*$ for symmetric groups,
 the sign representation is a direct summand of
$\mu\otimes \nu$ if and only if $\nu =\mu^t$.
Thus \eqref{polystand} implies
that $\sign$  first appears  in $\Delta_d(\mu)$
in the weight space
$$m+d(n(\mu)-n(\mu^t)) + n(\mu^t) = m+ dn({\mu}) - (d-1)n({\mu^t})
\qquad\text{where}\quad m=(n-1)/2.$$
 If $\lambda \leq \mu$ then $n({\lambda}) \geq n(\mu)$ by
 \cite[Theorem~B and Proposition~1.6]{Shi}.
 Moreover, as $\lambda^t\geq \mu^t$,
 we have $n(\lambda^t)\leq
n(\mu^t)$. Since $d\in \mathbb{R}_{\geq 1}$,
 $$m+dn({\lambda}) -(d-1)n({\lambda^t}) \geq m+dn({\mu}) -(d-1)n({\mu^t}) $$
 with equality if and only if $ \lambda=\mu$.

It follows that the copy of  $\sign$  appearing in the lowest possible weight
space of $\Delta_d(\mu)$ is never a weight of $\Delta_d(\lambda)$ for $\lambda
< \mu$. By  Corollary~\ref {poono}, this means that this
copy of $\sign$ is a weight for
$L_d(\mu)$ and hence that $e_-L_d(\mu)\not=0$. By
\eqref{subsec-4.2} this implies that $H_de_-H_d =H_d$, and so the
first equality of \eqref{morrat1} is proven.

It remains to show that $H_ceH_c=H_c$ for $c\in \R_{\geq 0}$.
 The argument of \eqref{subsec-4.2}
shows that we need to prove that $e$ does not annihilate any simple module
from $\mathcal{O}_c$.  In this case
$\triv$ appears in $\mu\otimes \nu$  if and only if $\nu=\mu$.
Therefore, \eqref{polystand} now implies that $\triv$ first
appears in $\Delta_c(\mu)$ in degree $m+c(n(\mu) -n(\mu^t))+n(\mu).$
Let $\lambda\leq \mu$. Then
$$m +c(n(\lambda)-n(\lambda^t))+n(\lambda)\
=\ m+(c+1)n(\lambda) -cn(\lambda^t) \
  \geq\  m+c(n(\mu)-n(\mu^t))+n(\mu),$$
with equality if and only if $\lambda=\mu$. This means that $\triv$ appears in
$\Delta_c(\lambda)$ in a higher degree than its first appearance in
 $\Delta_c(\mu)$.
In particular, the simple quotient $L_c(\mu)$ of $\Delta_c(\mu)$
 contains a copy
of $\triv$ and so it cannot be annihilated by $e$.
This therefore completes the proof of \eqref{morrat1} and \eqref{morrat11}
and hence proves the theorem.
  \end{proof}

  \subsection{General equivalences}\label{subsec-4.55}
  We now give the promised extension of Theorem~\ref{morrat}
  to more general values of   $c$.
 Since it requires no extra work, and it is put to crucial use in
  \cite{BFG}, we will also prove the
 result over more general base fields. Thus if $k$ is a
 subfield of $\C$, with $c\in k$, let $H(k)_c$
 denote the $k$-algebra defined by the
 generators and relations from \eqref{hc-defn}.
 We write $U(k)_c$, $Q(k)_c^{c+1}$,
 etc, for the corresponding
 objects defined over~$k$.

  \begin{hypothesis}\label{morrat-hyp}
  Set
  $\mathcal{C} = \{z: z=\frac{m}{d}\ \mathrm{where}\
  m, d \in \Z\text{ with }  2\leq d\leq n \text{ and } z\notin \Z\}.$
  Assume that $c\in \C$ is such that $c\notin \frac{1}{2} +
\mathbb{Z}$.  If $c$ is a rational number with  $-1<c<0$
assume further  that $c\not\in \mathcal{C}$.
\end{hypothesis}

 \subsection{Corollary}\label{morrat-cor} {\it Let $k\subseteq \C$ be a field
 and assume that $c\in k$ satisfies Hypothesis~\ref{morrat-hyp}.

 {\rm (1)}  $U(k)_c$    and $H(k)_{c}$
are   Morita equivalent. If $c \notin  (-2,-1)_{\mathcal{C}}
=\{z\in \mathcal{C} : -2<z<0\}$, then
$U(k)_c$ is Morita equivalent to $U(k)_{c+1}$.

{\rm (2)}  Let $a=-c$.
Then   $H(k)_{a}$ is
  Morita equivalent to $U(k)^-_{a} = e_-H(k)_{a}e_-$.
  If $a \notin  (1,2)_{\mathcal{C}}$, then
$U(k)^-_{a}$ is Morita equivalent to $U(k)^-_{a-1}$. }

\smallskip
\begin{proof}  (1) We start with the case $U_c=U(\C)_c$.
If $c\not\in \mathcal{C}$
then it follows from \cite[Theorem~8.1]{BEGqi}
and \cite[Theorem~4.3]{DJ2}
 that   $H_c$, $U_c$  and $U_c^-$
are simple, Morita equivalent rings
(see the introduction to \cite{BEG3} for the details).
Since this also  applies to $H_{c+1}$ the conditions \eqref{morrat1}
are trivially satisfied and the result follows.

Thus we may assume that
 $c\in \mathcal{C}$. If
 $c\geq -1$, then necessarily $c\geq 0$ and so the result follows from
 Proposition~\ref{morrat}. Otherwise $c\leq -1$. In this case
  the discussion before  \cite[Remark~2.2]{De} shows that
   there is an isomorphism
$\chi: H_c\to H_{-c}$ satisfying $\chi(e_-)=e$.
 Thus, for any $c$, \eqref{conj}
implies that  $U_c\cong U^-_{-c} \cong eH_{-c-1}e=U_{-c-1}$.
The result for $c\leq -1$ therefore follows from the cases already discussed.

Finally, let $k$ be an arbitrary subfield of $\C$ and consider $U(k)_c$.
In order to prove, for example, that $U(k)_c$ is Morita equivalent to
$U(k)_{c+1}$
we need to prove that $Q(k)P(k)=U(k)_{c+1}$ and $P(k)Q(k)=U(k)_c$.
By construction, $Q(\C)=Q(k)\otimes_k\C$, and similarly for $P(\C)$.
By the earlier part of the proof, $U(\C)_c/P(\C)Q(\C)=0$.
The faithful flatness of $U(\C)_c=U(k)_c\otimes_k\C$ as a
$U(k)_c$-module then ensures that $U(k)_c/P(k)Q(k)=0$, whence $PQ=0$.
All the other steps in the proof follow in exactly the same way.

(2) Using the identity $U_c\cong U^-_{-c}$, this follows from part (1).
\end{proof}

\subsection{Remarks}\label{morrat-cor-remark}
(1) The condition that $c\notin \frac{1}{2} + \mathbb{Z}$ is needed
in Theorem~\ref{morrat} and Corollary~\ref{morrat-cor} in order
to apply \eqref{kzhom} and may be unnecessary.
This is the case when   $n=2$ as $U_c$ is Morita equivalent to $U_{c+1} $ if and only if
$c\not= -\frac{3}{2}, -\frac{1}{2}$ (see, for example,
\cite[Proposition~8.2]{EG}). The point about the excluded cases  is that
$U_{-\frac{1}{2}}$ is simple but the two neighbouring algebras,
$U_{\frac{1}{2}}$, $U_{-\frac{3}{2}}$ are not.
 Combining \cite[Proposition~8.2]{EG}  with \cite[Theorem~B]{St}
  shows that $U_{-\frac{1}{2}}$ has infinite global dimension,
  and so the next
Corollary~\ref{gldim} also fails for this value of~$c$.

(2) This also shows that the  hypothesis
$c\notin (-2,0)_{\mathcal{C}}$ is  serious.
Indeed, for any $n\geq 2$,  let
$c=-m/n\in (-1,0)_{\mathcal C}$. Then one can prove that the
factor module $V_c =\Delta_c(\sign)/I_c$ considered in
\cite[Theorem~3.2]{CE}  does not contain a copy of the  $\WW$-module
$\triv$ (we thank Pavel Etingof for this fact).
In particular $eV_c=0$ and so {\it the functor $E_c$ is not an
equivalence}. If we further assume that $(m,n)=1$, then $V_c$ is
the unique irreducible finite dimensional $H_c$-module by
\cite[Corollary~3.3]{CE} and \cite[Theorem~1.2(ii)]{BEGfd}. Since
$U_c=\mathrm{End}_{H_c}(eH_c)$, this implies that
  $U_c$ has no finite dimensional modules. However, by
  Corollary~\ref{morrat-cor}(1) and  \cite[Theorem~1.2]{BEGfd}
  $U_{c\pm 1}\lmod$ does have such modules and so {\it there
is no equivalence between $U_c $ and $U_{c\pm 1} $}.

\subsection{Corollary}\label{gldim}
{\it Assume that $c\in \C$ satisfies Hypothesis~\ref{morrat-hyp}.
Then $H_c$ and $U_c$ have finite homological global dimension
and satisfy the Auslander-Gorenstein conditions and Cohen-Macaulay
conditions of \cite{lev}.}

\begin{proof} Since this result takes us a little far afield, the details of
the proof are left to the interested reader. Standard techniques show
that  $\ogr H_c\cong \cxy\ast\WW$ and hence $H_c$ have the given properties
 (see, for example, \cite[Theorem~4.4]{Br}). By
Corollary~\ref{morrat-cor}, $U_c$ is Morita equivalent to $H_c$  and it follows
that $U_c$ also has these  properties. \end{proof}

\subsection{The shift functor on $\mathcal{O}_c$}\label{subsec-4.6} Many
 computations for $\UU_c$ reduce to computations in category
 $\mathcal{O}$ and so it is important to know that, under the
 hypotheses of Theorem~\ref{morrat},  $S_c$ does provide an equivalence
 between the corresponding categories. This is the point of the next result.

\begin{prop}\label{shiftonO} Assume that $c\in \C$ satisfies
Hypothesis~\ref{morrat-hyp} and that $c\notin \Q_{\leq -1}$. Then the
shift functor $\widetilde{S}_c$ restricts to an
equivalence between $\mathcal{O}_{c}$ and $\mathcal{O}_{c+1}$ such
that $\widetilde{S}_c(\Delta_{c}(\lambda)) \cong \Delta_{c+1}(\lambda)$ for
all partitions $\lambda$ of $n$. Thus $S_c(e\Delta_{c}(\lambda))
=e\Delta_{c+1}(\lambda)$.\end{prop}

\begin{remark} By Corollary~\ref{morrat-cor}(2),  an analogue of the
proposition also holds when  $c\in \mathbb{Q}_{\leq -1}$,
provided that one  shifts in a negative
direction.
\end{remark}

\begin{proof}  The final assertion of the proposition
is immediate from the previous one
combined with Corollary~\ref{morrat-cor}(1).

 We begin by showing that $\widetilde{S}_c$ restricts to an
equivalence between $\mathcal{O}_{c}$ and $\mathcal{O}_{c+1}$.
Fix $M\in \mathcal{O}_{c}$.
Let $\mathcal{I}_t=\C[\h^*]^{\WW}_{\geq t}$ denote  the $\WW$-invariant
elements of $\C[\h^*]$ of degree at least $t$ and set
$I_t=\mathcal{I}_t\C[\h^\ast ]$, Then $\C[\h^\ast]/I_t$ is a finite dimensional
algebra and so all homogeneous elements of $\C[\h^*]$ of large degree
belong to $I_t$.
Thus it is enough to show that, if  $\widetilde{m}=he_-  \delta\otimes em \in
\widetilde{S}_c(M)= H_{c+1}e_-\delta \otimes_{\UU_c} eM$,
 for some $h\in H_{c+1}$ and $m\in M$,  then
 $\widetilde{m}$  is annihilated by $\mathcal{I}_t$ for $t\gg 0$.

Recall the $\EE$-grading on $H_c$ from \eqref{gradingsec}.
 Since $\C[\h^*]$ acts  locally
nilpotently on $M$, the PBW
isomorphism \eqref{PBW}  shows that any
homogeneous element of $H_c$ of sufficiently large
negative $\EE$-degree annihilates $m\in M$. Thus, assume that
$qm=0$ for all $q\in H_c$  with $\Edeg(q)\leq -t$ and let
$p\in \C[\h^*]^{\WW}_{\geq t}$.  Then
$$phe_-\delta \otimes em \ =\ [p,h]e_-\delta  \otimes em +
 h\delta \delta^{-1} pe_- \delta \otimes em \ =\ [p,h]e_-\delta \otimes em +
he_- \delta \otimes \delta^{-1} p \delta em.$$
Since $\Edeg \delta^{-1}p \delta=\Edeg p \leq -t$,
 we have $\delta^{-1} p\delta em = 0$ by the hypothesis on $t$.

 Therefore $p(he_- \delta \otimes em) = [p,h]e_- \delta \otimes em$
 for any such $p$.  Since the choice of $t$ was independent of $h$,
 this implies that $p^r (he_- \delta \otimes em) =
  \ad(p)^r(h)(e_- \delta \otimes em)$, for any $r>0$.
 Now, $p$ commutes with both $\C[\WW]$ and
 $\C[\h^*]$, and so  the defining relations
 of $H_{c+1}$ from \eqref{rcadef} ensure that
 the adjoint action of $p\in \C[\h^*]^{\WW}$ on $H_{c+1}$
 is locally nilpotent (see also \cite[Lemma~3.3(v)]{BEGqi}).
 Therefore a sufficiently large power of $p$ annihilates
 $he_-\delta \otimes em$.
 Thus $\widetilde{S}_c(M)\in \mathcal{O}_{c+1}$
  and $\widetilde{S}_c$ does restrict to the desired equivalence.

It remains to compute $\widetilde{S}_c(\Delta_{c}(\lambda))$ and we begin
with the analogous problem on $\reg{H_{c+1}}$.
 In the notation of \eqref{locdunk},
 $$\reg{\widetilde{S}_c(\Delta_{c}(\lambda))} =
\reg{H_{c+1}}e_-\delta \otimes_{\delta^{-1}\UU^-_{c+1}\delta }
e\Delta_c(\lambda).$$
By \eqref{locdunk},  $\reg{H_{c+1}} \cong A=D(\hr)\ast W$ and  so
$\reg{\widetilde{S}_c(\Delta_{c}(\lambda))}
\cong Ae_-\delta \otimes_{B} \reg{e\Delta_{c} (\lambda)},$
where $B=\delta^{-1} e_-Ae_-\delta $.
On the other hand, \eqref{conj} induces an isomorphism
$$\theta: Ae_-\delta
\otimes_{B} e\reg{\Delta_c(\lambda)}
 \longrightarrow Ae \otimes_{eAe} e\reg{\Delta_c(\lambda)};
\qquad ae_-\delta \otimes em\mapsto  a\delta e\otimes em.$$
Combined with the identity $H_ceH_c=H_c$ from Corollary~\ref{morrat-cor}(1),
this implies that
 \begin{equation}\label{nonzeromap1}
 \reg{\widetilde{S}_c(\Delta_{c}(\lambda))} \cong
 A e \otimes_{eA e} \reg{e\Delta_{c}(\lambda)}
  \cong \reg{\left(H_{c}e \otimes_{\UU_c} e\Delta_{c}(\lambda)\right)}
   \cong\reg{\Delta_{c}(\lambda)}\not=0.
   \end{equation}

If    $c\not\in \mathcal{C}$, we are done. Indeed, in this case
 \cite[Corollary~2.11]{BEGqi}
implies that $\Delta_{c+1}(\lambda)$,
$\Delta_{c}(\lambda)$ and hence
$\widetilde{S}_c(\Delta_{c}(\lambda))$ are  all simple modules.
The isomorphism \eqref{nonzeromap1} implies that
$\widetilde{S}_c(\Delta_{c}(\lambda))\hookrightarrow
 \reg{\Delta_{c+1}(\lambda)}$. Under this embedding,
$\widetilde{S}_c(\Delta_{c}(\lambda))\cap
 \Delta_{c+1}(\lambda)\not=0$ and hence $\widetilde{S}_c(\Delta_{c}(\lambda))=
 \Delta_{c+1}(\lambda)$.

 We may therefore assume that $c\in \mathcal{C}$, in which case
   Hypothesis~\ref{morrat-hyp} implies that $c\geq 0$ and
   we can use the KZ-functor from \eqref{KZ-defn}. By
   \eqref{nonzeromap1} and  \eqref{kzsp},
$
\KZ(\widetilde{S}_c(\Delta_{c}(\lambda)))
\cong \KZ(\Delta_{c}(\lambda)) \cong
Sp_q(\lambda)^{\ast}.$
By \eqref{kzhom} and \eqref{qsch} we therefore  have
\begin{equation}\label{nonzeromap}
\Hom_{H_{c+1}}(\widetilde{S}_c(\Delta_{c}(\lambda)),
\Delta_{c+1}(\lambda)) \cong \Hom_{S_q}(W_q(\lambda), W_q(\lambda)) = \CC.
\end{equation}
It follows from Corollary~\ref{poono}  that the composition
factors of $\Delta_{c+1}(\lambda)$ are of the form $L_{c+1}(\nu)$ with
$\nu\leq \lambda$ in the dominance ordering. We will show by an ascending
induction on this ordering that $\widetilde{S}_c(\Delta_{c}(\lambda)) \cong
\Delta_{c+1}(\lambda)$.

If $\lambda$ is minimal in the dominance ordering then
$\lambda=\sign$ and so  both $\Delta_{c+1}(\lambda)$ and
$\widetilde{S}_c(\Delta_{c}(\lambda))$ are simple
by Remark~\ref{poono}.  By \eqref{nonzeromap} there is a non-zero
map from $\widetilde{S}_c(\Delta_{c}(\lambda))$ to
$\Delta_{c+1}(\lambda)$ which therefore must
be an isomorphism. This begins the induction.

Let $\lambda$ be arbitrary and suppose that, for all
 $\nu < \lambda$ in the dominance ordering, we have
$\widetilde{S}_c(\Delta_{c}(\nu))\cong\Delta_{c+1}(\nu)$,
and hence  that $\widetilde{S}_c(L_{c}(\nu)) \cong L_{c+1}(\nu)$. Since
$\widetilde{S}_c$ is an equivalence, $\widetilde{S}_c(\Delta_{c}(\mu))$
 has simple head $\widetilde{S}_c(L_{c}(\mu))$ for each $\mu$.
By \eqref{nonzeromap} $\widetilde{S}_c(L_{c}(\lambda))$
 is therefore  isomorphic to a composition
factor of $\Delta_{c+1}(\lambda)$. But, by Corollary~\ref{poono}
and the remark thereafter,  the composition factors of
$\Delta_{c+1}(\lambda)$, except for  the head, are of the form $L_{c+1}(\nu)$ for
$\nu < \lambda$. Thus the non-zero map
\eqref{nonzeromap} must send the head of
$\widetilde{S}_c(\Delta_{c}(\lambda))$ to the head
 of $\Delta_{c+1}(\lambda)$ and so induce  an
isomorphism  $\widetilde{S}_c(\Delta_{c}(\lambda))\xrightarrow{\sim}
\Delta_{c+1}(\lambda)$. This completes the inductive step, and hence
the proof of  the proposition.
\end{proof}


\section{The Hilbert scheme}\label{sect-haiman}

\subsection{}\label{sect401}
Haiman's  work on Hilbert schemes
gives detailed information about their structure, in particular as ``Proj'' of
appropriate Rees rings.
The resulting formul\ae\ for the Poincar\'e series of
these rings will be crucial to the proof of the main theorem in
Section~\ref{sect-filt}. In this section, we briefly describe the
relevant results from the literature and  use this
 to derive the appropriate Poincar\'e series.

\subsection{}\label{subsec-5.1}
Let $\hin$\label{hin-defn} be {\it the Hilbert scheme of $n$
points on the plane}, which we realise as the set of ideals of
colength $n$ in the polynomial ring $\C[\C^2]$. Similarly,
identify the variety $S^n\C^2$ of \textit{unordered} $n$-tuples of
points in $\C^2$ with the categorical quotient $ \C^{2n}/\WW$
under the diagonal action of $\WW$ on $\C[\C^{2n}]$.
 Then the map\label{tau-defn}  $\tau : \hin \to S^n\C^2=\C^{2n}/\WW$
which sends an ideal to its support (counted with
multiplicity) is a resolution of singularities (see, for example,
\cite[Theorem~1.15]{Nak}).

 We will actually be interested in a resolution of singularities for
 $\C[\h\oplus \h^*]^\WW$ rather
than $\C[\C^{2n}]^\WW$, simply  because the
associated graded ring of $U_c$ is $\C[\h\oplus \h^*]^\WW$.
The results we need   follow easily from
the corresponding results on $\hin$ and  so we begin with the latter.

\subsection{The (isospectral) Hilbert scheme} \label{isospecsec} Following
 \cite[Definition~3.2.4]{hai3} {\it  the
 isospectral Hilbert
scheme}  $\XXX_n$\label{isospec-defn} is the reduced fibre product
$$\begin{CD} \XXX_n @> f_1  >> \C^{2n}
\\ @V \rho_1 VV @VVV \\ \hin @> \tau >> \C^{2n}/\WW. \end{CD}
$$
It is a highly
non-trivial fact (see \cite[Theorem~3.1 and the proof of
Proposition~3.7.4]{hai3}) that $\rho_1$ is a flat map of degree $n!$.

Haiman has given a description of both $\hin$ and
$\XXX_n$ as Proj of appropriate graded rings and we recall this description
 since
it will be extremely important to us.
Let $ \AAA^1=\C[\C^{2n}]^{\ep}$ \label{AAA-1-defn}
 be the space of $\WW$-alternating polynomials  in $\C[\C^{2n}]$ and write
  $  \JJJ^1= \C[\C^{2n}]\AAA^1$
for the ideal generated by    $\AAA^1$.
For $d\geq 1$ define $\AAA^d$ and $\JJJ^d$ to be the
respective  $d^{\text{th}}$ powers of $\AAA^1$ and $\JJJ^1$ using
 multiplication in $\C[\C^{2n}]$; thus
 \label{JJJ-defn}  $\JJJ^d=\C[\C^{2n}] \AAA^d$.
 Finally, set $\JJJ^0=\C[\C^{2n}]$,
 $\AAA^0=\C[\C^{2n}]^\WW$ and
 $\AAA=\bigoplus_{d\geq 0} \AAA^d\cong \AAA^0[t\AAA^1]$.
 Then   \cite[Proposition~2.6]{haidis} proves that
\begin{equation}\label{AAA-alg-defn}
\hin \cong \prj \AAA \mathrm{ \ as\ a\ scheme\ over\ }
\spc \AAA^0 =
\C^{2n}/\WW\end{equation}
  Similarly,    $\XXX_n \cong \prj  \dd$, where $\dd=\C[\C^{2n}][t\JJJ^1],$
is the blowup of
$\C^{2n}$ at $\JJJ^1$ \cite[Proposition~3.4.2]{hai3}.

\subsection{}\label{subsec-5.5} Observe that $\JJJ^d$ is generated by its
$\WW$-alternating or
$\WW$-invariant elements, respectively, depending on whether $d$
is odd or
even. Following Haiman we refer to these elements as having \textit{correct
parity}.

\begin{lem} \label{corpar} {\rm (1)} For any $d\geq 0$, $\AAA^d$ consists
of the elements of $\JJJ^d$ with the correct parity.

{\rm (2)} If $\C^n$ denotes the first copy of that space in $\C^{2n}$,
 then $\JJJ^d$ is a free module over both  $\C[\C^n]$ and $\C[\C^n]^{\WW}$.
 \end{lem}

\begin{proof} (1)
The statement is clearly true for $d=0,1$. Assume, by induction, that it is true
for $d-1$. We will suppose that $d$ is even, the argument in the odd case
 being similar. Since $\AAA^1$ generates the ideal $\JJJ^1$, any element
   $x\in \JJJ^d$ can be decomposed as $x= \sum_i p_i q_i$ where
   $p_i \in \JJJ^{d-1}$
  and $q_i \in \AAA^1$. Since $q_ie = e_-q_i$ we
  have $(p_iq_i)e = (p_i e_-)q_i$ for all $i$. If $x$ has the correct parity
  then $x = xe = \sum_i (p_iq_i)e = \sum_i (p_ie_-q_i).$
  But  $\JJJ^{d-1}e_{-}$ is the subset of
   $\WW$-alternating elements of $\JJJ^{d-1}$ and so
   $\JJJ^{d-1}e_{-}=\AAA^{d-1}$ by induction. Thus $x\in \AAA^{d-1}  \AAA^1
   =\AAA^d$.

   (2) By \cite[Proposition~4.11.1]{hai3} $\JJJ^d$ is a  projective  module
   over $\C[\C^n]$ and hence over $\C[\C^n]^\WW$.    Since $\C[\C^n]$ and
   $\C[\C^n]^\WW$  are polynomial rings, any such
   projective module is free by the Quillen-Suslin Theorem.
\end{proof}

\subsection{Geometric interpretation} \label{geomint}  There is a geometric
description of both $\AAA^d$ and $\JJJ^d$. Let $\BB_1$\label{tauto-defn} be the
{\it tautological rank $n$ vector bundle} on $\hin$  and  let
$\PP_1 = (\rho_1)_{\ast}\OO_{\XXX_n}$ \label{PP1-defn} denote the
{\it Procesi   bundle} of rank $n!$ arising from the map
$\rho_1 : \XXX_n \to \hin$. Write
$\LL_1 = \bigwedge^n \BB_1$\label{LLL-defn} for
the determinant bundle of $\BB_1$.  By \cite[Proposition~2.12]{haidis}
 $\LL_1$ is also the
canonical ample line bundle $\OO_{\hin}(1)$ associated to the presentation
$\hin\cong \prj \AAA$.

\subsection{}\label{subsec-5.7}
Set $l=dn$ for some $d\geq 1$  and write
 \label{RR-defn} $\RRR(n, l)
= H^0(\hin, \PP_1\otimes \BB_1^{l}).$
One should note that $\RRR(n,l)$ is defined in \cite{hai1}
to be the coordinate ring of the
polygraph $Z(n,l)$ but, by \cite[Theorem~2.1]{hai1},  it is also
isomorphic to the given ring of global sections. There is an action of
$\WW\times \WW^d$ on
$\PP_1\otimes \BB_1^{l}$,  with $\WW$ acting fibrewise on $\PP_1$ and
$\WW^d\subset \mathfrak{S}_{l}$
acting on $\BB_1^{l}$ by permutations.  By construction,
$(\PP_1)^{\WW} =\OO_{\hin}$  and  $(\BB_1^l)^{\ep_d} = \LL_1^d$,
where $\ep_d$\label{epsilon-defn} denotes the sign representation
of $\WW^d$.

The proof of \cite[Proposition~4.11.1]{hai3} shows that
 $\JJJ^d \cong \RRR(n, l)^{\ep_d}$.
 On the other hand, the action of $\WW^d$ is trivial on $\hin$, so
$\PP_1\otimes \BB_1^{l}$ is a
 direct sum of its isotypic components. Hence
\begin{equation} \label{cohj}
 \JJJ^d \cong \RRR(n, l)^{\ep_d} = H^0(\hin, (\PP_1\otimes \LL_1)^{\ep_d})
 =H^0(\hin, \PP_1\otimes \LL_1^d).
 \end{equation}
  It is not true, however, that the natural $\WW$-action on the two  sides
  agrees. Indeed, thanks to the proof of \cite[Proposition~4.2]{hai2}
  the isomorphism written $\WW$-equivariantly is
   \begin{equation}\label{cohj1}
   \JJJ^d \otimes \ep^{\otimes d} \cong H^0(\hin, \PP_1\otimes \LL_1^d).
   \end{equation}
   The reason for this is that the isomorphism in \eqref{cohj}
    is given by the $\C[\C^{2n}]$-module
    homomorphism sending (in the notation of \cite{hai2})
    the generators
     $\Delta_{t_1}(\mathbf{a},\mathbf{b})\cdots
     \Delta_{t_d}(\mathbf{a},\mathbf{b})$  on the right hand side
 to their evaluations on the left hand side: $\mathbf{a}\mapsto \mathbf{x}$,
     $\mathbf{b}\mapsto \mathbf{y}$. The
   element $\Delta_{t_j}(\mathbf{a},\mathbf{b})$ has a
      trivial $\WW$-action as no $\mathbf{x}$'s or
     $\mathbf{y}$'s are involved, whereas its specialisation
     has a $\WW$-action of $\ep^{\otimes d}$
     since that specialisation is the product
     of $d$ determinants.

As a result, \eqref{cohj1}  and Lemma~\ref{corpar} combine to prove:

 \begin{lem}\label{coha} There is an isomorphism of $\AAA^0$-modules
 $\AAA^d \cong \RRR(n,l)^{\WW\times \ep_d} = H^0 (\hin, \LL_1^d)
 .$\qed
 \end{lem}

\subsection{(Bi)graded characters}
\label{bigr1}
There is a bigrading\label{bigrad-defn} on
$\C[\C^{2n}]=\C[\mathbf{x},\mathbf{y}]$
 with $\deg x_i = (1,0)$ and
 $\deg y_j = (0,1)$ which, as in \cite[(12)]{hai1},  arises from the action of
  $\TT^2= (\C^*)^2$ on the plane $\C^2$ given by $(\alpha, \beta )\cdot (u,v)
  = (\alpha^{-1} u, \beta^{-1} v)$ for $(u,v)\in \C^2$.
   This action extends  to $\hin$, and the
   bundles $\PP_1$, $\BB_1$, $\LL_1$ are naturally $\TT^2$-equivariant.
   The isomorphisms from  \eqref{cohj} and Lemma~\ref{coha} respect the induced
    bigradings.
Of course, the sections $M$  of any one
 of these modules obtains an induced action of
$\TT^2$ and this is equivalent to a $\Z^2$-grading $M=\bigoplus M_{i,j}$;
explicitly,
 an element $f\in M$ is homogeneous of weight $(i,j)$ if
 $(\alpha,\beta)f=\alpha^i\beta^jf$.

The $\TT^2$--fixed points of $\hin$ are precisely the ideals $I_\mu$ that are
associated to partitions $\mu$ of $n$ by \label{I-mu-defn}
 $$I_{\mu} = \C \cdot \{ x^ry^s : (r,s)\notin
d(\mu)\} \subseteq \C [x,y],$$ see \cite[Proposition~3.1]{hai1}.
The set of monomials $\mathcal{B}_{\mu} = \{ x^ry^s : (r,s)\in
d(\mu) \}$ that are not in $I_{\mu}$ form a natural $\C$-basis of
$\C[x,y]/I_{\mu}$.

 \subsection{}\label{subsec-5.10}
For a bigraded space $V = \sum_{i,j}V_{i,j}$ with finite
dimensional weight spaces we define the {\it bigraded  Poincar\'{e} series}
\label{Poincare-defn} of $V$ to be
$$p(V,s,t) = \sum_{i,j} \dim(V_{i,j})s^it^j.$$

Haiman has calculated the bigraded Poincar\'{e} series of
$\RRR(n,l)$ and a similar calculation will allow us to find the
bigraded Poincar\'{e} series of $\JJJ^d$.  For a pair of
partitions $\lambda, \mu$ let $K_{\lambda \mu}(t,s)$
\label{kostka-defn} be the {\it Kostka--Macdonald coefficients}
defined in \cite[VI, (8.11)]{MacD}. Set
$$\Omega(\mu)=
\prod_{x\in d(\mu)}(1-s^{1+l(x)}t^{-a(x)})(1-s^{-l(x)}t^{1+a(x)})
\qquad\mathrm{and}\qquad P_{\mu}(s,t) = \sum_{\lambda}
s^{n(\mu)}K_{\lambda \mu}(t,s^{-1}) f_{\lambda}(1).$$ We remark
that many of the formul\ae\ we cite from Haiman's papers are given
in terms of Frobenius series $\mathcal{F}_M(z;s,t)$ but, as in
\cite[(6.5)]{hai2}, we can always specialise these to Hilbert
series $p(M,s,t)$ by specialising $s_\lambda(z)$ to
$f_\lambda(1)=\dim \lambda$.

\begin{prop} \label{bigr-hain} Under the $\TT^2$-bigraded structure,
the bigraded Poincar\'{e} series of $\JJJ^d$ is
$$p(\JJJ^d, s, t) = \sum_{\mu} P_{\mu}(s,t) \, \Omega(\mu)^{-1}
 s^{dn(\mu)}t^{dn(\mu^t)}.$$
\end{prop}
\begin{proof}
By \cite[Theorem~2.1]{hai1}
$H^i(\hin , \PP_1 \otimes \LL_1^d) = 0$ for $i>0$, while
$H^0(\hin , \PP_1 \otimes \LL_1^d) = \JJJ^d$ by \eqref{cohj}. Thus,
 in the notation of \cite[Section~3]{hai1},
$p(\JJJ^d, s, t)=\chi_{\PP_1 \otimes \LL_1^d}(s,t)$ and so, by
\cite[Proposition~3.2]{hai1},
\begin{equation}\label{bigr11}
p(\JJJ^d, s, t)=\sum_\mu p(\PP_1\otimes \LL_1^d(I_\mu),s,t)\,\,
\Omega(\mu)^{-1}
=\sum_\mu p(\PP_1(I_{\mu}),s,t) p(\LL_1(I_\mu),s,t)^d\,\,
\Omega(\mu)^{-1}.\end{equation}
Here we have used the fact that, as
$I_\mu$ defines a finite dimensional scheme,
 we can identify the sheaf $\PP_1\otimes \LL_1^d(I_\mu)$ with its global
 sections, and so  $p(\PP_1\otimes \LL_1^d(I_\mu),s,t)$
 is naturally defined.

  We now evaluate the right hand side of \eqref{bigr11}.
     It is proved in  \cite[(3.9)]{haidis},
   using the notation of  \cite[(1.9)]{haidis}, that
   $p(\LL_1(I_{\mu}),s, t) = \prod_{x\in d(\mu)} s^{l'(x)}t^{a'(x)}=
    s^{n(\mu)}t^{n(\mu^t)}$.
 On the other hand,  by \cite[Proposition~3.4]{hai1}
(which is proved in \cite[Section~3.9]{hai3} and uses the notation of
\cite[(46)]{hai1}), $p(\PP_1(I_{\mu}),s, t)
    = P_{\mu}(s,t)$.
    Substituting these observations into \eqref{bigr11} shows
    that $$ p(\JJJ^d, s, t)=\sum_\mu P_{\mu}(s,t)\, \Omega(\mu)^{-1}
  s^{dn(\mu)}t^{dn(\mu^t)},$$ as required.
   \end{proof}

\subsection{Blowing up $(\h\oplus \h^*)/\WW$}\label{hi-defn-sec}
All the results described so far have natural analogues for the
subvariety   $\h\oplus \h^*$ of $\C^{2n}$. Geometrically, this follows from
 the observation that the natural additive action of
$\C^2$ by translation  on $\hin$ gives a decomposition
 $\hin = \C^2\times \left( \hin\right)/\C^2$ into a product of varieties
 \cite[p.10]{Nak}. Unravelling the actions shows that $ \hin/\C^2$
 provides a resolution of singularities for  $\h\oplus \h^*$.
 However, since we need the algebraic consequences of Haiman's results, we will
 take a more algebraic approach.

   We emphasise that the embedding
$\h\oplus \h^* \hookrightarrow \C^{2n}$ is always given by  embedding $\h$
into the first
copy of $\C^n$ and $\h^*$ into the second copy. To fix notation, let $\h$ be
the hypersurface $\mathbf{z}=0$ in  $\C^{n}$ and similarly let $\h^*$ be the
hypersurface $\mathbf{z}^*=0$ in the  second copy of $\C^{n}$; thus
$\C[\C^{2n}]=\C[\h\oplus\h^*][\mathbf{z},\mathbf{z}^*]$. Since
$\mathbf{z},\mathbf{z}^*\in \C[\C^{2n}]^\WW$,
this induces the decomposition
$\C[\C^{2n}]^\WW =\C[\h\oplus\h^*]^\WW [\mathbf{z},\mathbf{z}^*]$. Following the
lead of \eqref{isospecsec}, we set
 $$A^1
= \cxy^{\ep}\ \subset \ \AAA^1=\C[\C^{2n}]^{\ep}\qquad\text{and}\qquad
J^1 = \cxy A^1\ \subset \ \JJJ^1= \C[\C^{2n}]\AAA^1.$$ \label{A-1-defn}
We then define $A^0=\cxy^\WW$, $J^0=\cxy$ and, for $d>1$,  take
$A^d=(A^1)^d$ and $J^d = (J^1)^d$ for  the respective $d^{\text{th}}$ powers
  using the  multiplication in $\cxy$. Finally, we write
 $$A=\bigoplus_{i\geq 0}A^i
 \cong A^0[A^1t]\qquad\text{and}\qquad S = \bigoplus_{i\geq 0}J^i
 \cong \cxy[J^1t]$$ for the corresponding Rees rings.
The next result is basic observation about these objects.

\begin{lem}\label{hi-basic-lem}  {\rm (1)}
For $d\geq 0$, $\AAA^d = A^d[\mathbf{z},\mathbf{z}^*]$
is the set of  polynomials with coefficients from $A^d$. Similarly,
 $\JJJ^d = J^d[\mathbf{z},\mathbf{z}^*]$.

{\rm (2)}
 Each $J^d$ is a free module over $\C[\h]$ and $\cx^{\WW}$.
\end{lem}

 \begin{proof}   (1)
 By definition,
 $\AAA^1=\left(\cxy[\mathbf{z},\mathbf{z}^*]\right)^\ep =
   \C[\h\oplus\h^*]^{\ep} [\mathbf{z},\mathbf{z}^*]
=A^1[\mathbf{z},\mathbf{z}^*]$ as polynomial extensions.
Thus $\AAA^d = (A^1[\mathbf{z},\mathbf{z}^*])^d = (A^1)^d[\mathbf{z},\mathbf{z}^*]
= A^d[\mathbf{z},\mathbf{z}^*]$
and  $\JJJ^d = A^d[\mathbf{z},\mathbf{z}^*]\cxy
= J^d[\mathbf{z},\mathbf{z}^*].$

(2)
By part (1) and Lemma~\ref{corpar},
$\JJJ^d = J^d[\mathbf{z},\mathbf{z}^*]$ is a free module over
$\cx[\mathbf{z}]$ and hence over $\cx$. Therefore, so is its
$\cx$-module summand $J^d$. \end{proof}

\subsection{}\label{hi-defn-sec2}
Recall the resolution of singularities  $\tau: \hin \to
\C^{2n}/\WW$  defined in \eqref{hin-defn} and
define $\hi=\tau^{-1} (\h\oplus\h^*/\WW)$ \label{hi-defn}, with the resulting
morphism $\tau: \hi\to \h\oplus\h^*/\WW$.  Using the identifications
 of \eqref{hi-defn-sec},
the basic properties of $\hi$ are easy to determine.

\begin{cor}\label{hi-basic-lem2} {\rm (1)} $\hi =\prj(A) $ and
 $\tau :\hi \to \h\oplus \h^*/\WW$
is  a resolution of  singularities.
\begin{enumerate}
\item[(2)] Moreover $\tau$ is a crepant resolution: that is
$\omega_{\hi} \cong \OO_{\hi}$.
 \item[(3)] Set $X_n=\prj(S)$. Then
$X_n$ is the reduced fibre product
 $$\begin{CD} X_n @> >> \h\oplus \h^*
\\ @V \rho VV @VVV \\ \hi @> \tau >> \h\oplus \h^*/\WW. \end{CD}
$$ and the map $\rho$ is flat of degree $n!$.
\end{enumerate}
\end{cor}

\begin{proof} (1)
 Recall from \eqref{AAA-alg-defn} that $\hin = \prj(\AAA)$.
By Lemma~\ref{hi-basic-lem}, $\AAA= A[\mathbf{z},\mathbf{z}^*]$.
The maps $A\hookrightarrow \AAA$ and
$\C[\mathbf{z},\mathbf{z}^*] \hookrightarrow \AAA$ give maps
$\hin\to \prj(A)$ and $\hin\to\mathrm{Spec}(\C[\mathbf{z},\mathbf{z}^*])
\cong  \C^2$ and hence, by universality,
a map    $\hin \to  \prj(A)\times \C^2$.
It is easy to check that this is an isomorphism locally and hence globally.
 The identification
 of $\h\oplus\h^*$ with the subvariety $\mathbf{z}=0=\mathbf{z}^*$
 of $\C^{2n}$ easily yields $\hi=\prj(A)$ and  so   $\hin = \hi\times \C^2$.
Since $\hin$ is a resolution of singularities of
$\C^2/\WW$, the result follows.

(2) By \cite[Exercise II.8.3(b)]{har} $\omega_{\hin} \cong
\omega_{\hi} \boxtimes \omega_{\C^2}$, the external tensor product
on $\hin = \hi \times \C^2$. Now (2) follows since $\omega_{\hin}
\cong \OO_{\hin}$ by \cite[Proposition 3.6.3]{hai3}.

(3) As in part (1), $\mathbb{S}= \bigoplus \JJJ^d
=S[\mathbf{z},\mathbf{z}^*]$ and $\prj(\mathbb{S}) \cong
\prj(S)\times \C^2$. The assertions of the corollary now follow
from the corresponding results for $\XXX=\prj(\mathbb{S})$ that
were stated in \eqref{isospecsec}.
\end{proof}

We also have analogues of $\PP_1$ and $\LL_1$ for $\hi$. These are defined in
the same way:
$\PP = \rho_{\ast}\OO_{X_n}$\label{PP-defn} is the
{\it Procesi bundle} on $\hi$  of rank $n!$ arising from the map $\rho : X_n \to \hi$
while  $\LL$\label{LL-A-defn} is  the
canonical ample line bundle $\OO_{\hi}(1)$ associated to the presentation $\hi
\cong \prj A$.

\subsection{}   Since
${\bf z}$ and ${\bf z}^*$ are bihomogeneous,
the bigradings of \eqref{bigrad-defn} to  pass $\hi$. Therefore,
  Lemma~\ref{hi-basic-lem}(1)  implies that
 $p(J^m, s, t)=(1-s)(1-t) p(\JJJ^m, s, t)$. Substituting this formula into
Proposition~\ref{bigr-hain} gives:
\begin{cor} \label{bigr} The bigraded Poincar\'{e} series of $J^d$ is
$$\qquad\qquad \qquad\qquad\qquad\qquad
p(J^d, s, t) = \sum_{\mu} P_{\mu}(s,t)(1-s)(1-t)\,
\Omega(\mu)^{-1} s^{dn(\mu)}t^{dn(\mu^t)}.
 \hfill \qquad\qquad\qquad \qquad\hfill\qed $$
\end{cor}

\subsection{}\label{babyverma} In Corollary~\ref{gr} we will give a singly
graded analogue of Corollary~\ref{bigr} that will be needed in the proof of the
Theorem~\ref{mainthm-intro}.
 In the proof we will  need the following combinatorial formul\ae\
for the fake degrees $f_\mu(v)$, as  defined in \eqref{fakedegrees}.

\begin{lem}  Let $\mu\in\irr{\WW}$. Then

{\rm (1)}
$ f_{\mu}(v) = v^N f_{\mu^t}(v^{-1})$, where $N=n(n-1)/2$,

{\rm (2)} $f_{\mu}(v)\prod_{x\in d(\mu)} (1 - v^{h(x)})
 = v^{n(\mu)}\prod_{i=1}^n (1-v^i),$ where $h(x)=1+a(x)+l(x)$ as in
 \eqref{d-mu-defn},

 {\rm (3)} $\sum_{\lambda} v^{n(\mu)}K_{\lambda \mu}(v^{-1},v^{-1})f_{\mu}(v^{-1})
 f_{\lambda}(1)
  =  \sum_{\lambda} f_{\lambda}(v^{-1}) f_{\mu}(1)f_{\lambda}(1).$
 \end{lem}

 \begin{proof} (1) This is a well-known formula (see, for example,
 \cite[p.453]{op}).

 (2,3) Up to a change of notation, these are both proved within the
 proof of  \cite[Theorem~6.4]{babyv}---see the displayed equations immediately
 after, respectively immediately before \cite[(18)]{babyv}.
 \end{proof}

\subsection{} \label{gr} The $\EE$-grading from \eqref{Euler-defn}
descends naturally to $\ogr D(\h)\cong \C[\h\oplus\h^*]$ and we will use
the same notation there;  thus $\Edeg \h^*=1$ and $\Edeg \h = -1$.
  For   an $\EE$-graded module (or, indeed, any $\Z$-graded module)
  $M=\bigoplus_{i\in\Z}M_i$, we write the corresponding Poincar\'e series as
  $p(M,v)=\sum v^i\dim_\C M_i$. Set
 \begin{equation} \label{factorial-defn}
 [n]_v! = \frac{\prod_{i=1}^n (1-v^i)}{(1-v)^n} . \end{equation}

\begin{cor} Under the $\EE$-grading,  the module \
$\overline{J^d} = {J^d}/{\cx^{\WW}_+J^d}$ \ has  Poincar\'e series

\begin{equation}
\label{Z-grading0}
 p(\overline{J^d}, v) =  \frac{\sum_{\mu}
f_{\mu}(1)f_{\mu}(v^{-1}) v^{-d(n(\mu) - n(\mu^t))}[n]_v!}{\prod_{i=2}^n
(1-v^{-i})}  . \end{equation}
\end{cor}

\vskip 10pt

\begin{proof} Since $\C[\h]^\WW_+$ is $\EE$-graded, so is $\overline{J^d}$,
 and so the result does make sense.  By Lemma~\ref{corpar}(2),
  the fundamental invariants of $\C[\h]^\WW$
form an r-sequence
in $J^d$ for any $d\geq 0$. Since these elements have degrees $2\leq r\leq n$,
 Corollary~\ref{bigr} implies that
  $\overline{J^d}$  has Poincar\'e series
\begin{equation}
\label{Z-grading15}
p(\overline{J^d}, v) =  \left( (1-t)\prod_{i=1}^n(1- s^i)\sum_{\mu}
 P_{\mu}(s,t)\, \Omega(\mu)^{-1}  s^{dn(\mu)}t^{dn(\mu^t)}
 \right)_{s=v, t=v^{-1}}
\end{equation} where $P_\mu$ and $\Omega(\mu)$ are defined in
\eqref{subsec-5.10}. Lemma~\ref{babyverma}(2) implies that
$$\label{Z-grading16}
 \Big(\Omega(\mu)\Big)_{s=v,t=v^{-1}}
 =f_{\mu}(v)^{-1}f_{\mu}(v^{-1})^{-1}
 \prod_{i=1}^n(1-v^i)(1-v^{-i}).$$
 This  gives
  \begin{equation}\label{Z-grading2}
   p(\overline{J^d}, v) =
 \frac{  \sum_{\mu} P_{\mu}(v,v^{-1})f_{\mu}(v)f_{\mu}(v^{-1})
  v^{dn(\mu)}v^{-dn(\mu^t)} }
  {\prod_{i=2}^n (1-v^{-i})}.
  \end{equation}
By Lemma~\ref{babyverma}(3) the numerator of this expression can
be described as
\begin{equation}
\label{numerator}
 \sum_{\mu} \left(\sum_{\lambda}
  f_{\lambda}(v^{-1}) f_{\lambda}(1)\right) f_{\mu}(1) f_{\mu}(v)
  v^{d(n(\mu)-n(\mu^t))}.
  \end{equation}
  Applying Lemma~\ref{babyverma}(1) and using  the equality
  $f_{\mu}(1) = f_{\mu^t}(1)$ from \eqref{fakedegrees2}
we find that   \eqref{numerator} equals
 \begin{equation} \label{rearrange} \sum_{\mu} \left(\sum_{\lambda}f_{\lambda}(v^{-1})
  f_{\lambda}(1)\right) f_{\mu^t}(1) f_{\mu^t}(v^{-1}) v^N
  v^{-d(n(\mu^t) - n(\mu))}.\end{equation}
The standard formula $\sum \dim
\C[\h]^{\text{co}\WW}_iv^{-i}=[n]_{v^{-1}}!$ shows that the fake
degrees satisfy the  identity
$$\sum_{\lambda} f_{\lambda}(v^{-1}) f_{\lambda}(1)=
   \frac{\prod_{i=1}^n (1-v^{-i})}{(1-v^{-1})^n} = [n]_{v^{-1}}!.$$
Applying this and    \eqref{rearrange}
 to  \eqref{Z-grading2} we find that
\begin{equation}
\label{weredone} p(\overline{J^d}, v) = \frac{\sum_{\mu}
f_{\mu^t}(1) f_{\mu^t}(v^{-1}) v^{-d(n(\mu^t) - n(\mu))}
v^N[n]_{v^{-1}}!}{\prod_{i=2}^n (1-v^{-i})}. \end{equation} After
changing the order of summation from $\mu$ to $\mu^t$ and using
the equality $$ v^N [n]_{v^{-1}} = v^N\frac{\prod_{i=1}^n
(1-v^{-i})}{(1-v^{-1})^n} = \frac{\prod_{i=1}^n (1-v^i)}{(1-v)^n}
= [n]_v!,$$ \eqref{weredone} becomes the required equality
\eqref{Z-grading0}, and so the corollary is proved.
\end{proof}


 \section{$\Z$--algebras} \label{zalg}

\subsection{}\label{Z-alg-defn}   Typically in noncommutative algebra---and
certainly in our case---one  cannot apply the Rees ring construction since
one is working with just right modules or homomorphism groups rather than
bimodules. One way round this is  to use $\Z$-algebras and in this section we
describe the basic properties that we need from this theory. The reader is
referred to \cite{BP} or \cite[Section~11]{SV} for the more general theory
and to \cite[Section~3]{bgs} for applications of $\Z$-algebras
 to Koszul duality.

 Throughout this paper a {\it $\Z$-algebra}
 will mean a {\it lower triangular $\Z$-algebra}. By definition, this is a
 (non-unital)  algebra $B=\bigoplus_{i\geq j\geq 0} B_{ij}$, where
 multiplication is defined in matrix fashion:
 $B_{ij}B_{jk}\subseteq B_{ik}$ for $i\geq j\geq k\geq 0$
 but $B_{ij}B_{\ell k}=0$ if $j\not=\ell$.
 Although $B$ cannot have a unit element, we do require that each
 subalgebra $B_{ii}$  has a unit element $1_i$ such that
 $1_ib_{ij}=b_{ij}=b_{ij}1_j$, for all  $b_{ij}\in B_{ij}$.

\subsection{}\label{Z-alg-defn2}
Let $B$ be a $\Z$-algebra. We define the category $B\lGr $   to be the category
of $\NN$-graded left  $B$-modules $M=\bigoplus_{i\in \NN} M_i$ such that
$B_{ij}M_j\subseteq M_i$ for all $i\geq j$  and $B_{ij}M_k=0$ if $k\not=j$.
Homomorphisms are defined to be graded homomorphisms of degree zero. The
subcategory of noetherian graded left $B$-modules will be denoted $B\lgr $. In
all  examples considered in this paper $B\lgr $ will consist precisely of the
finitely generated graded left  $B$-modules.

A module $M\in B\lGr $ is {\it bounded} if $M_n = 0$ for all but finitely many
$n\in \Z$ and {\it torsion} if it is a direct limit of bounded modules. We
let $B\lTors$ denote the full subcategory  of torsion modules in $B\lGr$
   and write  $B\ltors $
for the analogous subcategory of $B\lqgr$. The corresponding quotient
categories are written $B\lQgr = B\lGr /B\lTors $ and $B\lqgr  = B\lgr /B\ltors
$.\label{qgr-defn} Write $\pi(M)$ for the image in $B\lQgr $ of $M\in
B\lGr$.

\subsection{} \label{zalgex1} There are two basic examples
of $\Z$-algebras that will interest us. For the first,
suppose that $S= \bigoplus_{n\geq 0} S_n$ is an
$\NN$-graded algebra. As in \cite[Example~3.1.3]{bgs}
we  can canonically associate a $\Z$-algebra
   $\widehat{S}=\bigoplus_{i\geq j\geq 0}\widehat{S}_{ij}$ to $S$
   by setting $\widehat{S}_{ij} = S_{i-j}$ with multiplication
induced from  that in $S$. Define  categories
$S\lGr,\dots,S\lqgr$  in the usual manner. In particular,
 $S\lGr$ denotes the category
of $\Z$-graded $S$-modules, from which the other
definitions follow as in the last paragraph. We then let
$S\lGr_{\geq 0}$ denote the full subcategory of
$S\lGr$ consisting of $\NN$-graded $S$-modules
$M=\bigoplus_{i\in\NN}M_i$.
It is immediate from the definitions that
the identity map
$\iota: M=\bigoplus_{i\in\NN}M_i\mapsto M=\bigoplus_{i\in\NN}M_i$
gives equivalences of categories
$S\lGr_{\geq 0}\simeq \widehat{S}\lGr$ and
$S\lgr_{\geq 0}\simeq\widehat{S}\lgr$. For any
module $M\in S\lGr$, one has $\pi(M)=\pi(M_{\geq 0})$
in $S\lQgr$ and so $\iota$  induces category equivalences
\begin{equation}\label{zalgex11}
S\lQgr\simeq \widehat{S}\lQgr\qquad \text{and} \qquad
S\lqgr\simeq \widehat{S}\lqgr.
\end{equation}

\subsection{} \label{zalgex2} For the second class of examples,
suppose that we are given noetherian algebras
$R_n$ for $n\in \NN$ with $(R_i, R_j)$-bimodules
$R_{ij}$, for $i> j\geq 0$. Assume, moreover, that
there are   morphisms
$\theta_{ij}^{jk} : R_{ij}\otimes_{R_j}R_{jk}\to R_{ik}$ satisfying the
the obvious associativity conditions.
Then we can define a $\Z$-algebra
$R_\Z$ by $R_\Z=\bigoplus_{i\geq j\geq 0}R_{ij}$, where
 $R_{ii}=R_i$ for all $i$.

A particular example of this construction is the one that interests us.
Suppose that  $\{R_n : n\in \NN\}$ are Morita equivalent algebras,
with the equivalence induced from the progenerative $(R_{n+1},R_{n})$-bimodules
$P_n$. Define $R_{ij} =
P_{i-1}\otimes_{R_{i-1}}\otimes\cdots \otimes_{R_{j+2}} P_{j+1}\otimes_{R_{j+1}}P_j$
and $R_{jj}=R_j$, for $i>j\geq 0$.
Tensor products provide the isomorphisms $\theta_\bullet^\bullet$
and   associativity is  automatic.
The corresponding $\Z$-algebra $R_\Z=\bigoplus_{i\geq j\geq 0} R_{ij}$
will be called the
{\it Morita $\Z$-algebra associated to the data $\{R_n,P_n : n\in \NN\}$}.

\subsection{}
\label{Zalgequiv}
 Write $R\lmod$ for the category of finitely generated
left modules over a noetherian ring $R$.
Although easy, the next result provides the foundation
 for our approach to $U_c$: in order to study $U_c\lmod$ it suffices
 to study  $R_\Z\lqgr$, for any Morita
$\Z$-algebra $R_\Z$  with $R_0\cong U_c$.

\begin{lem} Suppose that $R_\Z$ is the Morita $\Z$-algebra
associated to the data $\{R_n,P_n : n\in \NN\}$, where $R_0$ is noetherian.
\begin{enumerate}
  \item[(1)] Each finitely generated graded left $R_\Z$-module is noetherian.
  \item[(2)] The
 association $\phi: M \mapsto \bigoplus_{n\in \NN}R_{n0}\otimes_{R_0}M$
 induces an equivalence of categories between
 $R_0\lmod$ and $R_\Z\lqgr$.
 \end{enumerate}
\end{lem}

\begin{proof} (1)
 Any finitely generated graded left $R_\Z$-module $M$ is a
graded image of
$\bigoplus_{a_i} \big(\bigoplus_{j\geq a_i} R_{ja_i}\big)
\otimes_{R_{a_i}} R_{a_i},$ for some $a_i\in \NN$
and so we may assume that
$M=  \bigoplus_{j\geq a} R_{ja},$ for some $a\geq 0$.
Let $L\subseteq M$ be a graded submodule  and write $R^*_{ij}$ for the dual of
the progenerator $R_{ij}$.  Then
$$X(j)=R^*_{ja}\otimes_{R_j}L_j \subseteq
R^*_{ja}\otimes_{R_j}M_j =
R^*_{ja}\otimes R_{ja}
\xrightarrow{\sim} R_a,\qquad \text{for}\ j\geq a.$$
As $R_a$ is Morita equivalent to $R_0$, it is noetherian and so $\sum_{j\geq
a}X(j)=\sum_{i=a}^b X(i)$, for some $b\geq a$. Now,
$$L_k = R_{ka}X(k) \subseteq \sum_{i=a}^b R_{ka}X(i) =
\sum_{i=a}^b R_{ki}R_{ia}X(i) =\sum_{i=a}^b R_{ki}L_i
\qquad\text{for}\ k\geq a.$$
Thus $L$ is generated by $L_j$ for $b\geq j\geq a$.
Finally, as each $L_i$ is a submodule of the noetherian left
$R_i$-module $R_{ia}$,
it is finitely generated  and hence so is $L$.

 (2) Certainly
 $\phi(M)\in R_\Z\lGr $ and, as $\phi(M)$  is finitely generated by the generators
 of ${}_{R_0}M$, one has $\phi(M)\in R_\Z\lgr $.
Thus $\Phi (M)=\pi\phi(M)\in R_{\Z}\lqgr$. Since $\Phi$ sends
$R_0$-module homomorphisms to graded $R_\Z$-module homomorphisms, $\Phi$ is a
functor.

Conversely, suppose that $\widetilde{N}\in R_\Z\lqgr$ and pick a preimage $N\in
R_\Z\lgr$. Then $N$ is generated by $\bigoplus_{i=0}^a N_i$, for some $a$,
and so $N_j=R_{ja}N_a$, for all $j\geq a$.
For $j\geq i\geq a$ we have natural maps of $R_a$-modules
$$\theta_{ji}: R^*_{ia}\otimes N_i \cong
 R^*_{ia}\otimes R^*_{ji}\otimes R_{ji}\otimes N_i \cong
 R^*_{ja}\otimes (R_{ji}\otimes N_i )
 \twoheadrightarrow R^*_{ja}\otimes N_j,$$
 where the tensor products are over the appropriate $R_k$.
 By the associativity of tensor products,
 $\theta_{ki}=\theta_{kj}\theta_{ji}$, for all $k\geq j\geq
 i\geq a$. Since each $N_i$ is a noetherian $R_i$-module,
  each $R^*_{ia}\otimes N_i$ is a noetherian $R_a$-module
  and so $\theta_{ji}$ is an isomorphism for
 all $j\geq i\gg 0$. Equivalently, $N_j\cong R_{ji}\otimes N_i$
 for all such $j\geq i$.

 Set $\Theta(\widetilde{N})
 =R^*_{j0}\otimes N_j\in R_0\lmod$ for some $j\gg 0$.
 Since any two preimages of $\widetilde{N}$ in $R_\Z\lgr $
 agree in high degree, $\Theta(\widetilde{N})$ is independent of the choice of
 $N$.  Moreover, as $R^*_{j0}=R^*_{k0}R_{kj}$,
 $$\phi(\Theta(\widetilde{N}))_{\geq j}
 \cong \bigoplus_{k\geq j} R_{k0}\otimes R^*_{j0}\otimes N_j
 \cong \bigoplus_{k\geq j} R_{kj}\otimes N_j = \bigoplus_{k\geq j}N_k,$$ and so
 $\Phi\Theta(\widetilde{N}) = \widetilde{N}.$ Checking
 that $\Theta$ and $\Phi$ are inverse equivalences is now routine.
 \end{proof}

\subsection{}  We remark that many of the standard techniques
and results concerned with associated graded modules for unital algebras
 extend routinely to $\Z$-algebras. These only appear in peripheral
 ways in this paper and so we refer the reader to \cite{GS2}
 for a discussion of these results.

\section{The main theorem}\label{sect-filt}

 \subsection{}\label{sect601}
In this section we prove the main theorem of
 the paper by proving Theorem~\ref{mainthm-intro} from the introduction.
Indeed we will prove more generally that a version of that theorem
holds for all values of $c\in \C$ that satisfy Hypothesis~\ref{morrat-hyp}.
As was true with Corollary~\ref{morrat-cor} and
Proposition~\ref{shiftonO}, the theorem will take slightly different forms
depending on whether $c\in \mathbb{Q}_{\leq -1}$ or not,
so it is convenient to separate
the cases with

\subsection{Hypothesis}\label{main-hyp}
{\it The element $c\in\C$  satisfies Hypothesis~\ref{morrat-hyp}
  but $c\not\in \mathbb{Q}_{\leq -1}$. }

\subsection{}\label{subsec-6.1}  Assume that Hypothesis~\ref{main-hyp} holds.
 By Corollary~\ref{morrat-cor} there is a Morita
equivalence $S_c: \UU_c\lmod \to  \UU_{c+1}\lmod$
given by $S_c(M) =  Q_c^{c+1} \otimes_{\UU_c} M$, where
 $Q_c^{c+1}=eH_{c+1}e_-\delta\subset D(\hr)\ast\WW$ is
considered as a right $\UU_c$-module via \eqref{conj}.
Following \eqref{zalgex2} we can therefore define a Morita $\Z$-algebra
  $B(c)= B=U_\Z$\label{B-ring-defn}  associated
 to the data
 $\{U_{c+i},\,  Q_{c+i}^{c+i+1}; i\in \NN\}$; thus
 $B=\bigoplus_{i\geq j\geq 0}B_{ij}$
 where,  for integers $i>j\geq 0$,
   \begin{equation}\label{Mij-defn}
    B_{jj}=\UU_{c+j} \qquad\text{and}\qquad
   B_{ij} =   \ Q_{c+i-1}^{c+i}Q_{c+i-2}^{c+i-1}\cdots
  Q_{c+j}^{c+j+1},
   \end{equation}
   where the multiplication
in taken in $D(\hr)\ast W$.
Note that, by Corollary~\ref{morrat-cor}, we have a natural isomorphism
\begin{equation}\label{tpdef}
B_{ij}\  \cong\
   Q_{c+i-1}^{c+i}\otimes_{\UU_{c+i-1}}Q_{c+i-2}^{c+i-1}
   \otimes_{\UU_{c+i-2}}\cdots\otimes_{\UU_{c+j+1}} Q_{c+j}^{c+j+1},
   \end{equation}
and so this does accord with the definition in \eqref{zalgex2}.

\subsection{The Main Theorem}\label{subsec-6.6}
Assume  that $c\in \C$ satisfies Hypothesis~\ref{main-hyp}.
The differential operator filtration $\ord$ on
$D(\hr)\ast\WW$, as defined in \eqref{order-filt-defn}, induces filtrations on
the subspaces  $B_{ij}$ and hence on $B$, which we will again write as $\ord$.
The fact that these filtrations are induced from that of
$D(\hr)\ast\WW$ ensures that  the associated graded object
$$\ogr B = \bigoplus_{i\geq j\geq 0}\ogr B_{ij}$$ is also a $\Z$-algebra.
Similarly,  recall from \eqref{A-1-defn}  the $\NN$-graded algebra
$A=\bigoplus_{i\geq 0} A^{i}$ associated to $\hi$. In this section it is
more convenient to use the isomorphic algebra
$A=\bigoplus_{i\geq 0} A^{i}\delta^i$ to which
  we  canonically associate
 the  $\Z$-algebra   \label{Aij-defn}
 $\widehat{A} = \bigoplus_{i\geq j \geq 0} A^{i-j}\delta^{i-j},$
 in the notation of  \eqref{zalgex1}.

\begin{thm} \label{main} Assume that $c\in \C$
 satisfies Hypothesis~\ref{main-hyp} and define $B$ and  $\widehat{A}$
  as above. Then:
\begin{enumerate}
\item There is an equivalence of  categories \
$\UU_c\lmod\ \xrightarrow{\sim}\ B\lqgr$.
\item There is an equality $\ogr B = e \widehat{A}e$  and hence
a graded $\Z$-algebra
isomorphism $\ogr B \cong \widehat{A}$.
\item
$\ogr B\lqgr\simeq\coh \hi$.
\end{enumerate} \end{thm}

Combining Theorem~\ref{main} with Corollary~\ref{morrat-cor} and
the isomorphism $U_c\cong U_{-c-1}$ from the proof of that result
  gives:
\begin{cor}\label{main-cor} {\rm (1)}
  Assume that $c\in \C$ satisfies Hypothesis~\ref{morrat-hyp}.
Then there exists a $\Z$-algebra $B'$ such that
  $U_c\lmod \simeq B'\lqgr$ and
$\ogr B\cong \widehat{A}$. Thus $\ogr B'\lqgr\simeq \coh(\hi)$.

{\rm (2)} If $c\in \C $ with $c\not\in \frac{1}{2}+\Z$,
then $H_c\lmod \simeq B''\lqgr$ and  $\ogr B''\lqgr\simeq \coh(\hi)$
for some $\Z$-algebra~$B''$.
 \qed
\end{cor}

\subsection{}\label{app-to-main}
Analogues of Theorem~\ref{main} also hold for certain important
 $U_{c+k}$-modules
and we will derive the theorem from one of these. The module in question
is the $(U_{c+k},H_c)$-bimodule $N(k)=B_{k0}eH_c$
with the induced $\ord$ filtration coming from the inclusion
$N(k)\subset D(\reg{\h})\ast \WW$.
Recall the definition of $J^d$ from \eqref{A-1-defn}.

\begin{prop}\label{pre-cohh}
Assume that $c\in\C$ satisfies Hypothesis~\ref{main-hyp} and let $k\in \NN$.
Then $\ogr N(k)=e J^k\delta^k$ as submodules of $\ogr D(\reg{\h})\ast \WW
=\C[\h\oplus \h^*]\ast \WW$.
\end{prop}

\subsection{Outline of the proof of the theorem and proposition}
\label{surjstrat} {\it For the rest of the section, we will assume that $c\in
\C$ satisfies Hypothesis~\ref{main-hyp}.}  Thus the notation from
\eqref{subsec-6.1} and \eqref{subsec-6.6} is available and, by
Corollary~\ref{morrat-cor}, $N(k)\cong B_{k0}\otimes_{U_c} eH_c$ is a
progenerative $(U_{c+k},\,H_c)$-bimodule. As will be  shown in  \eqref{subsec-6.21},
Theorem~\ref{main} follows easily from  Proposition~\ref{pre-cohh}, so we need
only discuss the proof of the latter result.  This  is nontrivial and  will
take most of the section but, in outline, is as follows.

It is easy to see that $eJ^k\delta^k\subseteq \ogr N(k)$
(see Lemma~\ref{thetainjA}). The other inclusion is considerably harder.
  The philosophy behind  the proof
 is to note that we can grade
both $J^{k}\delta^{k} $ and $ N(k)$  by the $\EE$-gradation.
This is not immediately
  useful since the graded pieces of the two sides are infinite dimensional
but both sides  have factor modules
for which the graded pieces are finite dimensional. For
$eJ^{k}\delta^{k}\cong J^k\delta^k$
 the factor is the module $\overline{J^{k}}\delta^{k}$  described by
Corollary~\ref{gr}, while the analogous factor $\overline{N(k)}$ of $\ogr N(k)$
 is described in \eqref{B-freeA} and Corollary~\ref{poincare-S2A} and
 is related to the   standard modules $\Delta_{c+k}(\mu)$.
The key observation is that
  these factors have the same Poincar\'e series and so  they are naturally
 isomorphic as graded vector spaces. The proof of the theorem then amounts
 to lifting this isomorphism to give the desired
 equality $eJ^k\delta^k = \ogr N(k)$.

This also shows that the result has to be non-trivial. Indeed,
an alternative proof of the proposition (or the theorem)
would also provide an alternative proof to  a number of the
results from \cite{hai3}.

\subsection{}\label{abstract-products}
 We start with two elementary observations that will be used frequently.
 If $R=\bigcup F^iR$ is a filtered ring and $r\in F^mR\smallsetminus F^{m-1}R$,
  we write $\sigma(r)=[r+F^{m-1}R]\in \gr_F^mR$ for the
  {\it principal symbol}\label{princ-symbol-defn} of $r$.

 \begin{lem} Let $R=\bigcup F^iR$ be a filtered $k$-algebra, for a field $k$.

 {\rm (1)} Let $A$, $B$ be subspaces of   $R $ and
 give $A$, $B$ and $AB$ the induced filtration $F$. Then
 $(\gr_FA)(\gr_FB)\subseteq \gr_FAB$, as subspaces of $\gr_FR$. Indeed,
 if $a\in A$ and $b\in B$ satisfy  $\sigma(a)\sigma(b)\not=0$,
  then $\sigma(a)\sigma(b)= \sigma(ab)$.

 {\rm (2)} Suppose that $A=\bigcup F^iA$ is a filtered right $R$-module
 and that $B=\bigcup F^iB$ is a filtered left $R$-module and give the
vector space $A\otimes_RB$ the \emph{tensor product filtration}:
$F^n(A\otimes B)=\sum_j F^jA\otimes F^{n-j}B$. Then there is a
natural surjection $\gr_FA\otimes_{\mathrm{gr}\,R}\gr_FB\twoheadrightarrow
\gr_F(A\otimes_RB)$.
 \end{lem}

\begin{proof} (1) Identify $\gr_FA=\bigoplus
(F^nA+F^{n-1}R)/F^{n-1}R\subseteq \gr_FR$ so that the result makes sense.
Suppose that $\bar{a}\in \gr_F^nA$ and $\bar{b}\in \gr_F^mB$ are such that
$\bar{a}\bar{b}\not=0$ in $\gr_FR$. Lift $\bar{a}$ to $a\in F^nA$ and $\bar{b}$
to $b\in F^mB$.
Then, as elements of $\gr_FR$, one has
$\bar{a}\bar{b} = [a+F^{n-1}R][b+F^{m-1}R]
\subseteq [ab+F^{n+m-1}R]$. Since $\bar{a}\bar{b}\not=0$,
$ab\in F^{n+m}R\smallsetminus F^{n+m-1}R$, whence
$\bar{a}\bar{b}=\sigma(ab)$ is the image of $ab$ in $\gr_F(AB)$.

(2) Define a map $\rho: \gr_FA\times \gr_FB \to \gr_F(A\otimes_RB)$
by $\rho(\bar{a}, \bar{b} ) = [a\otimes b + F^{n+m-1}(A\otimes B)]$,
for  $\bar{a}\in \gr_F^nA$, $\bar{b}\in \gr_F^mB$ and where  the rest of the
notation
is the same as for part (1). This   clearly defines a $\C$-bilinear map
that is $\gr_F R$-balanced in the sense that
$\rho(\bar{a}\bar{r}, \bar{b}) = \rho(\bar{a}, \bar{r}  \bar{b})$
for $\bar{r}\in \gr^s_FR$. By universality, $\rho$ therefore   induces a map
$ \gr_FA\otimes_{\gr R}\gr_FB\to \gr_F(A\otimes_RB)$.
It is surjective since
$F^{n+m}(A\otimes B)/ F^{n+m-1}(A\otimes B)$ is spanned by elements of the
given form  $[a\otimes b + F^{n+m-1}(A\otimes B)].$
\end{proof}

\subsection{Lemma}\label{grade-elements}
Let $R=\bigcup_{i\geq 0}F^iR$  be a filtered ring, pick $r\in R$
and let $ I$  be a subset of $R$.
Under  the induced filtrations, $\gr_F(rI) = \sigma(r)\gr_F(I)$  in the
following cases:
\begin{enumerate}
\item[(1)] $\sigma(r)$ is regular in $\gr_FR$;
\item[(2)] $ r=r^2 \in F^0(R)$ and $rI\subseteq I$.
\end{enumerate}

\begin{proof}   Assume that $r\in F^sR\smallsetminus F^{s-1}R$.
 We claim that, in both cases, it suffices to prove that
$F^n(rI) = rF^{n-s}I$ for   all   $n\geq s$.
 Indeed, if this is true then the identity  $F^m(rI) = rI\cap F^mR$
 implies that the $n^{\mathrm{th}}$ summand of $\gr(rI)$ equals
$$\frac{F^n(rI)}{F^{n-1}(rI)} \ = \  \frac{F^n(rI)}{F^n(rI)\cap F^{n-1}R}
      \  \cong\  \frac{ F^n(rI) + F^{n-1}R}{F^{n-1}R}
        \ = \  \frac{rF^{n-s}I + F^{n-1}R}{F^{n-1}R},$$
 which is the $n^{\mathrm{th}}$  summand of $\sigma(r)gr(I)$.

(1) In this case,
$rt\in F^n(rI) = rI\cap F^n(R) \Leftrightarrow
 t\in I\  \mathrm{and}\ t\in F^{n-s}R,$
 as required.

(2) Here,
$rF^{n}I \subseteq F^n(rI)$
whence $rF^{n}I = r^2F^{n}I \subseteq  rF^n(rI)
\subseteq rF^nI$. Since $rF^n(rI) = F^n(rI) $ this implies that
 $rF^n(I)=F^n(rI)$.
\end{proof}

\noindent
{\bf Example}. It is easy to check that some hypotheses
are required for the lemma to hold.
For example, filter the polynomial ring $R=\C[x,y]$ by
$x,xy\in F^0R$ but $y\in F^1R$.
Then  $x,xy\in F^0(xR)$, yet
$xy\not\in \sigma(x)\gr_FR$.

\subsection{} We now turn to the proof of Proposition~\ref{pre-cohh}.
As was mentioned in \eqref{surjstrat} the inclusion
$J^k\delta^k e \subseteq \ogr N(k)$ is easy.

\begin{lem}\label{thetainjA}
{\rm (1)}   For $i\geq j\geq 0$ we have
$e(A^{i-j}\delta^{i-j})e
\subseteq \ogr B_{ij}$.

{\rm (2)} The inclusion of part (1)  is an equality for $i=j$ and
for  $i=j+1$.

{\rm (3)}
For $k\geq 0$ there is an inclusion
$eJ^{k}\delta^{k}
\subseteq \ogr N(k)$ of left $eA^0e$-modules. This is an equality for $k=0$.
\end{lem}

\begin{proof}  (2)
By the PBW Theorem~\ref{PBW},
 $\ogr B_{ii} = e(\cxy \ast \WW)e$ and so the claim holds for $i=j$.
Similarly,   since $ e,\delta\in\ord^0(D(\reg{\h})\ast\WW)$
and  $\delta$ is regular in $\ogr(D(\reg{\h})\ast\WW)$,
Lemma~\ref{grade-elements} implies that
$$
\ogr B_{j+1,j} = \ogr( eH_{c+j+1}e_- \delta)  =
\ogr(eH_{c+j+1}e_-)\delta=
 e(\ogr H_{c+j+1})e_-\delta = e(\cxy\ast\WW)e_-\delta = eA^1\delta e.$$

(1) Combining part (2) with  Lemma~\ref{abstract-products}(1) and induction
shows that
$$
 (eA^1\delta^1e)^{i-j} \ = \
 \ogr B_{i,i-1}\ogr B_{i-1,i-2}\cdots \ogr B_{j+1,j}
\ \subseteq\  \ogr \left( B_{i,i-1} \cdots   B_{j+1,j}  \right)
\ = \ \ogr B_{ij}.$$

(3) When $k=0$, the assertion   $eJ^k\delta^{k}
= \ogr N(k)$ is just the statement that
 $e\cxy = e(\cxy\ast \WW) $.
 When  $k>0$,  part (i) and
 Lemma~\ref{abstract-products}
give $eJ^k\delta^k = eA^k\delta^ke\C[\h\oplus\h^*]\ast \WW
 \subseteq \ogr B_{k0}\ogr(eH_c)\subseteq \ogr N(k).$
\end{proof}

\subsection{}\label{h-defn}
The next several results will be aimed at getting a more detailed understanding
of the bimodule structure of $N(k)$ and its factors.
For the most part we are interested in their graded structure
for which the actions of the elements   $\hh_{c+t}\in H_{c+t}$ from
\eqref{hdefn} are particularly useful.
 Given an $(\UU_{c+s},\UU_{c+t})$-bimodule $M$, define
$$\hh \bigdot m = \hh_{c+s} m - m\hh_{c+t} \text{ for any $m\in M$}.$$
When $s=t$ this is just  the adjoint action of $\hh_{c+s}$ on $M$.

\begin{lem} \label{diaggrad} {\rm (1)} $ e\hh_{c+t-1} e =
\delta^{-1}e_-\hh_{c+t} e_-\delta$.

{\rm (2)} The action of $\hh$ is diagonalisable
on the modules  $N(i)$, $B_{ij}$ and $M(i)=H_{c+i}eB_{i0}$,
for any $i\geq j\geq 0$.
\end{lem}

\begin{proof} (1) Use the first paragraph of the proof of
 \cite[Theorem~4.10]{gordc}.

(2) We start with the $B_{ij}$.  If $b_1\in B_{i\ell}$ and
$b_2\in B_{\ell j}$, then
$\hh\bigdot (b_1b_2) = (\hh\bigdot b_1)b_2 + b_1(\hh\bigdot b_2)$. Thus,
by induction, it suffices to prove the result for each $B_{t,t-1}
= e H_{c+t}\delta e$.
Clearly $e\hh=\hh e$. Thus, by part (1),
  for any  $m\in H_{c+t}$ we have
 \begin{equation}\label{diaggrad1}
\hh \bigdot e m \delta e  \  =\
  \hh_{c+t}e m \delta e - e m \delta e \hh_{c+t-1}  \ =
\   e \hh_{c+t} m \delta e - e m \delta (\delta^{-1}e_-\hh_{c+t}e_-\delta)
 \   =\  e ([\hh_{c+t}, m])  \delta e.
   \end{equation}
By \eqref{grading} $H_{c+t}$ is diagonalisable under the
adjoint $\hh_{c+t}$-action and so   the result for $B_{ij}$ follows.
The same argument works for the modules $N(i)$ and $M(i)$ if one uses
 the decompositions $N(i)=(B_{i0})(eH_c)$ and
$M(i)=(H_{c+i}e)(B_{i0})$.  \end{proof}

\subsection{}\label{B-freeA}
 The factors of $N(k)$ that most interest us are defined as follows.
 Since $N(k)$ is a
$(\UU_{c+k},H_c)$-bimodule, the embeddings $\C[\h]^{\WW}\hookrightarrow
\UU_{c+k}$ and $\C[\h^*]\hookrightarrow H_c$ make $N(k)$ into a
$(\C[\h]^{\WW},\, \C[\h^*])$-bimodule.
Let $\C$ be the trivial module over either $\C[\h]^{\WW}$
or $\C[\h^*]$ and set
$ \overline{N(k)} = \C\otimes_{\C[\h]^{\WW}} N(k)$ and $
\underline{N(k)} = N(k) \otimes_{\C[\h^*]} \C.$
As $\C$ is a graded  $\hh$-module,
the adjoint action of $\hh$ on $N(k)$ from Lemma~\ref{diaggrad}
 induces a $\Z$-grading,
again called the $\hh$-grading, on both
$\overline{N(k)}$  and $\underline{N(k)}$.   If an element
$b$ from any of these three modules has
degree $n$ in this grading we write $\hdeg(b)=n$.\label{hdeg-defn}
The reader should note  that, as will be explained in
\eqref{poincare-S2A}, this is {\it not} the same as the $\EE$-gradation
on these modules.

 The next result gives the elementary properties of these  modules.

\begin{lem}\label{Bbar-freeA} {\rm(1)} For any $i\geq j\geq 0$,
 $B_{ij}\subseteq U_{c+i}\cap U_{c+j}$.
\begin{enumerate}
\item[(2)]  For $k\geq 0$, both $N(k)$ and $U_{c+k}$ are free left
$\C[\h]^\WW$-modules, while $N(k)$ is a free right $\C[\h^*]$-module.
\item[(3)]   $\underline{N(k)}$ is a  finitely generated, free left
$\C[\h]^\WW$-module.
\item[(4)]  Similarly,
$\overline{N(k)}$ is a finitely generated, free
right $\C[\h^*]$-module.
\end{enumerate}
\end{lem}

\begin{proof} We will use frequently and without comment the fact that
$\C[\h^*]$ is a free $\C[\h^*]^\WW$-module. Moreover, as
$\C[\h^*]^\WW$ is a polynomial ring, any projective
$\C[\h^*]^\WW$ is free by the Quillen-Suslin Theorem.

(1) By induction, we may assume that $i=j+1$.
The inclusion $B_{ij}=eH_{c+i}\delta e\subseteq U_{c+i}$ is immediate.
If  $p\in H_{c+i}$ then, by \eqref{conj},
 $$epe_-\delta\ =\ e\delta^{-1}\delta p
e_- \delta\ =\ \delta^{-1} e_- \delta p e_-\delta\ \in\ \delta^{-1}
e_-H_{c+i}e_-\delta\ =\ U_{c+j}.$$

 (2) By the PBW Theorem~\ref{PBW}, each $H_d$ is free as a left $\C[\h]$-module
  and  as a right  $\C[\h^*]$-module. Therefore, $H_d$ is a free left
  $\C[\h]^W$-module as is its summand
 $H_de$. Under the left action of $\WW$, $(H_de)^\WW=eH_de$
since, if $fe\in (H_{d}e)^W$,
then $fe =|\WW|^{-1}\sum_{w\in \WW} wfe = efe$.
But $ (H_{d}e)^W$ is a $\WW$-module summand of $H_{d}e$, while
the actions of $\WW$ and $\C[\h]^\WW$ commute. Thus
$U_{d}=(H_{d}e)^\WW$ is  a $\C[\h]^\WW$-module summand of
$H_{d}e$ and hence is  free.
By Corollary~\ref{morrat-cor}, $N(k)\cong B_{k0}\otimes_{U_c}eH_c$
is a projective left $U_{c+k}$-module
and hence a free  left $\C[\h]^\WW$-module.

On the other hand, $N(k)$ is a projective right $H_c$-module and
hence   a projective right $\C[\h^*]$-module.

(3)
Set $X=H_{c}\otimes_{\C[\h^*]} \C$.
Clearly $X\in \mathcal{O}_{c}$ in the sense of \eqref{standard-defn}
and,
by   \eqref{PBW},
$X  \cong \C[\h]\otimes_{\C} \C \WW$  as   left $\C[\h]\ast \WW$-modules.
Thus $X$ is a finitely generated
free left $\C[\h]$-module
and so, by \cite[Proposition~2.21]{GGOR}, $X$ has a filtration
whose factors are standard modules.

 By definition, $\underline{N(k)} = eM$
where $M = \widetilde{S}_{c+k-1}\circ \cdots \circ \widetilde{S}_{c}(X)$,
in the notation of \eqref{shift-defn}.  By
Proposition~\ref{shiftonO}  $M$
also has a finite filtration by standard modules and so
\cite[Proposition~2.21]{GGOR} shows that $M$ is a
finitely generated free module over   $\C[\h]$ and hence over
 $\C[\h]^{\WW}$. Thus, so is its summand $eM$.

(4) We first show that $N(k)$ is a finitely generated right module
over $R=(\C[\h]^\WW)^{\mathrm{op}} \otimes_\C \C[\h^*]$.
By part~(1), $B_{k0}\subseteq U_{c}$ and so $N(k)\subseteq eH_c$. Thus
$\ogr N(k) \subseteq \ogr H_{c}
= \C[\h\oplus \h^*]\ast \WW$, which is certainly a noetherian
 $\C[\h]^\WW\otimes \C[\h^*]$-module.
Since the $\ord$ filtration on $N(k)$
  is the one induced from  $D(\reg{\h})\ast \WW$, the actions of
 $\C[\h]^\WW$ and $\C[\h^*]$ on $\ogr N(k)$ are  the natural ones
 induced from the actions of those rings on
 $N(k)\subset D(\reg{\h})\ast\WW$. In other words,
  the given $R$-module structure of
$\ogr N(k) $ is   the one induced from the
$R$-module structure of $N(k)$.
Since the former module is finitely generated, so is the latter.

Let $y_1,\ldots ,y_{n-1}$ be the generators of $\C[\h^*]$ and let
$q_1, \ldots ,q_{n-1}$ be the
fundamental invariants of $\C[\h]^{\WW}$.
By (2), the $\{y_j\}$ form an r-sequence in $N(k)$, while (3)
implies that the $\{q_j\}$ form an r-sequence in the factor
$\underline{N(k)}=N(k)/\sum N(k)y_j$ as a module over
$\C[\h]^\WW=R/\sum y_jR$. Thus
 $\Sigma = \{y_\ell ,q_{m} : 1\leq \ell, m\leq n-1\}$ is a regular sequence
for the right $R$-module
$N(k)$.
In particular, if $\mathfrak{n}=\sum y_iR+q_jR$, then
$\Sigma$ is  an r-sequence for the
 $R_{\mathfrak{n}}$-module $N(k)_\mathfrak{n}$.
By the Auslander-Buchsbaum formula
\cite[Ex.~4, p.114]{matsumura},
 $N(k)_\mathfrak{n}$ is therefore free as a
$R_{\mathfrak{n}}$-module.

Finally, consider $\overline{N(k)}=N(k)/\sum q_jN(k)$.
   Under the induced $\hh$-grading,  $\overline{N(k)}$ is a finitely
generated, graded   $\C[\h^*]$-module and  so corresponds to a
$\C^*$-equivariant coherent sheaf on $\h^*$. As a result the locus where
$\overline{N(k)}$ is not  free  is a $\C^*$-stable closed subvariety of
$\h^*$. If this locus  is non-empty it must  contain the unique
$\C^*$-fixed point $\mathfrak{p}=(y_1,\dots,y_{n-1})$ for this expanding
 $\C^*$-action. But then
$(\overline{N(k)})_{\mathfrak{p}}$ would not be free, contradicting the
conclusion of the last paragraph.\end{proof}

\subsection{}\label{poincare-start}
We next need to understand the graded structure of
the modules $\overline{N(k)}$ and $\underline{N(k)}$ under the
$\hh$-grading. To do this,  we express
 $\underline{N(0)}$  as a weighted sum of standard modules in the Grothendieck
 group $G_0(U_c)$ and then to use Proposition~\ref{shiftonO} to write
$\underline{N(k)}=B_{k0}\otimes \underline{N(0)}$ in a similar manner.
This is quite delicate since there are some subtle shifts involved and we
first want to understand these shifts  for $B_{ij}\otimes \Delta_c(\mu)$.

We will need to  work with
 the following  graded version  $\widetilde{\OO}_d$\label{cat-O-gr-defn}
  of $\OO_d$ constructed in  \cite[Section~2.4]{GGOR}.
 The objects $M$ in $\widetilde{\OO}_d$ are
finitely generated $H_d$-modules on which $\C[\h^*]$ acts  locally nilpotently and
 which come equipped with a
$\mathbb{Z}$-grading $M = \bigoplus_{r\in \Z} M_{r}$ such that
$p M_{r} \subseteq M_{r+\ell}$ for each
$p\in H_d$ with   $\Edeg(p)=\ell$.
The morphisms
are homogeneous $H_d$-module homomorphisms of degree
zero. A {\it  graded standard module}
\label{graded-standard-defn}
$\widetilde{\Delta}_d(\mu)$, isomorphic to
$\Delta_d(\mu)$ as an ungraded module,  is given by setting
$\widetilde{\Delta}_d(\mu)_r = \C[\h]_r\otimes \mu$.
By local nilpotence and finite generation, each weight space of a module
 $M\in \widetilde{\OO}_d$ is finite dimensional and so   $M$ has
 a well-defined Poincar\'{e} series. There is a degree shift functor\label{module-shift-defn}
 $[1]$ in $\widetilde{\OO}_d$
 defined by $M[1]_r = M_{r-1}$. By abuse of notation, $\widetilde{\OO}_d$
will also denote the corresponding category of graded $U_d$-modules.

\begin{lem}\label{standAAA}
 Fix  $i\geq j\geq 0$ and $\mu\in\irr{\WW}$.
 Give $B_{ij}$ the adjoint $\hh$-grading and let
 $B_{ij}\otimes_{U_{c+j}}e\widetilde{\Delta}_{c+j}(\mu)$
 have the grading this induces. Then
 $B_{ij}\otimes_{U_{c+j}}e\widetilde{\Delta}_{c+j}(\mu)\in \widetilde{\OO}_{c+i}$
 and, as elements of that category,
 $$B_{ij}\otimes_{U_{c+j}}e\widetilde{\Delta}_{c+j}(\mu)
 \cong e\widetilde{\Delta}_{c+i}[(i-j)(n(\mu)-n(\mu^t))].$$
 \end{lem}

\begin{proof} Write $\nabla = B_{ij}\otimes_{U_{c+j}}e\widetilde{\Delta}_{c+j}(\mu)$
 and let $\deg_{c+u} $ denote the degree function in  $\widetilde{\OO}_{c+u}$.
By hypothesis, the graded structure of an element $b\otimes v\in \nabla$
is given by $\deg (b\otimes v) = \hdeg(b)+\deg_{c+j}(v)$.
Proposition~\ref{shiftonO} implies
that (as ungraded modules)
\begin{equation}\label{grot3}
\nabla
=S_{c+i}\circ \cdots \circ S_{c+j+1}(e\Delta_{c+j}(\mu)) \cong
e\Delta_{c+i}(\mu).
\end{equation} Thus, under its given grading,  $\nabla \in \widetilde{\OO}_{c+i}$.

Unfortunately, it is not easy  to write the generator $e\otimes \mu$
of $e\Delta_{c+i}(\mu)$ as an element of $\nabla$ and for this reason the
shift in the grading in \eqref{grot3} is subtle.
In order to understand this we will use the
canonical grading  from \eqref{cangrad} and we
write the corresponding degree function as $\deg_{\mathrm{can}}$.
The advantage of this grading is that it is simply given
by the left multiplication of $\hh_{c+i}$.  Thus, as  \eqref{grot3} is an
 isomorphism of left $U_{c+i}$-modules and hence of left
 $\C[\hh_{c+i}]$-modules, it is automatically a graded isomorphism
  under the canonical grading.

Since $\h^*$ has $\EE$-degree $1$,
the canonical grading on $\Delta_d(\mu)$, for any $d\in \C$,  is a shift of the
grading on $\widetilde{\Delta}_d(\mu)$. The shift is easy to compute.
 By definition, the generator $1\otimes \mu$ of $\widetilde{\Delta}_d(\mu)$  has
$\deg_d(1\otimes \mu)=0$ whereas, by
Proposition~\ref{subsec-3.10},
the generator $1\otimes\mu$ of $\Delta_d(\mu)$
has $$\deg_{\mathrm{can}}(1\otimes \mu)=
D(d,\mu) =(n-1)/2+ d(n(\mu)-n({\mu^t})).$$
We may therefore regard $\Delta_d(\mu)$ as being in
 $\widetilde{\OO}_d$, in which case
  \begin{equation} \label{shft}
 \Delta_d(\mu)= \widetilde{\Delta}_d(\mu)[D(d,\mu)].\end{equation}

  Let $b\in B_{ij}$ with $\hdeg(b)=r$
 and  suppose that $v\in e\Delta_{c+j}(\mu)$ has $\deg_{\mathrm{can}}(v)=s$.
 Then
 $$\hh_{c+i}\cdot b\otimes v = (\hh\bigdot b)\otimes v +b\,\hh_{c+j}\otimes v
\  =\  (\hh\bigdot b)\otimes v +b\otimes \hh_{c+j}v
=(r+s)b\otimes v.$$
Thus $\deg_{\mathrm{can}}(b\otimes v) =\hdeg(b)+
\deg_{\mathrm{can}}(v)$.
Finally,  \eqref{shft} implies that
\begin{eqnarray*}
\deg_{c+i}(b\otimes v) &=& \deg_{\mathrm{can}}(b\otimes v) - D(c+i, \mu) \ = \
 \hdeg (b) + \deg_{\mathrm{can}}(v)  - D(c+i,\mu) \\ &=&
  \hdeg(b) + \deg_{c+j}(v) + D(c+j,\mu) - D(c+i,\mu) \\
&=& \deg (b\otimes v) + (j-i)(n(\mu) - n(\mu^t)),
\end{eqnarray*}
as required.
\end{proof}

\subsection{}\label{poincare-sa-sect}
Given a $\Z$-graded complex vector space $M =
\bigoplus_{r\in\Z}M_r$ such that $\dim_{\C} M_r$ is finite for all $r$
then, as in \eqref{gr},  we define
the  Poincar\'{e} series\label{formal-Poincare-defn}
 of $M$ to be $p(M,v) = \sum v^r \dim_{\C} M_r.$
Each $N(k)$ is graded via the adjoint $\hh$ action from
\eqref{hdeg-defn}, although
of course the summands are infinite dimensional. Thus in order to understand
the more  detailed structure of $N(k)$ and $\ogr N(k)$
we will consider the  Poincar\'e series of
the  factor modules $\overline{N(k)}$ and $\underline{N(k)}$.

\begin{prop}\label{poincare-SA}
 If $\overline{N(k)}$ as graded via
the adjoint $\hh$ action on $N(k)$, then
its Poincar\'{e} series
 is
 \begin{equation}\label{poincare-SA00}p(\overline{N(k)}, v) =
\frac{\sum_{\mu} f_{\mu}(1)f_{\mu}(v^{-1}) v^{-k(n(\mu)
- n(\mu^t))}[n]_v!}{\prod_{i=2}^n
(1-v^{-i})}.
\end{equation}
\end{prop}

\begin{proof}
We first calculate the Poincar\'{e} series for $\underline{N(k)}$,
and we begin with $\underline{N(0)}$.
 As in the proof of
 Lemma~\ref{Bbar-freeA}(3),   $X = H_{c}  \otimes_{\C[\h^*]} \C $
 is an object of
$\widetilde{\OO}_{c}$, where the grading is the natural one defined by
$\deg(1\otimes 1)=0$.   By construction,  $eX \cong \underline{N(0)}$
and this is a {\it graded} isomorphism since the adjoint  $\hh$-graded structure of
$\underline{N(0)}=U_c/I$ is simply defined by $\hdeg(e)=0$.
Thus, as elements of the Grothendieck group  $G_0(\widetilde{\OO}_{c})$,
we can write $ [X] =
\sum_{\mu} p_{\mu} [\widetilde{\Delta}_{c}(\mu)]$ for some
$p_{\mu} \in \Z
[v,v^{-1}]$. By \eqref{PBW} we have a graded isomorphism
 $X\cong \C[\h]\otimes \C\WW$.
  Applying $(\C\otimes_{\C[\h]}-)$ to the formula  $ [X] =
\sum_\mu p_\mu  [\widetilde{\Delta}_{c}(\mu)]$ therefore yields $\C \WW
= \sum_{\mu} p_{\mu} [\mu].$ It follows from \eqref{fakedegrees2} that
$p_{\mu} = f_{\mu}(1)$ and so $
[\,\underline{N(0)}\,] = \sum_{\mu} f_{\mu}(1) [e\widetilde{\Delta}_{c}(\mu)] .
$
Combining this formula with Lemma~\ref{standAAA} shows that
\begin{equation}\label{verygrot}
[\,\underline{N(k)}\,] =
\sum_{\mu} f_{\mu}(1)v^{k(n(\mu)-n(\mu^t))} [e\widetilde{\Delta}_{c+k}(\mu)] .
\end{equation}

The Poincar\'e series of $N(k)$ is now easy to compute. First,
 in the {\it canonical grading},   \eqref{polystand} shows that
$$p(\Delta_d(\mu), v, W) = v^{D(d,\mu)}\frac{\sum_{\lambda}
f_{\lambda}(v) [\lambda\otimes \mu]}{\prod_{i=2}^n(1-v^i)}
\qquad\mathrm{and\ so}\qquad  p(e\Delta_d(\mu), v)=
 v^{D(d,\mu)} \frac{f_{\mu}(v)}{\prod_{i=2}^n (1-v^i)}
 $$ for any $d\in \C$.
 Therefore,  \eqref{shft} implies that
 $p(e\widetilde{\Delta}_d(\mu), v)=
  f_{\mu}(v) \prod_{i=2}^n (1-v^i)^{-1}
 $ in the graded category $\widetilde{\OO}_{d}$.
Combined with \eqref{verygrot} this shows that
\begin{equation}\label{wrongsideformulaA}
p(\underline{N(k)},v)  \ = \
\frac{\sum_{\mu} f_{\mu}(1)f_{\mu}(v) v^{k(n(\mu) -
n(\mu^t))}}{\prod_{i=2}^n
(1-v^{i})}.
\end{equation}

Finally, we  calculate the Poincar\'{e} series of $\overline{N(k)}$.
By Lemma~\ref{Bbar-freeA}(2,3), an
 $\hh$-homogeneous basis for this module is given by lifting a homogeneous
 $\C$-basis from $\overline{N(k)}\otimes_{\C[\h^*]}
  \C = \C \otimes_{\C[\h]^{\WW}} {\underline{N(k)}}.$
 Thus, combining
     \eqref{wrongsideformulaA} with the formul\ae\
    $p(\C[\h]^{\WW}, v) =  \prod_{i=2}^n (1-v^i)^{-1}$ and
 $ p(\C[\h^*], v) = (1-v^{-1})^{n-1}$
 gives
\begin{equation}\label{wrongsideformula2A}
p(\overline{N(k)}, v) =
\frac{\sum_{\mu} f_{\mu}(1)f_{\mu}(v) v^{k(n(\mu) -
n(\mu^t))}}{(1-v^{-1})^{n-1}}.
\end{equation} This needs to be adjusted to yield  \eqref{poincare-SA00}.
  Set $N=n(n-1)/2$.  Then  Lemma~\ref{babyverma}(1)
and \eqref{fakedegrees2} combine  to show that
\begin{eqnarray*} \sum_{\mu} f_{\mu}(1)f_{\mu}(v) v^{k(n(\mu) - n(\mu^t))} &=&
\sum_{\mu} f_{\mu^t}(1)f_{\mu^t}(v^{-1}) v^{k(n(\mu) - n(\mu^t))}\\
 &=&
v^N\sum_{\lambda} f_{\lambda}(1)f_{\lambda}(v^{-1})
v^{k(n(\lambda^t) - n(\lambda))}.\end{eqnarray*} Moreover,
rearranging  \eqref{factorial-defn} gives
$$ [n]_v!
\ = \  \frac{\prod_{i=1}^n (1-v^i)}{(1-v)^n} \ =\
v^N\frac{\prod_{i=1}^n (1-v^{-i})}{(1-v^{-1})^n}.$$
Combining these formul\ae\ with \eqref{wrongsideformula2A} gives
\eqref{poincare-SA00}.
\end{proof}

\subsection{}\label{poincare-S2A}
  Recall the Euler gradation $\Edeg$ on $D(\reg{\h})\ast \WW$ and its subrings
  from  \eqref{gradingsec}.
  Since $e$, $e_-$ and $\delta$ are homogeneous under this action,
each $Q_{c+\ell}^{c+\ell+1}$ and hence each $B_{ij}$ and $N(k)$
is also graded under this action.
As in \eqref{gradingsec}, this induces a graded
structure, again called $\Edeg$,  on $\ogr B_{ij}$ and $\ogr N(k)$.
 Since the fundamental invariants of $\C[\h]^\WW$ are
 $\EE$-homogeneous,  the $\EE$-grading on $N(k)$  descends to gradings
  on $\overline{N(k)}$ and $\underline{N(k)}$.
   Similarly,  each $A^{u}\delta^{u}$ and $J^u\delta^u$
  has an $\EE$-grading induced from that on $\C[\h\oplus\h^*]$
and hence so does $A=\bigoplus_{u\geq 0}A^u\delta^u.$

However, the $\EE$-grading on $B_{k0}$ and hence on $N(k)$ is {\it not}
equal to the adjoint  $\hh$-grading.
The problem is that, in \eqref{diaggrad1}, the adjoint $\hh$ action does
not ``see''
the element $\delta$. Thus if we wish to relate the Poincar\'e series of
$N(k)$ to that of   $J^{k}\delta^{k}$ we need the following  slight
modification of   Proposition~\ref{poincare-SA}.

\begin{cor} Let $k\geq 0$,  set $N=n(n-1)/2$ and  write  $K=kN$.
\begin{enumerate}
\item{} If $b\in B_{ij}$ for $i\geq j\geq 0$  is homogeneous under the
$\hh$-grading then it is homogeneous in the $\EE$-grading and
$\Edeg b = (i-j)N + \hdeg b.$
\item{}
Under the $\EE$-grading,
$\overline{N(k)}$ has  Poincar\'{e} series
$$p(\overline{N(k)}, v) =  v^{K}
\frac{\sum_{\mu} f_{\mu}(1)f_{\mu}(v^{-1}) v^{-k(n(\mu)
- n(\mu^t))}[n]_v!}{\prod_{i=2}^n
(1-v^{-i})}.
$$\vskip 5pt \noindent
while $\underline{N(k)}$ has Poincare series
 $p(\underline{N(k)},v) = v^{K}\displaystyle
\frac{\sum_{\mu} f_{\mu}(1)f_{\mu}(v) v^{k(n(\mu) -
n(\mu^t))}}{\prod_{i=2}^n
(1-v^{i})}.$\end{enumerate}\end{cor}

\medskip
\begin{proof} (1) If $b_1\in B_{ik}$ and $b_2\in B_{kj}$ then
$\hh\bigdot (b_1b_2) = (\hh\bigdot b_1)b_2 + b_1(\hh \bigdot b_2)$
and  $[\EE, b_1b_2] = [\EE,b_1]b_2 + b_1[\EE,b_2]$.
 By induction, it therefore suffices to  prove the result
 when $b=em\delta e\in B_{k,k-1}=
 eH_{c+k} \delta e$, for some $k>0$.
By \eqref{diaggrad1} we see that
$\hh\bigdot b = e[\hh_{c+k}, m]\delta e$ whereas
$[\EE, b ] = e[\EE, m] \delta e + em[\EE, \delta] e$.
By \eqref{gradingsec},  $ [\hh_{c+k}, m]=[\EE, m]$
and so the two gradings
differ by $\Edeg \delta = N$.

(2)  This
follows from part (1) combined with
 Proposition~\ref{poincare-SA}, respectively
 \eqref{wrongsideformulaA}.\end{proof}

\subsection{}\label{filter-injA}  Fix $k\geq 0$ and for notational simplicity
write $\mathcal{J}=e J^k\delta^k$ and $\mathcal{N}=N(k)$.
The final step in the proof of
Proposition~\ref{pre-cohh}  is to show that the
inclusion $\Theta:  \mathcal{J} \hookrightarrow \ogr \mathcal{N}$
from Lemma~\ref{thetainjA}(3)  is surjective. In order to effectively use
 Corollary~\ref{poincare-S2A}, we do this by lifting $\Theta$ to a
 $\C[\h]^{\WW}$-module map $\theta: \mathcal{J}\to \mathcal{N}$.

The order filtration
on $D(\hr)\ast \WW$  induces a
graded structure on $\ogr D(\hr)\ast \WW\cong
\C[\hr\oplus\h^*]\ast\WW$ and hence on
$\ogr \caln$, which we call the \textit{order gradation}; thus
  $\deg_{\ord} (\C[\h]\ast \WW)=0$,
 while $\deg_{\ord } \h=1$.
We will use the same terminology for the induced grading on the
rings $A^0=\C[\h\oplus \h^*]^W$ and $A$ and the module $\calj$.

Let $\caln^m=\ord^m\caln$ denote the elements in $\caln$ of order $\leq m$.
Similarly, write $\calj = \bigoplus_{m\geq 0} \ogr^m \calj$ for
the graded structure of $\calj$ under the $\ord$ gradation
and write  the induced order filtration  as $\calj = \bigcup \calj^m$, for
$\calj^m = \ord^m \calj= \bigoplus_{0\leq i\leq m} \ogr^i\calj$.

\begin{lem} There exists  an injective map $\theta :
\calj\hookrightarrow \caln$ of left $\mathbb C[\h]^W$-modules
such that:
\begin{enumerate}
\item $\theta$  is a graded homomorphism under the $\EE$-gradation
and is a filtered homomorphism  under the order filtration.
\item The
associated graded map
$\ogr\, \theta: \calj\to \ogr \caln$ induced by $\theta$
is precisely $\ogr \theta = \Theta$.
\end{enumerate}
\end{lem}

\begin{proof}
Trivially, $\Theta$ is  an $\EE$-graded map (by
 which we always mean a graded map of degree zero), as well
 as being graded under the $\ord$  gradation.
For any $m$,  $\ogr^m\calj$ is   an $\EE$-graded $\mathbb C[\h]^W$-module.
By Corollary~\ref{hi-basic-lem}(2)
  $\calj$ is a free left $\mathbb C[\h]^W$-module, and hence  so is
 each summand $\ogr^m\calj$.
 Thus we may pick an $\EE$-homogeneous free basis $\{a_{jm}\}$
 for $\ogr^m\calj$.
 Now $a_{jm}=\Theta(a_{jm})\in \ogr^m \caln={\caln}^m/\caln^{m-1}$ and  the surjection
 $\pi_m: \caln^m\to  {\caln}^m/\caln^{m-1}$ is an $\EE$-graded
 surjection. Thus,   for each ${j,m}$ we can pick an
 $\EE$-homogeneous preimage $\theta(a_{jm})\in \caln^m$
 of $\Theta(a_{jm})$.

 Define $\theta$ to be  the
 $\mathbb C[\h]^W$-module map induced by the map
$ a_{jm}\mapsto \theta(a_{jm})$ on
 basis elements. Since $\pi_m$ is a left $\C[\h]^W$-module map,
 a straightforward induction on orders of elements
 ensures that the $\theta(a_{jm})\in \caln^m$ are a
 free basis for the module they generate. The other conclusions of the lemma
 follow automatically from the construction of $\theta$.
  \end{proof}

\subsection{}\label{step-1} As happens with many questions about
$\WW$-invariants, it is easy to prove that  $\Theta$ is surjective on $\hr$.
Given a left $\C[\h]^\WW$-module $M$, we will write  $M[\delta^{-2}]$ for the
localisation $\C[\h]^\WW[\delta^{-2}]\otimes_{\C[\h]^\WW}M$.  Clearly, when $M$
is a left $\C[\h]$-module, $M[\delta^{-2}]$  is naturally isomorphic to
$\C[\h][\delta^{-1}]\otimes_{\C[\h]}M$.

\begin{lem} (1) The inclusion $ \Theta[\delta^{-2}] : \calj[\theta^{-2}]
\hookrightarrow (\ogr \caln)[\delta^{-2}]$
  is an equality.

 (2)  The induced map $ \theta[\delta^{-2}] : \calj[\theta^{-2}]
\to  \caln[\delta^{-2}]$
   is an
  isomorphism.
  This map is graded under   the $\EE$-grading and  is a filtered isomorphism
  under the order filtration, in the sense that
  $\theta[\delta^{-2}]$ maps $ \ord^n\calj[\delta^{-2}] $
   isomorphically to  $\ord^n\caln[\delta^{-2}]$
  for each $n$.
\end{lem}
\begin{proof}  (1) By  \eqref{locdunk}
 $B_{k,k-1}[\delta^{-2}] = eH_{c+k}\delta[\delta^{-2}]e= e (D(\hr)\ast \WW )e,$
 for any $k\in\C$.
 Repeated application of this shows that
 $B_{ij}[\delta^{-2}]  = e (D(\hr)\ast \WW)e$
 and hence, by Corollary~\ref{morrat-cor},
  that $\caln[\delta^{-2}] =  e (D(\hr)\ast \WW )eH_c=e(D(\hr)\ast \WW) .$
 Since $\ord(\delta^2)=0$, we
 deduce that $(\ogr \caln)[\delta^{-2}] =
  e(\C[\hr \oplus \h^{\ast}]\ast \WW).$
On the other hand,
since $\delta^{2k}\in J^k\delta^k\subseteq \C[\h\oplus\h^*]$,
certainly
$\calj[\delta^{-2}] = e\C[\hr \oplus \h^{\ast}]=
 e( \C[\hr \oplus \h^{\ast}]\ast \WW)$.
Since $\Theta$ is given by  inclusion,
$ \Theta[\delta^{-2}]$ is therefore an isomorphism.

(2) By Lemma~\ref{filter-injA}, $\theta $ and hence $ \theta[\delta^{-2}] $
are graded maps under the $\EE$-gradation and filtered under the order
filtration.  Since $\gr( \theta[\delta^{-2}] )= \Theta[\delta^{-2}] $ is an
isomorphism,  necessarily $\theta [\delta^{-2}] $ is a filtered isomorphism.
\end{proof}

\subsection{Notation}\label{eqpoi-sect}  As in \eqref{filter-injA}, set
$\calj=eJ^k\delta^k$,  $\caln=N(k)$ and write
$\theta(\calj)^m=\ord^m \theta(\calj)=\theta(\calj)\cap \caln^m$
for all $m\geq 0$.
We rewrite the $\C[\h]^\WW$-basis of $\theta(\calj)$ constructed in the proof of
Lemma~\ref{filter-injA} as $\{a_{g\ell m}\}$,  where each
$a_{g\ell m}$ is $g$-homogeneous under the $\EE$-gradation
 and has order exactly $\ell$.
Since these were induced from the bases $\{ a_{c\ell }\}$ of $\ogr^\ell \calj$,
the set $\{a_{g\ell m} : \ell\leq t\}$ does give a basis of $\theta(\calj)^t$.

 By  Lemma~\ref{Bbar-freeA}(2), $\caln$ is a free left $\C[\h]^\WW$-module
 and it is certainly graded. Thus, by Theorem~\ref{graded-proj-thm}, it is
  graded-free.
  We may therefore pick a $
\C[\h]^\WW$-basis $\{ b_{g u}\}$ of $\caln$ where, again, each
$b_{g u} $ is $\EE$-homogeneous of degree $g$ but of unspecified  order.
 This basis is far
from unique and  one
cannot expect that  $\{b_{gu} : b_{gu} \in \caln^m\}$
 forms a basis of $\caln^m$;
indeed at this stage we do not even know that  $\caln^m$ is a free $\mathbb
C[\h]^W$-module.

\subsection{}\label{subsec-step42}
We are now ready to put these observations
together to prove the hard part of Proposition~\ref{pre-cohh}.

\begin{prop} \label{grsameA}  Fix $k\geq 0$ and set
$\calj = eJ^k\delta^k$ and $\caln=N(k)$.  Then the  map $\theta:\calj\to\caln$
is an isomorphism. \end{prop}

\begin{proof} Set $\mathfrak{m}= \C[\h]^{W}_+$
and note that $\caln/\mathfrak{m}\caln=\overline{N(k)}$. On the other hand,
in the notation
of Corollary~\ref{gr},
   $\calj/\mathfrak{m}\calj \cong \overline{J^{k}}[K]$
   is the shift of $\overline{J^k}$ by
   $\deg \delta^{k} = K=kn(n-1)/2$.
    By Corollaries~\ref{gr} and  \ref{poincare-S2A},
  we therefore have an equality
of Poincar\'{e} series under the $\EE$-gradation:
 \begin{equation} \label{eqpoiA}
p( \calj/\mathfrak{m}\calj, v) =
v^{K}\frac{\sum_{\mu}
f_{\mu}(1)f_{\mu}(v^{-1})v^{-k(n(\mu) - n(\mu^t))}[n]_v!}
{\prod_{i=2}^n (1-v^{-i})} = p( {\caln/\mathfrak{m}\caln}, v).
\end{equation}

Keep the $\C[\h]^\WW$-bases of $\theta(\calj)\cong \calj$ and $\caln$ described in
Notation~\ref{eqpoi-sect}.
We write ${a(g\ell m)} =
g$ whenever  $a_{g\ell m}$ exists for
that choice of $g,\ell,m$;
thus $\sum_{g\ell m} v^{a(g\ell m)}$ denotes the sum
$\sum v^g$, where one has one copy of $v^g$ for each
 $\ell,m$ for which $a_{g\ell m}$ exists.
Define ${b(gu)}$ analogously.
Since the bases $\{a_{g\ell m}\}$ and $\{b_{gu}\}$ induce
$\C$-bases of $\calj/\mathfrak{m}\calj$, respectively $\overline{N(k)}$,
\eqref{eqpoiA}  can be reinterpreted as
 \begin{eqnarray} \label{eqpoi2}
 \sum_{g,\ell,m} v^{a(g\ell m)}  &= &
 v^{K} \frac{\sum_{\mu} f_{\mu}(v^{-1})f_{\mu}(v)v^{-k(n(\mu)-
 n(\mu^t))}[n]_v!}
 {\prod_{i=2}^n (1-v^{-i})}  =
  \sum_ {g,u}  v^{b(gu)}.\end{eqnarray}

We note that \eqref{eqpoi2} has several consequences for
the $a(g\ell m)$ and $b(gu)$.
\begin{enumerate}
 \item[($\dagger$1)]  For fixed $g$, there exist only finitely many elements
 $a_{g\ell m}$ and $b_{gu}$.
 This is because the
 middle expression in \eqref{eqpoi2} is a well-defined series.

 \item[($\dagger$2)] There exists a universal upper bound
$a(g\ell m)\leq T$.  This is because the numerator in the middle expression in
 (\ref{eqpoi2}) is a finite sum of polynomials.
  However, there is no universal lower bound.

\item[($\dagger$3)] For any   $g_0$, the number of
  $a_{g\ell m}$ with $g=g_0$ equals the number of
 $b_{gu}$ with $g=g_0$.
  This is simply because  $\sum v^{a(g\ell m)} = \sum v^{b(gu)}$
  and the numbers are finite by ($\dagger$1).
\end{enumerate}

We aim to adjust the basis
$\{b_{gu}\}$
to be equal to the basis $\{a_{g\ell m}\}$, and we achieve this by a downwards
induction on $g$.  The induction starts since, by ($\dagger$3),
 there are no basis elements
$b_{gu}$ with $g>T$.

Let $-\infty < G\leq T$ and, by induction,  suppose that
 $\{b_{gu} : u\in \Z\}= \{a_{g\ell m} : \ell,m \in \Z\}$
 for all $g>G$.
Suppose that there exists a basis element
$b_{Gw} \not\in \{a_{G\ell m}\}$.
By Lemma~\ref{step-1}(2), $\theta(\calj)[\delta^{-2}]=\caln[\delta^{-2}]$ and so
there exists a homogeneous element $\mathbf{x}^m\in \C[\h]^\WW$ of
$\EE$-degree $m$ such that
$\mathbf{x}^mb_{Gw}\in \theta(\calj)$. Thus we have
the $\EE$-homogeneous equation
\begin{equation}\label{bigger-kk}
\mathbf{x}^mb_{Gw} = \sum_{g<G} c_{gfh}a_{gfh}
+\sum  c_{Gfh}a_{Gfh}+
\sum_{g>G} c'_{gz}b_{gz},\end{equation}
where $c_{gfh}, c'_{gz}\in \mathbb C[\h]^W$ and summation over
 $f,h,z$ is suppressed.
Since $\theta(\calj)\subseteq \caln$, we may write
each $a_{gfh}$ as an $\EE$-homogeneous sum
$a_{gfh} =\sum d_{\bullet}b_{uz}$ for some $d_\bullet = d_{fghuz} \in \C[\h]^W$
and obtain
\begin{equation}\label{bigger-k}
\mathbf{x}^mb_{Gw} = \sum_{g<G} c_{gfh}d_{\bullet}b_{uz}
+\sum  c_{Gfh}d_{\bullet}b_{uz} +
\sum_{g>G} c'_{gz}b_{gz}.
\end{equation} Both the last two displayed equations  are
 $\EE$-homogeneous of $\EE$-degree $G+m$
and so, by \eqref{bigger-kk},
  each element $c_{gfh}$ must have $\EE$-degree $\geq m$.
Thus the $b_{uz} $ appearing in the first two terms on the right hand side of
 \eqref{bigger-k} must have $\EE$-degree $\leq G$. Thus the only
appearance of   $b_{gz}$ with $g>G$ is in the third sum.
Since the $b_{uz}$ are a $\C[\h]^\WW$-basis of $\caln$,
 that third term   $\sum_{g>G} c'_{gz}b_{gz}$ is actually zero.

Now consider where the specific term
$b_{Gw}$ appears on the right hand side of
\eqref{bigger-k}. For $g<G$, \eqref{bigger-kk} implies that
  $\Edeg c_{gfh}>m$ for each $f, h$  and so
$b_{Gw}$ cannot appear in the first sum.
Thus  it must
appear nontrivially in some term $c_{Gf'h'}d'b_{Gw}$
in the second sum.
In this case, \eqref{bigger-kk} implies that
$\Edeg c_{Gf'h'} =m$. Hence
$d'\in\mathbb C\smallsetminus\{0\}$ and
$$a_{Gf'h'} = d'b_{Gw} +\sum_{(uz)\not= (Gw)}
d_{uz}'' b_{uz}.$$
Thus we can replace $b_{Gw}$ by $a_{Gf'h'}$ in our basis
for $\caln$. By ($\dagger$3), the sets
$\{a_{G\ell m} : \ell,m\in \Z\}$ and $\{b_{Gu} : u\in \Z\}$
have equal finite cardinality.
After a finite number of steps we therefore have $\{b_{Gu}\} \subseteq
\{a_{G\ell m}\}$ and hence $\{b_{Gu}\}=
\{a_{G\ell m}\}$.
This completes the inductive step and hence the proof of the lemma.
\end{proof}

We can now pull everything together and prove both Theorem~\ref{main} and
Proposition~\ref{pre-cohh}.
\subsection{Proof of Proposition~\ref{pre-cohh}} \label{subsec-6.21A}
 Recall from Lemma~\ref{thetainjA} that
 $\Theta : eJ^k\delta^{k}\to \ogr
N(k)$ is the natural inclusion.
On the other hand, for any $k \geq 0$, Proposition~\ref{grsameA}
implies that the map $\theta: eJ^k\delta^{k}\to
N(k)$ is an isomorphism.  Lemma~\ref{filter-injA}(2) therefore implies that
$\gr_{\Lambda} N(k) = \ogr \theta(eJ^k\delta^{k})  =
 \Theta(eJ^k\delta^k) = eJ^k\delta^k$.  \qed

\subsection{Proof of Theorem~\ref{main}}\label{subsec-6.21}
(1) This is immediate from Corollary~\ref{morrat-cor}(1)
 and  Lemma~\ref{Zalgequiv}.

 (2)  Fix $i\geq j\geq 0$. Since $c+j$ still satisfies
 Hypothesis~\ref{main-hyp}, Proposition~\ref{pre-cohh}  implies that
 $\ogr B_{ij}eH_{c+j} =eJ^{i-j}\delta^{i-j}$.
   Multiplying  this identity  on the right by $e$ and applying
  Lemma~\ref{grade-elements} and Corollary~\ref{morrat-cor}(1) gives
$$eJ^{i-j}\delta^{i-j} e    = \ogr(B _{ij}eH_{c+j})e
      =\ogr(B_{ij} eH_{c+j}e) =\ogr B_{ij}   .$$
Since $\delta$ transforms under $\WW$ by the sign representation,
Lemma \ref{corpar}(1) shows that
$eJ^{i-j}\delta^{i-j} e= eA^{i-j}\delta^{i-j}e$.
 Combining these observations gives
 $\ogr B_{ij} = eA^{i-j}\delta^{i-j} e$.
 Therefore, $\ogr B= \bigoplus \ogr B_{ij} = e\widehat{A}e\cong \widehat{A}$,
 as graded vector spaces. In order to ensure that this is an isomorphism
 of graded $\Z$-algebras  we need to check that the multiplication
  in $\ogr B$ coming from
 the tensor product multiplication in $B$ is the same as the natural
  multiplication in $\widehat{A}$.
  This  follows from
  Lemma~\ref{abstract-products}(1).

 (3) The equivalences $\ogr(B)\lqgr \simeq A\lqgr\simeq
  \coh(\hi)$ follow from (2) combined with
 \eqref{zalgex1}, respectively Corollary~\ref{hi-basic-lem2}(1).\qed

\subsection{Corollary}\label{order-free} {\it Assume that $c\in \C$ satisfies
Hypothesis~\ref{main-hyp} and pick $i\geq j\geq 0$. Then, for $m\geq 0$,
 each of the modules $\ord^mN(i)$, $\ogr^mN(i)$, $\ord^m B_{ij}$ and
 $\ogr^m B_{ij}$ is  free as a  left $\C[\h]^\WW$-module.}

 \begin{proof} By construction and Proposition~\ref{pre-cohh},
 the map $\Theta: \ogr N(i)\to  eJ^{i}\delta^i$ is an isomorphism of
 $\ord$-graded modules. Thus $\ogr^mN(i)\cong \ogr^m  eJ^{i}\delta^i$
 is a free $\C[\h]^\WW$-module by Lemma~\ref{hi-basic-lem}.
 By induction on $m$, it follows that $\ord^mN(i)$ is also free.
 The analogous results for $B_{ij}$ follow by multiplying everything
 on the right by~$e$.
 \end{proof}

\subsection{}\label{cohh-subsect} We end the section by noting that
Proposition~\ref{pre-cohh} provides an
interesting connection between $H_c$-modules and the isospectral scheme
$X_n$ defined in \eqref{hi-defn-sec2}.
 Adjusting   to the conventions of this section,
  we identify $\hi = \prj \widetilde{A}$, for
 $\widetilde{A}=\bigoplus A^k\delta^k$. By construction,
  the Procesi
bundle $\PP=\rho_*\OO_{X_n}$  from \eqref{hi-defn-sec2}
is then just the image in $\coh \, \hi$
   of the $\widetilde{A}$-module  $\bigoplus J^k\delta^k$.
 Thus the next  result is an immediate consequence
of  Proposition~\ref{pre-cohh}.

\begin{cor} Assume that $c\in \C$ satisfies Hypothesis~\ref{main-hyp}.
Let
 $e\widetilde{H}_c= \bigoplus_{k\geq 0}B_{k0}\otimes_{U_c} e\widetilde{H}_c$
 be the $B$-module associated to the $U_c$-module $eH_c$
and filter each $B_{k0}\otimes_{U_c}eH_c \cong B_{k0}eH_c$ by the $\ord$
filtration.  Set $\ogr e\widetilde{H}_c = \bigoplus \ogr
B_{k0}eH_c$. Then the sheaf   associated to $\ogr e\widetilde{H}_c$
in   $\coh\hi $ is the Procesi bundle $\PP$.
\qed\end{cor}

\subsection{}\label{cohh-subsect-chat}
Just as Theorem~\ref{main} can be
interpreted as saying that $U_c$ provides a noncommutative model for $\hi$, so
Corollary~\ref{cohh-subsect} can be interpreted as saying  that  the algebra
$H_c$ provides a  noncommutative  model for $X_n$. Here is one aspect of this
analogy. It follows from \cite{BKR} and \cite{hai1} that there is an
equivalence $\xi$ of derived categories  between $\h\oplus\h^*/\WW$ and
$\hi$ that is induced by a Fourier-Mukai transform over $\PP$.  Now pass to the
noncommutative situation, replacing $\h\oplus\h^*/ \WW$, $\hi$ and $\PP$
by  $H_c\lmod$,   $B\lqgr$ and $eH_c$, respectively. Then
Corollary~\ref{morrat-cor}  shows that $eH_c$ still  induces a derived
equivalence between the two categories. Indeed, it is even a  equivalence of
categories. The fact that derived equivalences in the commutative case  can
become full equivalences in the noncommutative case happens elsewhere  and is
in accord with the philosophy behind  \cite[Conjecture~1.6]{GK} (see
\cite[Remark~1.7]{GK}).

  As will be justified in \cite{GS2}, Corollary~\ref{cohh-subsect}  therefore
``sees'' the equivalence $\xi$ and this provides some intriguing connections
between sheaves on $\hi$ and modules over $H_c$.

\subsection{} If one considers Cherednik algebras in characteristic $p>0$,
 where $H_c$ is a finite module over its centre, then the relationship between
 $H_c$ and $\hi$ becomes closer still. For example, \cite{BFG} shows that
  there is even a derived equivalence between $H_c$ and an Azumaya algebra over
 a Frobenius twist of $\hi$. Similarly in characteristic zero, symplectic
 reflection algebras with parameter $t=0$ are finite modules over their
 centre, and \cite[Theorem 1.2]{GSm} shows that there are often derived
 equivalences between these algebras and varieties that deform
 Hilbert schemes.


\section{Tensor product filtrations}\label{sect7}

\subsection{}\label{sect701}
The tensor product decomposition \eqref{tpdef} of the $B_{ij}$  can be used to
give a second filtration on that module by inducing a filtration on $B_{ij}$
from the $\ord$ filtration on the tensorands.  It turns out that the main
theorem is essentially equivalence to the assertion that  the two filtrations
are equal. In this short section we give the details  behind this assertion.
Analogues of this result also hold for the module $N(k)$ defined  in
\eqref{app-to-main} and the module $M(k)=H_{c+k}eB_{k0} =H_{c+k}\delta e
B_{k-1,0}$ defined in \eqref{app-c-1} and so we begin by giving a general
context for all three results.

\subsection{}\label{tens-defn-sect}
 For fixed $i\geq j\geq 0$ we are interested in the following
 tensor product decompositions
\begin{equation}\label{tensor-1}
B_{ij}\cong Q_{c+i-1}^{c+i}\otimes Q_{c+i-2}^{c+i-1}\otimes\cdots \otimes
Q_{c+j}^{c+j+1},
\end{equation}
\begin{equation}\label{tensor-101}
N(i)\cong  Q_{c+i-1}^{c+i}\otimes \cdots \otimes
Q_{c}^{c+1}\otimes eH_c\qquad\mathrm{or}\qquad
N(i) \cong B_{i0}\otimes eH_c
\end{equation}
and
\begin{equation}\label{app-c-cor}
 M(i)\cong H_{c+i}\delta e\otimes_{U_{c+i-1}} B_{i-1,i-2}\otimes
\cdots\otimes_{U_{c+1}} B_{10}
\qquad\mathrm{or}\qquad
M(i) \cong H_{c+i}\delta e\otimes_{U_{c+i-1}}B_{i-1,0}
\end{equation}
where the tensor products are over the appropriate rings $U_k$.
Corresponding to these decompositions we have the
 {\it tensor product filtration} $\ten$ defined by
 the following convention: Given a module $C=C_1\otimes\cdots \otimes C_r$,
 where each $C_j$ is filtered by the $\ord$ filtration, define
 \begin{equation}\label{tensor-2}
\ten^n(C) = \Big\{\sum c_{1}\otimes \cdots\otimes  c_r,
\ \mathrm{where}\ c_m\in \ord^{\ell(m)}(C_m) \ \mathrm{with}\
\sum_{m=1}^{r} \ell(m)\leq n\Big\}.
\end{equation}
As usual, we will write the associated graded module as $\tgr C = \bigoplus
\ten^n C/\ten^{n-1} C$.

\begin{lem}\label{ord-tens}  Assume that $c\in \C$ satisfies
Hypothesis~\ref{main-hyp}. Let $C$ denote one of the objects
$B_{ij}$, $N(i)$ or $M(i)$
and consider the tensor product filtrations induced from one of
the tensor product decompositions (\ref{tensor-1}--\ref{app-c-cor}).
Then $\ord^mC=\ten^mC$, for all $m\geq 0$.
\end{lem}

\begin{proof} We will  prove the result for the decomposition
\eqref{tensor-1} and the first decomposition in each of
\eqref{tensor-101} and \eqref{app-c-cor}. The proof in the remaining cases
is left to the reader as it uses essentially the same argument, although one
needs to use the conclusion of the lemma  for \eqref{tensor-1}.

 In each of the three cases we are given a decomposition
$C=C_1\otimes\cdots \otimes C_r$, say with $\ogr C_j=D_j$
and $\ogr C=D$. Moreover, by Theorem~\ref{main}, respectively
Proposition~\ref{pre-cohh} combined with Lemma~\ref{thetainjA}, respectively
  Proposition~\ref{app-c-prop} combined with Lemma~\ref{thetainjC},
there is an equality
$D_1\cdots D_r=D$ given by  multiplication
in $D(\hr)\ast\WW$.
Equivalently,  the natural multiplication map
$\chi: D_1\otimes \cdots \otimes D_r\to D$ is surjective.
Consider the graded map $\chi$ in more detail.
Given elements $\bar{\alpha}_j\in \ogr^{m(j)} D_j$,
with $m = \sum m(j)$, lift the $\bar{\alpha}_j$ to elements
$\alpha_j\in ord^{m(j)}C_j$. Then $\chi$ is defined
by
$$\chi(\bar{\alpha}_1\otimes\cdots\otimes \bar{\alpha}_r)
=\left(\alpha_1\cdots\alpha_r + \ord^{m - 1}C\right)/\ord^{m - 1}C.
$$
By the definition of
the $\ten$ filtration, this says that
 image of $\chi$ is contained in (and indeed equal to)
$\bigoplus_{m} \bigl(\ten^m C + \ord^{m-1} C\bigr)/\ord^{m-1} C$.
But $\chi$ is surjective. By induction on
$m$ we therefore have $\ord^mC= \ten^m C + \ord^{m-1} C = \ten^m C$.
\end{proof}

\subsection{}\label{ord-tens-chat}
The equality of filtrations given by Lemma~\ref{ord-tens}
is not merely a formality; indeed the result for $B_{ij}$ is essentially the
 same result as Theorem~\ref{main}. To see this, suppose that
 $\ogr B_{ij}=\tgr  B_{ij}$ for all $i\geq j\geq0$.
  As Lemma~\ref{thetainjA}(2) shows,
 $\ogr B_{\ell+1,\ell} = A^1\delta$ for each $\ell$ and
 so, by Lemma~\ref{abstract-products}(2), we get a surjection $\chi$ from
 $E=(A^1\delta)^{\otimes(i-j)}$ onto $\tgr B_{ij}=\ogr B_{ij}$.

 The multiplication map $\phi: E\to (A^1\delta)^{i-j}$ is surjective
 and its kernel is the largest torsion $A^0$-submodule of $(A^1\delta)^{i-j}$.
 On the other hand $\ogr B_{ij}\subseteq e\C[\h\oplus\h^*]^\WW$
 is a torsion-free $A^0$-module and so
 $\mathrm{ker}(\phi) \subseteq \mathrm{ker}(\chi)$.
 Thus $\ogr B_{ij} = E/\mathrm{ker}(\chi) $ is a homomorphic image of
   $ (A^1\delta)^{i-j} $. Since  $ (A^1\delta)^{i-j} $
 is a right ideal of the domain $A^0$, any proper factor
 of  $ (A^1\delta)^{i-j} $  will be torsion. Thus  $\mathrm{ker}(\phi)=
  \mathrm{ker}(\chi)$ and  $\ogr B_{ij} \cong (A^1\delta)^{i-j}$.

\subsection{}\label{order-counter} The observation in \eqref{ord-tens-chat}
suggests that Lemma~\ref{ord-tens} will only hold for very special
decompositions and this is indeed the case. In essence,
Theorem~\ref{main} says that the identity
$B_{ij}\cong B_{i,i-1}\otimes\cdots\otimes B_{j+1,j}$ is a filtered isomorphism.
On the other hand, an identity like $H_c\cong H_ce\otimes_{U_c} eH_c$
from Theorem~\ref{morrat} is clearly not filtered;
in writing the element $1$ as an element of $H_ce \otimes eH_c$
an easy computation shows that
one needs to use commutators of elements from $\C[\h]$ and $\C[\h^*]$
and so $1\notin \ten^0(H_c)$. However,
$ge=ge\cdot 1 \in \ten^0(H_c)$ for any  $0\not=g\in \C[\h]^\WW$ and
so  $\sigma(ge)\sigma(1)=0$  in $\tgr H_c$.
On the other hand, as $1$ is a regular element of
$\ogr H_c$, no such equation is possible $\ogr H_c$.
Thus $\ten H_c\not\cong \ogr H_c$.

As a second example, it is easy to check
 that Lemma~\ref{ord-tens} will  fail for
 $M(i)$ if one introduces one more tensor product,
$M(i) \cong H_{c+i}e \otimes_{U_{c+i}} B_{i0}$.
Indeed, Lemma~\ref{thetainjC} implies that
$ \ogr M(1)=\C[\h\oplus \h^*]\delta e$. On the other hand,
for the given  decomposition Lemmas~\ref{thetainjA} and \ref{abstract-products}
imply that  $\tgr H_c$ is a homomorphic image of
$ T=\ogr H_{c+1}e\otimes_{U_{c+1}}\ogr Q_c^{c+1}\cong
\C[\h\oplus \h^*]e\otimes_{A^0} A^1\delta e.$
Clearly the image of $T$ in $\ogr M(1)$ is just
 $\C[\h\oplus \h^*]e A^1\delta e=J^1\delta e$.
By the argument of the second paragraph of
\eqref{ord-tens-chat}, this is also the image of $\tgr M(1)$ in $\ogr M(1)$.

\appendix
\section{Graded projective modules}\label{app-a}

\subsection{} The aim of this appendix is to prove the following graded
analogue
of a well-know result of Kaplansky \cite[Theorem~2]{Kap}, for which we do
not
know a reference.

\begin{thm}\label{graded-proj-thm}
Let $A=\bigoplus_{i\geq 0} A_i$ be a connected $\NN$-graded $k$-algebra
(thus $A_0=k$).
Let $P$ be a
right $A$-module that is both graded and projective. Then
$P$ is a \emph{graded-free} $A$-module in the sense that
$P$ has a free basis of homogeneous elements.
\end{thm}

\begin{proof}
Throughout this proof all graded maps are graded maps of degree zero.
We will write the degree of a homogeneous element $x\in P$ as $|x|$.

An observation of Eilenberg \cite[Section~1]{Eil}  shows that $P$ is
graded projective in the sense that there is a graded isomorphism $F\cong
P\oplus Q$, for some $A$-module $Q$ and graded-free $A$-module $F$.
We need a minor variant on this result, so we give the proof.
Take a graded surjection
$\phi: F=\bigoplus f_iA\twoheadrightarrow P$ and an
ungraded splitting $\theta: P\to F$. If $p_i=\phi(f_i)$, then write
$\theta(p_i)=g_i+h_i$, where $g_i$ is the homogeneous component of
$\theta(p_i)$  with   $|g_i|=|p_i|$. Then check that
 the map $p_i\mapsto g_i$ also splits
$\phi$. This proof also shows that, if $P$ is countably generated, then
we can take $F$ to be a countably generated graded-free module.

The heart of the proof of the theorem  is contained in the next two
sublemmas.

\subsection{Sublemma}\label{graded-proj-sublemma1} {\it
Under the hypotheses of the
theorem, $P$ is a graded direct sum of countably generated $A$-modules.}

\begin{proof} The proof of
\cite[Theorem~1]{Kap} also works in the category of graded modules.
\end{proof}

\subsection{Sublemma}\label{graded-proj-sublemma}  {\it
Keep the hypotheses of the theorem and assume that  $P$ is
countably generated. If $x\in
P$ then there exists a graded-free direct summand $G$ of $P$ such that
$x\in G$. }

\begin{proof}  By  the result of Eilenberg  described above,
we may pick    a graded isomorphism $F\cong
P\oplus Q$, for some $A$-module $Q$ and countably generated
graded-free $A$-module $F$.
Select a homogeneous basis $\{u_i : i\in \mathbb N\}$
for $F$ such that there is a graded expression $x=\sum_{i=1}^n u_ia_i$,
with
$a_i\in A$ and $n$ as small as possible.

We first claim that
no $a_j$ can be written as a {\it left}  linear combination of the other
$a_\ell$. Indeed, suppose that $a_n=\sum_{i=1}^{n-1} r_ia_i$, for some
$r_i\in A$. By taking the appropriate component we may assume that each
$r_i$ is
homogeneous with
$|r_i|=|a_n|-|a_i|$. It follows that $|u_nr_i|=|u_i|$ and hence that
 $u_i'=u_i+u_nr_i$
is homogeneous. However
$$\sum_{i=1}^{n-1}  u_i'a_i
=\sum_{i=1}^{n-1} u_ia_i + u_n(\sum_{i=1}^{n-1} r_ia_i)
=x.$$
This contradicts the minimality of $n$ and proves the claim.

Reorder the basis $\{u_\ell\}$
so that $|u_i|\leq |u_{i+1}|$ for   $1\leq i\leq n$ and
write $u_i=p_i+q_i$, for $p_i\in P$, $q_i\in Q$, all of the same degree.
Notice that
$P\ni x=\sum u_ia_i=\sum p_ia_i + \sum q_ia_i$
and so $\sum q_ia_i\in P\cap Q=0$. Hence
\begin{equation}\label{kap-1}
x= \sum_{i=1}^n u_ia_i=\sum_{i=1}^n p_ia_i
\end{equation}
Next write each $p_i$ as a homogeneous sum
$p_i=\sum_{j=1}^nu_jc_{ji} + t_i$, where
$t_i\in \sum_{i>n}u_iA$. Then
$$x = \sum_{i=1}^n u_ia_i=\sum p_ia_i =
\sum_{i,j=1}^n u_jc_{ji}a_i  +\sum_{i=1}^n t_i a_i.$$
Since $\{u_i\}$ is a basis,
\begin{equation}\label{kap2}
a_j=\sum_{i=1}^n c_{ji}a_i\qquad\text{for}\quad 1\leq j\leq n.
\end{equation}

We claim that $c_{ji}=0$ for $i<j$ and that
$|c_{ji}|>0$ whenever  $i>j$ (and $c_{ji}\not=0$).
 Since $|u_i|\leq |u_{i+1}|$,
we have $|a_i|\geq |a_{i+1}|$ for each $i$.
Also $|c_{ji}|=|u_i|-|u_j|$ for all $i,j$
and so $c_{ji}=0$ if $|u_i|<|u_j|$.
Thus both parts of the claim are clear
when $|u_i|\not= |u_{j}|$; equivalently, when $|a_i|\not= |a_{j}|$.
So, suppose that $|a_i|=|a_j|$, for some $i\not =j$ and that
$c_{ji}\not=0$.
Then $c_{ji}\in k^*$ and so  \eqref{kap2} expresses $a_i$ as a
 left linear combination of the other
$a_\ell$. This contradicts the initial minimality assumption on
 $n$ and proves the claim.
Note that $c_{jj}=1$ for all $j$, since otherwise \eqref{kap2} would
express
$a_j$ as a left linear combination of the other
$a_\ell$.

The last paragraph implies that $C=(c_{ji})$ is an upper triangular
matrix,
with  units on the diagonal and so  it is invertible. In particular,
$\{p_1,\dots,p_n\}\cup \{u_{n+\ell} : \ell >0\}$ is  a basis for $F$. Thus
$G=\sum_{i=1}^n p_iA$ is a graded-free direct summand of $F$ contained in
$P$.
Thus $G$ is also a graded-free direct summand of $P$
which, by \eqref{kap-1}, contains $x$. \end{proof}

\subsection{} The proof of the theorem follows from the sublemmas by an
easy induction.
By Sublemma~\ref{graded-proj-sublemma1} we may assume that $P$ is
countably generated,
say  by homogeneous elements $z_i $ for $i\in \mathbb N$.
By induction, suppose that there is a graded decomposition
 $P=Q_1\oplus\cdots \oplus Q_n\oplus R_n$, where each
$Q_i$ is graded-free and $z_i\in Q_1\oplus\cdots \oplus Q_i$, for $1\leq
i\leq n$.
By Sublemma~\ref{graded-proj-sublemma} this does hold when $n=1$.
Write  $z_{n+1} = q+r$ as a homogeneous sum, where $
q\in \sum Q_j$ and $r\in R_n$. Since $R_n$
also satisfies the hypotheses of  Sublemma~\ref{graded-proj-sublemma},
$R_n$ has a
graded-free summand $Q_{n+1}$ containing $r$, completing the inductive
step.
Finally,
$$\widetilde{P}\ = \ \lim_{n\to \infty} \big( Q_1\oplus\cdots \oplus
Q_n\big)
\ \cong \ \bigoplus_{i= 1}^\infty Q_i$$
is a graded-free submodule of  $P$ that contains each $z_i$.
Therefore $P=\widetilde{P}$.
\end{proof}


\section{Another module}
\label{C}

\subsection{}\label{app-c-1}
Fix $c\in \C$ that satisfies Hypothesis~\ref{main-hyp}
and an integer $k\geq 0$.
For applications in \cite{GS2} we will need an analogue of
Proposition~\ref{pre-cohh} for the left  $H_{c+k}$-module
$M(k) = H_{c+k}eB_{k0}\subseteq D(\hr)\ast \WW$.
As before, we filter $M(k)$ by the induced order filtration $\ord$,
so that $\ogr M(k)\subseteq \ogr D(\hr)\ast\WW = \C[\reg{\h}\oplus \h^*]\ast
\WW.$ The aim of this appendix is then to prove:

\begin{prop}\label{app-c-prop}
The left $H_{c+k}$-module $M(k) = H_{c+k}eB_{k0}$ satisfies
$\ogr M(k) = J^{k-1}\delta^ke.$
\end{prop}

\smallskip
Recall that  Proposition~\ref{pre-cohh} showed that the module
$N(k)=B_{k0}\otimes eH_c$ had associated graded ring $eJ^k\delta^k$.
In a sense, Proposition~\ref{app-c-prop} is just a left-right analogue
of that result and so much of the present proof is formally very similar to that
 of Proposition~\ref{pre-cohh}.

 We should first explain  why the two results involve different powers of $J^1$.
 The reason is that one can write $M(k) = H_{c+k}eH_{c+k}\delta eB_{k-1,0}$. By
  Corollary~\ref{morrat-cor} and \eqref{morrat11}
  the left hand end of this expression
collapses to give
$ M(k) = H_{c+k}\delta eB_{k-1,0}.$
In particular, $M(1) = H_{c+1}\delta e$.
A routine  computation using Lemmas~\ref{abstract-products} and
\ref{grade-elements}  then gives:

\subsection{Lemma}
\label{thetainjC}
{\it $\ogr M(1) = \cxy \delta e$ while  $J^{k-1}\delta^{k}e
\subseteq  \ogr M(k)$ for all $k\geq 1$.} \qed

\medskip
It  takes considerably
  more work to show that
$J^{k-1}\delta^{k}e $ actually  equals $   \ogr M(k)$ for $k>1$. The
proofs of the  first few steps in this argument are   very similar to those of
 Lemmas~\ref{Bbar-freeA}, \ref{filter-injA} and \ref{step-1}
 in the proof of Proposition~\ref{pre-cohh}  and so
 we will just indicate how to modify the earlier proofs to work here.

\subsection{}\label{B-freeC}
 Since $M(k)$ is a
$(H_{c+k},\UU_c)$-bimodule, the embeddings $\C[\h]\hookrightarrow
H_{c+k}$ and $\C[\h^*]^\WW\hookrightarrow \UU_c$ make $M(k)$ into a
$(\C[\h],\, \C[\h^*]^\WW)$-bimodule.
Let $\C$ be the trivial module over either $\C[\h]$
or $\C[\h^*]^\WW$ and set $
 \overline{M(k)} = \C\otimes_{\C[\h]} M(k)$ and $
\underline{M(k)} = M(k) \otimes_{\C[\h^*]^\WW} \C.$

\begin{lem}\label{Bbar-freeC}
{\rm (1) } $M(k)$ is
 free  as a  left $\C[\h]$-module
and a right  $\C[\h^*]^\WW$-module.
\begin{enumerate}
\item[(2)] $\underline{M(k)}$ is a  finitely generated, free left
$\C[\h]$-module.

\item[(3)]  Analogously,
$\overline{M(k)}$ is a finitely generated, free
right $\C[\h^*]^\WW$-module.
\end{enumerate}
\end{lem}

\begin{proof}
 (1) By Corollary~\ref{morrat-cor}, $M(k)$ is projective as a left
 $H_{c+k}$-module
 and as a right $\UU_c$-module. By
 \eqref{PBW}, $H_{c+k}$  and hence
 $M(k)$ is   free as a left $\C[\h]$-module.
 Similarly, the argument of
  Lemma~\ref{Bbar-freeA}(2) shows that
 $U_c$  and hence $M(k)$ are free right $\C[\h^*]^\WW$-modules.

(2) This is contained in the proof of Lemma~\ref{Bbar-freeA}(3).

(3) Mimic the proof of Lemma~\ref{Bbar-freeA}(4).
\end{proof}

\subsection{}\label{filter-injC}
Using the conventions from \eqref{poincare-S2A},
each $M(k)$ and $J^{k-1}\delta^{k}e$
is   $\EE$-graded. Since  $\C[\h]_+$  is  $\EE$-graded,
the $\EE$-grading on $M(k)$ descends to one on $\overline{M(k)}$.
Similarly, $J^{k-1}\delta^{k}e$  has the order grading
 $\ogr$  from \eqref{filter-injA}.
Write $\Theta: J^{k-1}\delta^ke\hookrightarrow \ogr M(k)$
for the inclusion from Lemma~\ref{thetainjC}.

\begin{lem} There exists  an injective map $\theta :
J^{k-1}\delta^ke\hookrightarrow M(k)$ of left $\mathbb C[\h]$-modules
such that:
\begin{enumerate}
\item $\theta$  is an $\EE$-graded homomorphism
and is a filtered homomorphism  under the order filtration.

\item The
associated graded map
$\ogr \theta: J^{k-1}\delta^ke \to \ogr M(k)$ induced by $\theta$
is precisely $\ogr \theta = \Theta$.

\item In the notation of \eqref{step-1}, the inclusion
 $ \theta[\delta^{-2}] :
  (J^{k-1}\delta^k e)[\delta^{-2}] \to       M(k)[\delta^{-2}]   $
   is an   isomorphism.
  This map is $\EE$-graded  and  is a filtered isomorphism
  under the order filtration.
\end{enumerate}
\end{lem}

\begin{proof} (1,2) As in the proof of Lemma~\ref{filter-injA}, one
constructs $\theta$  by  lifting  a $\EE$-homogeneous basis of the
free $\C[\h]$-module $\ogr^n(J^{k-1}\delta^k)e$ to a set of
$\EE$-homogeneous elements in $\ord^n M(k)$.

(3)  This is essentially the same as the proof of Lemma~\ref{step-1}.
  \end{proof}

\subsection{}
By Lemma~\ref{diaggrad}, $M(k)$ is graded under the
 adjoint  $\hh$-action and, as both copies of $\C$ are $\hh$-graded modules,
 this grading restricts to one on  $\overline{M(k)} $ and $
\underline{M(k)}$. In each case,  we call this {\it the  $\hh$-grading}.
For the reasons given in \eqref{poincare-S2A}, this does not equal the
$\EE$-grading.
\begin{prop}\label{poincare-SC}  If $\overline{M(k)}$ is graded via
the adjoint $\hh$ action, then
it has  Poincar\'{e} series
  $$p(\overline{M(k)}, v) =
\frac{\sum_{\mu} f_{\mu}(1)f_{\mu}(v^{-1}) v^{-(k-1)(n(\mu)
- n(\mu^t))}}{\prod_{i=2}^n
(1-v^{-i})}.
$$\end{prop}

\begin{proof}  This is similar to the proof of Proposition~\ref{poincare-SA}
except that we use the module $Y=H_ce\otimes_{R}\C$, where $R=e\C[\h^*]^\WW e$,
in place of $X=H_c\otimes_{\C[\h^*]}\C$. As in that proposition,
$Y$ is an object in $\widetilde{\OO}_c$ and so we can write
$ [Y] =
\sum_{\mu} p_{\mu} [\widetilde{\Delta}_{c}(\mu)]$ for some
$p_{\mu} \in \Z
[v,v^{-1}]$. To calculate the $p_{\mu}$  note that, by
\eqref{PBW},  $Y\cong
\C[\h]\otimes \C[\h^*]^{\text{co}\WW}$.
Applying $(\C\otimes_{\C[\h]}-)$ to
the equation $ [Y] =\sum p_{\mu} [\widetilde{\Delta}_{c}(\mu)]$
therefore yields $[\C[\h^*]^{\text{co}\WW}]
= \sum_{\mu} p_{\mu} [\mu].$ Thus \eqref{fakedegrees} implies that
$p_{\mu} = f_{\mu}(v^{-1})$ (this is a polynomial in $v^{-1}$ rather than $v$
since $\C[\h^*]$ is negatively $\EE$-graded) and so, as an element of
 $G_0(\widetilde{\OO}_{c})$,
\begin{equation}\label{grot22}
[Y] = \sum_{\mu} f_{\mu}(v^{-1}) [\widetilde{\Delta}_{c}(\mu)] .
\end{equation}

Now consider $\underline{M(k)}$, which we can write as
$H_{c+k}e\otimes_{U_{c+k}} B_{k0}\otimes_{U_c} eY$.
By \eqref{morrat11} and Corollary~\ref{morrat-cor},
 $H_{c+k}e\otimes_{U_{c+k}} e\widetilde{\Delta}_{c+k}(\lambda)
\cong \widetilde{\Delta}_{c+k}(\lambda)$. Thus \eqref{grot22} and
  Lemma~\ref{standAAA} combine to show that
$$
[\underline{M(k)} ] =
 \sum_{\mu} f_{\mu}(v^{-1})v^{k(n(\mu)-n(\mu^t)}
  [\widetilde{\Delta}_{c+k}(\mu)] .$$
 As graded vector spaces,
  $\widetilde{\Delta}_{c+k}(\mu) \cong \C[\h]\otimes \mu$
  and so
  $p(\widetilde{\Delta}_{c+k}(\mu), v) = f_\mu(1) (1-v)^{-(n-1)}$
   by \eqref{fakedegrees2}.
Therefore,
\begin{equation}
\label{wrongsideformulaC}
p(\underline{M(k)}, v) =
\frac{\sum_{\mu} f_\mu(1)f_{\mu}(v^{-1})v^{k(n(\mu)-n(\mu^t)}}
{(1-v)^{(n-1)}}.
\end{equation}

By parts (2) and (3) of  Lemma~\ref{Bbar-freeC}, a
 homogeneous basis for   $\overline{M(k)}$ is given by lifting a homogeneous
  $\C$-basis for $\overline{M(k)}\otimes_{\C[\h^*]^\WW}
  \C = \C \otimes_{\C[\h]} {\underline{M(k)}}.$
 Thus, combining
     \eqref{wrongsideformulaC} with the formul\ae\
 $p(\C[\h^*]^\WW, v) =  \prod_{i=2}^n (1-v^{-i})^{-1}$ and
 $ p(\C[\h], v) = (1-v)^{n-1}$ gives
\begin{equation}\label{wrongsideformula2C}
p(\overline{M(k)}, v) =
\frac{\sum_{\mu} f_{\mu}(v^{-1})f_{\mu}(1) v^{k(n(\mu) -
n(\mu^t))}}{\prod_{i=2}^n (1-v^{-i})}.
\end{equation}
By \cite[Theorem 8]{op} the fake degrees satisfy
 $f_{\mu}(v^{-1}) = f_{\mu^t}(v^{-1})v^{n(\mu^t)-n(\mu)}.$
Combined with \eqref{fakedegrees2} this implies that
 $$f_{\mu}(v^{-1}) f_{\mu}(1) v^{k(n(\mu)-n(\mu^t))}
=
  f_{\mu^t}(v^{-1})f_{\mu^t}(1)
  v^{-(k-1)(n(\mu^t)-n(\mu))}. $$
Substituting this into \eqref{wrongsideformula2C} gives the stated formula for
$p(\overline{M(k)},\, v)$.
\end{proof}

\subsection{}\label{poincare-S2C}
 As was true
 for Corollary~\ref{poincare-S2A}, we need to slightly modify
 Proposition~\ref{poincare-SC} in order to compute the
  Poincar\'e series for $\overline{M(k)}$ under the $\EE$-grading.

\begin{cor} Set $K=kn(n-1)/2$ and $\mathfrak{n}=\C[\h]_+$.
Under the $\EE$-grading there is an equality
 of  Poincar\'{e} series
 \begin{equation}\label{poincare-sss}
 \phantom{\frac{\displaystyle \int}{\displaystyle \int}}
 p(\overline{M(k)}, v) =  v^{K}
 \frac{\sum_{\mu}   f_{\mu}(1)f_{\mu}(v^{-1}) v^{-(k-1)(n(\mu)
- n(\mu^t))}}{\prod_{i=2}^n
(1-v^{-i})}= p(J^{k-1}\delta^k/\mathfrak{n}J^{k-1}\delta^k,\,v).
\end{equation}
\end{cor}

\begin{proof} Equation~\ref{diaggrad1} continues to hold if we
replace $em\delta e$ by $m\delta e$. Thus the argument of
Corollary~\ref{poincare-S2A}(1) combined with   Proposition~\ref{poincare-SC}
and the  formula $M(k)=H_{c+k}\delta e B_{k-1,0}$ gives the first equality of
\eqref{poincare-sss}.

In order to obtain the second equality in \eqref{poincare-sss}, note that
   $p(J^{k-1}\delta^k/\mathfrak{n}J^{k-1}\delta^k,\,v)
= v^Kp(J^{k-1} /\mathfrak{n}J^{k-1},\,v)$. Set
$p(v)=p(J^{k-1} /\mathfrak{n}J^{k-1},\,v)$ and
  $q(v)=p(J^{k-1}/\mathfrak{m}J^{k-1},\,v)$,
  where $\mathfrak{m}=\C[\h]^\WW_+$ The Poincar\'e series $q(v)$
 has been computed in Corollary~\ref{gr}.
 Since that series was obtained by specialising the bigraded Poincar\'e series
  $p(J^d, s, t)$ from  Corollary~\ref{bigr}, it follows immediately
  that
  $$
  p( v)\ = \ \frac{p(\C[\h],v)}{p(\C[\h]^\WW,v)}
 \,  \,q(v) \ = \  \frac{(1-v)^{n-1}}{\prod_{i=2}^{n}(1-v^i)} \, \, q(v)
      \ = \ \frac{q(v)}{[n]_v!}$$
where the final equality  uses  \eqref{factorial-defn}.
 Substituting  these observations into
  Corollary~\ref{gr} gives the second equality in \eqref{poincare-sss}.
 \end{proof}

\subsection{Proof of proposition~\ref{app-c-prop}} \label{subsec-6.21C}
 We first  show that  the map
 $\theta:J^{k-1}\delta^ke \to M(k)$ is an isomorphism for all $k\geq 1$.
This  is analogue of Proposition~\ref{grsameA}.
In that case, a purely formal argument
showed that Proposition~\ref{grsameA} followed from \eqref{eqpoiA}.
The same argument  can be used,
essentially without change, to show that  the bijectivity of $\theta$ follows
from \eqref{poincare-sss}.

 Combined with  Lemma~\ref{filter-injC}(ii) this says that
$\ogr M(k) = \ogr \theta(J^{k-1}\delta^{k}e)  =
J^{k-1}\delta^ke$, as required. \qed


 \clearpage

\section*{Index of Notation}\label{index}
\begin{multicols}{2}
{\small  \baselineskip 14pt

 $\AAA^1$,  $A^1$, alternating polynomials,
 \hfill\eqref{AAA-1-defn}\eqref{A-1-defn}

$\AAA=\bigoplus \AAA^i$, $A=\bigoplus A^i$, \hfill\eqref{AAA-1-defn},\eqref{A-1-defn}

$\widehat{A} = \bigoplus_{i\geq j \geq 0} A^{i-j}$, \hfill\eqref{Aij-defn}

$\BB_1$, the tautological rank $n$  bundle, \hfill\eqref{tauto-defn}

$B=\bigoplus B_{ij}$ for $B_{ij}= \prod_{a=u}^{v-1} Q_{a}^{a+1}$,
\hfill\eqref{B-ring-defn}

canonical grading $W_\alpha $, \hfill\eqref{cangrad-defn}

$d(\mu)= \{ (i,j)\in \NN\times \NN : j <
\mu_{i+1}\}$, \hfill\eqref{d-mu-defn}

$\delta=\prod_{s\in \mathcal{S}} \alpha_s$, \hfill\eqref{delta-defn}

 $\Delta_c(\mu) $, the standard module, \hfill\eqref{standard-defn}

$\widehat{\Delta}_c(\mu) $, the graded standard module,
 \hfill\eqref{graded-standard-defn}

dominance ordering on $\irr\WW$,  \hfill\eqref{dominance-defn}

Dunkl-Cherednik representation $\theta_c$, \hfill\eqref{theta-defn}

$\EE=\sum x_i\delta_i$, the Euler operator, \hfill\eqref{Euler-defn}

$\Edeg$, the Euler grading, \hfill\eqref{Euler-deg-defn}

$e, e_-$, trivial and sign idempotents, \hfill\eqref{e-defn}

fake degrees $f_\mu$, \hfill\eqref{fake-defn}

 $H_c$, the rational Cherednik algebra, $\h,\h^*$,  \hfill\eqref{hc-defn}

$\hr$, \hfill \eqref{delta-defn}

$\mathbf{h} = \mathbf{h}_c=
\frac{1}{2} \sum_{i=1}^{n-1} x_iy_i + y_ix_i \in H_c$, \hfill\eqref{boldh-defn}

$\hdeg$, the $\hh$-grading, \hfill\eqref{hdeg-defn}

Hecke algebra $\hec{q}$, \hfill\eqref{hecke-defn}

Hilbert schemes  $\hin$,\  $\hi$,   \hfill\eqref{hin-defn},\eqref{hi-defn}

$I_\mu$, monomial ideal for a partition $\mu$,  \hfill\eqref{I-mu-defn}

 $\JJJ^1= \C[\C^{2n}]\AAA^1$, $J^1= \C[\h\oplus\h^*]A^1$,
 \hfill\eqref{JJJ-defn}, \eqref{A-1-defn}



$L_c(\mu)$, simple factor of $\Delta_c(\mu)$, \hfill\eqref{L-defn}

$\LL_1 =\OO_{\hin}(1)$, \ $\LL=\OO_{\hi}(1)$,
\hfill\eqref{LLL-defn},\eqref{PP-defn}

 $\mathcal{O}_c$, category $\mathcal{O}$ for $H_c$,
  \hfill\eqref{cat-O-defn}

  $\widetilde{\OO}_c$, graded category $\mathcal{O}$ for $H_c$,
  \hfill\eqref{cat-O-gr-defn}

 $[n]_v! = (1-v)^{-n}\prod_{i=1}^n (1-v^i)$, \hfill \eqref{gr}



$N(k)=B_{k0}eH_c$, \hfill\eqref{app-to-main}

$\overline{N(k)} = \C\otimes  N(k)$, \
$\underline{N(k)} = N(k)\otimes \C$, \hfill\eqref{B-freeA}

 $\ord,\ogr$, order filtration and order gradation,
 \hfill\eqref{order-filt-defn}

$\PP_1$, $\PP$, the  rank $n!$ Procesi bundles,
\hfill\eqref{PP1-defn},\eqref{PP-defn}

$p(M,v)$,   Poincar\'{e} series, \hfill\eqref{gr}

$p(V,s,t)$, bigraded Poincar\'{e} series, \hfill\eqref{Poincare-defn}

$p(M,v,\WW)$,  $\WW$-graded Poincar\'{e} series,  \hfill\eqref{W-poincare}

$\qgr$, $\QGr$, quotient categories,  \hfill\eqref{qgr-defn}

$Q^{c+1}_{c} = eH_{c+1}e_-\delta =e H_{c+1}\delta e$, \hfill\eqref{Q-defn}

$\mathbb{R}(n,l)= {\rm H}^0(\hin, \PP_1\otimes \BB_1^{l})$,
\hfill\eqref{RR-defn}

$\rho_1: \XXX_n\to \hin$, $\rho: X_n\to \hi$,
 \hfill\eqref{isospec-defn},\eqref{hi-defn-sec2}

$\mathbb{S} = \oplus \mathbb{J}^i$, $S=\oplus {J}^i$,
\hfill\eqref{JJJ-defn},\eqref{A-1-defn}

$S_q = S_q(n,n)$, $q$-Schur algebra, \hfill\eqref{schur-defn}

$\mathcal{S}$,  the reflections in $\WW$,
\hfill\eqref{involution-defn}

$\sigma(r)$, the principal symbol of $r$, \hfill\eqref{princ-symbol-defn}

$\sign$, the sign representation of $\WW$, \hfill\eqref{sign-defn}

 Specht module $Sp_q(\mu)$, \hfill\eqref{specht-defn}


 $ \tau: \hin\rightarrow  \C^{2n}/\WW$,   \hfill\eqref{tau-defn}

 $ \tau: \hi\rightarrow \h\oplus\h^*/\WW$,
\hfill\eqref{hi-defn}

$\triv$, the trivial representation of $\WW$, \hfill\eqref{triv-defn}

$U_c=eH_{c}e$, the spherical subalgebra, \hfill\eqref{spherical-defn}

$U_c^-=e_-H_ce_-$, the anti-spherical subalgebra, \hfill\eqref{spherical-defn}


$\WW=\mathfrak{S}_n$, the symmetric group, \hfill\eqref{symmetric-defn}

$\mathbb{X}_n$, $X_n$, isospectral Hilbert schemes,
 \hfill\eqref{isospec-defn}, \eqref{hi-defn}

}
\end{multicols}


\end{document}